%% file: zlst.tex
\documentclass [a4paper, 12pt] {amsart}

\include {lwen}

\begin {document}

\title [Linear selections of superlinear set-valued maps] {Linear selections\\of superlinear set-valued maps\\with some applications to analysis}
\author {D.~V.~Rutsky}
\email {rutsky@pdmi.ras.ru}
\date {\today}
\address {St.Petersburg Department
of Steklov Mathematical Institute RAS
27, Fontanka
191023 St.Petersburg
Russia}

\thanks {This research is supported by the Chebyshev Laboratory
(Department of Mathematics and Mechanics, St. Petersburg State University) under RF Government grant 11.G34.31.0026}


\begin {abstract}
A. Ya. Zaslavskii's results \cite {zaslavsky1981en} on the existence of a linear (affine) selection for a linear (affine)
or superlinear (convex) map $\Phi : K \to 2^Y$ defined on a convex cone (convex set) $K$ having the interpolation property are extended.
We prove that they hold true under more general conditions on the values of the mapping and study some other properties of the selections.
This leads to a characterization of Choquet simplexes in terms of the existence of continuous affine selections for
arbitrary continuous convex maps.  A few applications to analysis are given, including a construction that leads to
the existence of a (not necessarily bounded)
solution for the corona problem in polydisk $\mathbb D^n$ with radial boundary values that are bounded almost everywhere on $\mathbb T^n$.
\end {abstract}

\maketitle

\tableofcontents

\setcounter {section} {-1}

\section {Introduction}

\label {anintroduction}

We begin by a cursory introduction to the subject of this paper before going into formal details.
Suppose that $T : K \to \mathcal V \subset 2^Y$ is a superadditive set-valued map,
i. e. that $T (x) + T (y) \subset T (x + y)$ for all $x$ and $y$ (such maps are also known as convex processes; see, e. g., \cite {rockafellar1970}).
When does $T$ have a linear selection $S$, meaning that $S (x) \in T (x)$ for all $x$, and what is the set of all such selections?
There is a closely related question about affine selections of convex maps, i.~e.
maps satisfying $(1 - \theta) T (x) + \theta T (y) \subset T ((1 - \theta)x + \theta y)$ for all $0 < \theta < 1$ and $x$, $y$.
This question appears to be more intuitive than its counterpart about superlinear maps (see Figure \ref {figconvexmap} below).
Superlinear set-valued maps appear naturally in certain problems of economics among other things
(see, e. g., \cite {rubinov1980}, \cite {aubinfrankowska1990}).

The aim of this paper is to give these questions proper treatment
and show how it leads to some useful tools that
can be applied to certain problems in analysis.
Our starting point is a remarkable work \cite {zaslavsky1981en} by A. Ya. Zaslavskii in which he proved that
a superlinear map $T$ has a linear selection if the following conditions are sufficient: $T$ is positive homogeneous,
$K$ is a cone with the so-called interpolation property (also known as the Riesz decomposition property)
and $\mathcal V$ is a collection of convex sets that has certain compactness properties.
Similar results were obtained independently in \cite {smajdor1990} and \cite {smajdor1996}
without positive homogeneity of $T$ but under the assumption that $K$ is a cone with a cone-basis,
which is in general considerably weaker than the interpolation property.
There are some results for sublinear and concave maps as well as superlinear maps; see, e.~g.,
\cite {protasov2011} and references contained therein.
V.~V.~Gorokhovik obtained in \cite {gorokhovik2008} some interesting results establishing (under some restrictions)
a natural relationship between an affine set-valued map and the set of its linear selections (see Figure~\ref {figconvexmap}).

\begin {figure} [tbph]
\includegraphics [width=13cm] {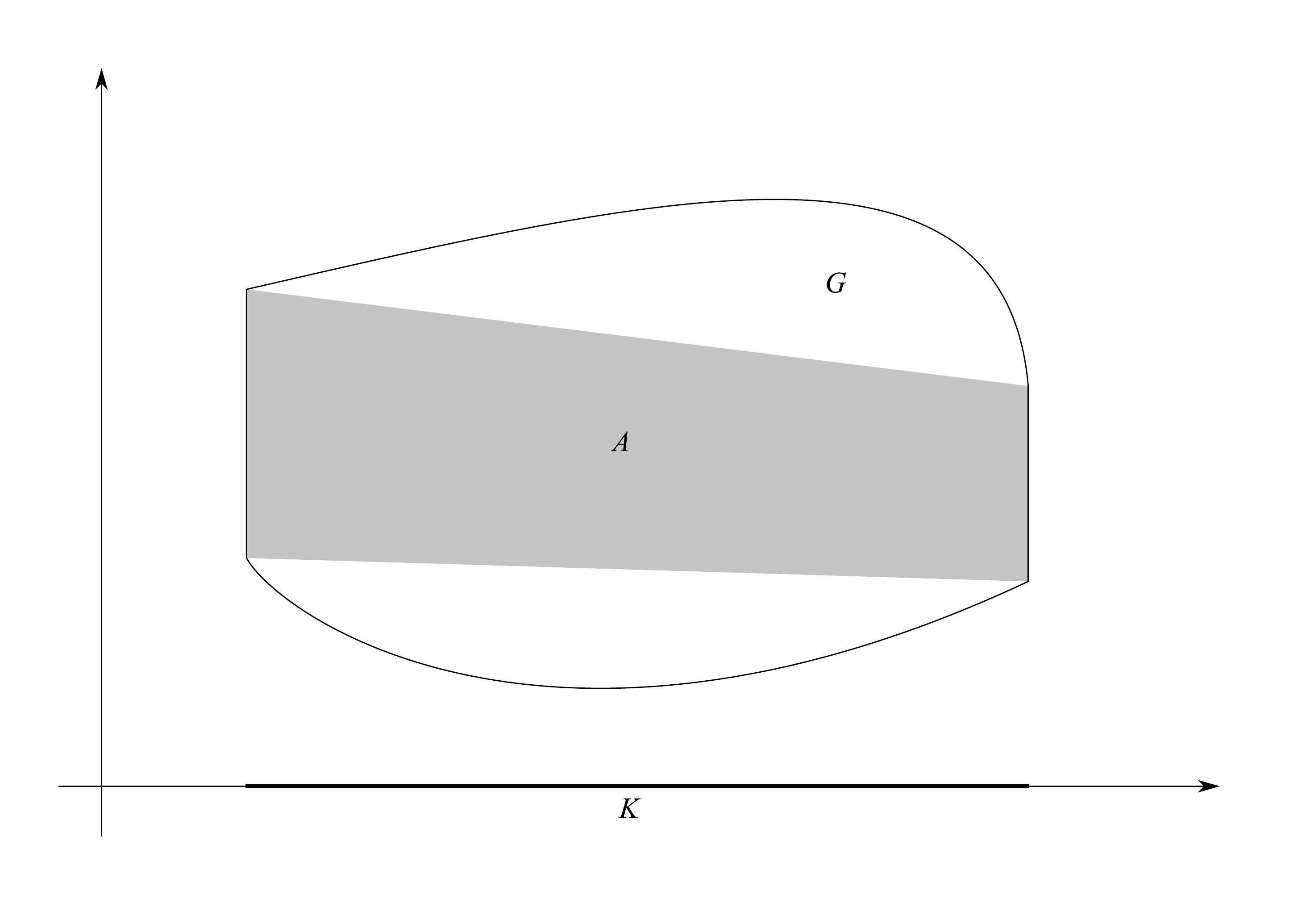}
\caption {An example of the graph $G$ of a convex map defined on a segment $K \subset \mathbb R$.
Gray area $A$ is the graph of the maximal affine submap of this map.  The graph of any affine selection of $G$ is contained within $A$ and
the collection of all such graphs covers the set $A$.}
\label {figconvexmap}
\end {figure}

All of the results related to linear selections of superlinear maps and convex selections of affine maps mentioned up to this point
are obtained under some typical restrictions that allow fairly easy proofs,
leaving the necessity of these restrictions largely unexplored.
Although it appears to be difficult to characterize the exact necessary and sufficient conditions in the general case
for results of this sort,
in this paper we will study this question to an extent for various cases
and obtain some generalizations.
Perhaps the most significant improvement is that we consider maps taking values in \emph {sets of compact type}
(i.~e. in a collection of sets closed under intersection of centered families;
see Section \ref {lsolm} for the definition) instead of the usual compact sets,
which allows us to include maps taking values in convex closed sets in a Banach ideal lattice on $\mathbb Z$
having the Fatou property that are closed in measure.
Moderate refinement of the method of \cite {zaslavsky1981en} leads to
an exact condition on points $x$ and $y$ that guarantees the existence of a linear selection $S$
satisfying $S x = y$ for a given superlinear set-valued map $T$.
In Section \ref {acocs} we further characterize the sets of affine selections of a convex map under some restrictions
and show, as a fairly easy consequence of the results presented in this paper,
that given a compact convex set $K$ in a locally convex linear topological space the following conditions are equivalent.
\begin {enumerate}
\item
Every convex map on $K$ taking values in bounded segments of the line $\mathbb R$ has an affine selection.
\item
$K$ is a Choquet simplex.
\end {enumerate}
The compact type condition cannot be dropped in general, as will be evident in Section \ref {sotctc}.
However, in the case of a cone having a cone-basis (see definition in Section \ref {conebases}),
which in finite-dimensional case is any cone having a simplex as a base, it is not needed.
Although the assumption of the existence of a cone-basis is much stronger than the Riesz decomposition property,
this result has an interesting by-product: in Section \ref {sotctc} we show that for no infinite-dimensional normed
Banach ideal lattices of measurable functions the cone of nonnegative functions has a cone-basis.

The paper is organized as follows.
In Section \ref {introduction} we introduce basic definitions and discuss some simple properties.
In Section \ref {asosas} Zaslavskii's result concerning additive submaps of superadditive maps is presented.
Section \ref {cfac} contains some definitions and facts mostly of topological nature that will be used only sparingly, mostly for topological refinements
of various results.
In Section \ref {lsolm} main results about linear selections of superlinear and linear set-valued maps are presented.
Proofs of the main theorem concerning linear selections of linear maps are given in Section \ref {lslm}.
In Section \ref {ascm} we establish the link between linear selections of superlinear maps and affine selections of convex maps,
which allows us to transfer the results of Section \ref {lsolm} directly to the latter case.
Section \ref {acocs} is devoted to a characterization of Choquet simplexes in terms of the existence of certain linear selections.
In Section~\ref {mosom} we provide
a description of all continuous affine selections of an upper semicontinuous convex map defined on a metrizable Bauer simplex and
acting into closed subsets of a Banach space and show that this method essentially characterizes only linear selections of superlinear maps acting
on the set of Baire measures on a metrizable compact set.
In Section~\ref {scm} we show that the main results of this paper are applicable to maps taking values in
convex sets that are closed in measure in a Banach ideal lattice on $\mathbb Z$ satisfying the Fatou property.
Section \ref {conebases} contains simple versions of the main results for cones having a cone-basis.  Although it builds on the
terminology introduced up to that point, this section is essentially independent of the other sections.
In Section \ref {sotctc} we give an example showing that in a core result about the existence of linear selections for linear
set-valued maps the condition that the values of the map have compact type cannot be dropped
and obtain a characterization of normed Banach ideal lattices of measurable functions such that the cone of positive functions admits a cone-basis.
In Section~\ref {ssvm} we extend main results of this paper to certain cases of maps acting on the cone of positive functions from $\lclassg {1}$.
In Section \ref {addandlin} we discuss briefly the relationship between additivity and linearity.
The rest of the paper is devoted to applications.  In fact, in order to understand them only familiarity with the main results of
\cite {zaslavsky1981en} is required; the necessary results are covered in 
Sections~\ref {introduction}, \ref {asosas} and \ref {lsolm}.
In Section~\ref {ationbo} we show how linear selections of superlinear maps arise in questions related to
the existence of a right inverse for a linear operator taking values in $\lclassg {1}$ and interpret the main results
in this setting.
In Section~\ref {atass} we give a rather superficial treatment to a similar question whether
an estimate for a best approximation implies the existence of a linear operator that realizes it.
Finally, in Section \ref {aafpd} we sharpen a result from the theory of analytic functions in the polydisk
with the help of linear selections of superlinear maps,
which allows us to
obtain a weak version of the Corona Theorem on the polydisk in Section \ref {atcp}.
The last Section \ref {concrem} contains some remarks to various parts of the paper.

\section {Semilinear spaces}

\label {introduction}

Let us now introduce the subjects of this paper in a most general form suitable for stating the main results.
A set $\mathcal S$ equipped with a binary operation of addition $+$ and multiplication by positive real scalars is called
a \emph {semilinear space} (over the ordered field of real numbers $\mathbb R$) \cite [\S 1.3] {kutateladzerubinov1971} if $\mathcal S$
is an Abelian semigroup with respect to addition $+$
and the axioms of a semilinear space
$\lambda (x + y) = \lambda x + \lambda y$, $1 x = x$,
$\lambda (\mu x) = (\lambda \mu) x$, $\lambda x + \mu x = (\lambda + \mu) x$ hold true for all $x, y \in \mathcal S$ and $\lambda, \mu > 0$.
We say that a semilinear space is \emph {monoid} if it contains an identity $\mathbf 0$ with respect to addition.
Any semilinear non-monoid space $\mathcal S$ can be extended to a monoid semilinear space by adding the identity
$\mathbf 0$ and extending the definition of the operations as follows:
$x + \mathbf 0 = x$ and $\lambda \mathbf 0 = \mathbf 0$ for all $x \in \mathcal S$ and $\lambda > 0$.
Every linear space $Y$ is a monoid semilinear space.
Moveover, every convex cone $K$ in a semilinear space $\mathcal S$ is a semilinear space.
However, not every semilinear space can be identified with a cone in a linear space.
For example, the simplest nontrivial
Boolean algebra $\{0, 1\}$ with operations $0 + 0 = 0$, $1 + 0 = 0 + 1 = 1 + 1 = 1$, $\lambda 0 = 0$ and $\lambda 1 = 1$ for all $\lambda > 0$
is a two-point semilinear space, but there exists no two-point cone in a linear space (in this paper
we work with the linear spaces over the field of reals $\mathbb R$ only).
In fact, any ring of sets $\mathcal R$ is a semilinear space with operations defined by $A + B = A \cup B$ and $\lambda A = A$
for all $A, B \in \mathcal R$ and $\lambda > 0$.

For every semilinear space $\mathcal S$ the set $\mathcal C (\mathcal S)$ of nonempty convex sets $A, B \subset S$ equipped with the usual
multiplication $\lambda A = \{\lambda x \mid x \in A\}$ and Minkowski sum $A + B = \{ x + y \mid x \in A, y \in B \}$ is a semilinear space.
It is easy to see that $\mathcal C (\mathcal S)$ is monoid if and only if $\mathcal S$ is monoid.

A relation $\preccurlyeq$ on a semilinear space $\mathcal S$ is said to be compatible with the semilinear structure
of $\mathcal S$ if $\preccurlyeq$ is invariant under
addition and multiplication by positive constants, that is,
$x \preccurlyeq y$ for some $x, y \in \mathcal S$ implies $\lambda x \preccurlyeq \lambda y$ and
$x + z \preccurlyeq y + z$ for all $\lambda > 0$ and $z \in \mathcal S$.
There is a natural transitive relation $\preccurlyeq$
associated with every semilinear space $\mathcal S$: $x \preccurlyeq y$ for $x, y \in \mathcal S$ if and only if
$y - x \in \mathcal S$; that is, $x \preccurlyeq y$ means that there exists $z \in \mathcal S$ such that $y = x + z$.
This relation is compatible with the semilinear structure of $\mathcal S$.  We call $\preccurlyeq$ the \emph {intrinsic order}
of the semilinear space $\mathcal S$.
Relation $\preccurlyeq$ is a partial order if and only if $\mathcal S$ is monoid and has no element invertible in $+$ except $0$;
in this case a relationship $x_1 - y_1 = x_2 - y_2 \Leftrightarrow x_1 + y_2 = x_2 + y_1$ on pairs $x - y$ for $x, y \in \mathcal S$
is an equivalence relationship, and
the set $\mathcal S - \mathcal S = \left\{ x - y \mid x, y \in \mathcal S\right\}$
with linear operations induced by $\mathcal S$ is a vector space having $\mathcal S$ as a generating cone.
Note that $\mathcal C (\mathcal S)$ is usually ordered by inclusion $\subset$, which is a partial order
compatible with the semilinear structure of $\mathcal S$.
It is easy to see that in general there is no relation between the two orders $\preccurlyeq$ and $\subset$.
A semilinear space $\mathcal S$ is embedded canonically into $\mathcal C (\mathcal S)$ by the map $x \mapsto \{x\}$.

Suppose that $\mathcal S$ and $\mathcal V$ are semilinear spaces, $\mathcal V$ is equipped with a partial order $\leqslant$
and $T : K \to \mathcal V$
is a map defined on a convex set $K \subset \mathcal S$.
We say that $T$ is \emph {superadditive} if $K$ is a cone and $T (x) + T (y) \leqslant T (x + y)$ for all $x, y \in K$.
$T$ is \emph {superlinear} if $T$ is superadditive and positive homogeneous (i.~e. $T (\lambda x) = \lambda T x$ for all $x \in K$ and $\lambda > 0$).
$T$ is \emph {additive} if $K$ is a cone and $T (x) + T (y) = T (x + y)$ for all $x, y \in K$.
$T$ is \emph {linear\footnote {\emph {quasilinear} in \cite {zaslavsky1981en}; this term looks somewhat clumsy, although there probably is a certain
point in distinguishing the combination of additivity and positive homogeneity from linearity in some situations.}}
if $T$ is additive and positive homogeneous, 
or, equivalently, if $T$ is affine and positive homogeneous.

A map $S : K \to \mathcal V$ is a \emph {submap} of $T$ if $S \leqslant T$ pointwise; that is, $S (x) \leqslant T (x)$ for all $x \in K$.
Now suppose that $T : K \to \mathcal C (\mathcal V)$ is a set-valued map.
A map $S : K \to Y$ is called a \emph {selection} of $T$ if $S (x) \in T (x)$ for all $x \in K$; that is, $S$ is single-valued and $S \leqslant T$.
It is natural to ask whether a given superlinear map has linear selections,
and if it does, how the set of all such selections can be parametrized and what additional requirements like continuity or passage
through a given point in the graph can be imposed on selections.
Let us establish some relevant notions before proceeding to the answers in the following sections.

Suppose that $\mathcal V$ is a semilinear space with a partial order $\leqslant$ compatible with the semilinear structure of $\mathcal V$.
Note that even if $\mathcal V$ is a monoid we do not require that $0 \leqslant x$ for all $x \in \mathcal V$, and indeed this is not generally
the case for $\mathcal V = \mathcal C (\mathcal S)$ ordered by inclusion which is our main example.
The following two properties are crucial for the main argument of \cite {zaslavsky1981en}.
A set $D$ with a partial order $\leqslant$ is called \emph {lower directed} if every couple of elements $x, y \in D$ has a lower bound $z \in D$,
$z \leqslant x$, $z \leqslant y$.
We say that the order $\leqslant$ is \emph {lower complete} in $\mathcal V$ if every lower directed set $D \subset \mathcal V$ has a greatest lower bound
$\inf D \in \mathcal V$.
We say that a semilinear space $\mathcal V$ is \emph {order regular} with respect to a partial order $\leqslant$ if
it is lower complete with respect to $\leqslant$ and for any lower directed sets $A, B \subset \mathcal V$ and $z \in \mathcal V$ such that
$z \leqslant x + y$ for all $x \in A$ and $y \in B$ it is also true that $z \leqslant \inf A + \inf B$.
A prominent example of an order regular semilinear space is the space of all nonempty
compact convex sets of a locally convex Hausdorff linear topological space
ordered by inclusion; we will treat it in Section \ref {lsolm} below (see Proposition \ref {csgs}).

Let $\mathcal S$ be a semilinear space equipped with its intrinsic transitive relation $\preccurlyeq$.  This relation is extended naturally
onto sets $A, B \subset \mathcal S$ as follows: $A \preccurlyeq B$ if and only if $x \preccurlyeq y$ for all $x \in A$ and $y \in B$.
We use the Minkowski sum $A + B = \{x + y \mid x \in A, y \in B\}$ for $A, B \subset \mathcal S$.
The initial segments $[0, x]$ are defined by
$[0, x] = \{ y \in \mathcal S \mid y \preccurlyeq x \}$ for $x \in \mathcal S$.
We say that $\mathcal S$ has the \emph {interpolation property} if for every finite sets $A, B \subset \mathcal S$
such that $A \preccurlyeq B$ there exists an intermediate point $x \in \mathcal S$ satisfying $A \preccurlyeq x \preccurlyeq B$;
$\mathcal S$ has the \emph {Riesz decomposition property} if $[0, x] + [0, y] = [0, x + y]$
for any $x, y \in \mathcal S$.
A well-known theorem of F.~Riesz shows (see e. g. \cite [\S 1.8, Theorem 1.54] {conesandduality}; although account in the book is given for ordered
vector spaces only, it is easy to verify that the proof remains valid in the semilinear setting)
that these two properties are equivalent to one another and equivalent
to the following convenient purely algebraic property: for any $z \in \mathcal S$ and
$\{x_j\}_{j = 1}^m, \{y_k\}_{k = 1}^n \subset \mathcal S$ such that
$z = \sum_{j = 1}^m x_j = \sum_{k = 1}^n y_k$ there exist some $\{z_{jk}\}_{j = 1, k = 1}^{m, n} \subset \mathcal S$ satisfying
$x_j = \sum_{k = 1}^n z_{jk}$ and $y_k = \sum_{j = 1}^m z_{jk}$ for all $j$ and $k$.
We say that a semilinear space $\mathcal S$ is a \emph {lattice} or a \emph {Riesz space} if
for any $x, y \in \mathcal S$ there exists a least upper bound $x \vee y$ and a greatest lower bound $x \wedge y$ in $\mathcal S$.
It is obvious that every semilinear space which is a lattice
satisfies the interpolation property, but the converse is not true (see, e. g., \cite [\S 1.8, Example~1.58] {conesandduality}).

We now present a simple but crucial pattern that will be used repeatedly in applications in Sections \ref {ationbo}--\ref {aafpd} below,
mostly with $X = \lclassg {1}$ and $\mathcal S = \left[ \lclassg {1} \right]_+ = \{ f \in \lclassg {1} \mid f \geqslant 0 \text { a. e.}\}$.
We say that $\|\cdot\|_X$ is a \emph {quasiseminorm} on a semilinear space $X$ if it is a real-valued positive homogeneous function on $Z$ satisfying
$\|f + g\|_Z \leqslant A (\|f\|_Z + \|g\|_Z)$ for all $f, g \in Z$
with a constant $A$ independent of $f$ and $g$.  Surely the norm of a normed space is a quasiseminorm.
\begin {proposition}[see {\cite [Proposition 8] {zaslavsky1981en}}]
\label {gccont}
Suppose that $X$ is a linear space and $\mathcal S \subset X$ is a generating cone for $X$.
Then every linear operator $T : \mathcal S \to Y$ into a linear space $Y$ can be extended to all $X$ by
$T (x - y) = T x - T y$ for $x, y \in \mathcal S$.
Suppose also that $Y$ is a normed space, $X$ is a quasiseminormed vector lattice,
the quasiseminorm $\|\cdot\|_X$ is additive on the cone $\mathcal S = X_+$ of all nonnegative elements of $X$ and $T$ is bounded.
Then the extension of $T$ described above is bounded with the same norm.
\end {proposition}
Indeed, since $T$ is linear the extension depends only on $x - y$: if $x - y = x_1 - y_1$, then $x + y_1 = x_1 + y$,
$T (x) + T (y_1) = T (x_1) + T (y)$, and therefore $T (x) - T (y) = T (x_1) - T (y_1)$.
Now suppose that the boundedness assumptions of Proposition \ref {gccont} are satisfied.
Since $X$ is a lattice, the cone $X_+$ is generating for $X$, so for arbitrary $z \in X$
we can take $x = z \vee 0$ and $y = z - x$.
Then
\begin {multline*}
\|T z\|_Y = \|T (x) - T (y)\|_Y \leqslant
\\
\|T\|_{X \to Y} (\|z \vee 0\|_X + \|-(z \wedge 0)\|_X) = 
\\
\|T\|_{X \to Y}  \| |z| \|_X = \|T\|_{X \to Y}  \| z \|_X
\end {multline*}
since the norm $\|\cdot\|_X$ is assumed to be additive on $X_+$.
The proof of Proposition \ref {gccont} is complete.

\section {Additive submaps of superadditive maps}

\label {asosas}

The starting point of the present paper is the following result
obtained\footnote {Strictly speaking, in \cite {zaslavsky1981en} it was stated for superlinear maps,
but the proof works for additive maps just as well.}
by A. Ya. Zaslavskii \cite [Theorem 1] {zaslavsky1981en} that we are going to discuss in this section.
\begin {theorem}
\label {linselg}
Suppose that $\mathcal S$ is a semilinear space having the interpolation property, $\mathcal V$ is a semilinear space
equipped with a partial order $\leqslant$ compatible with the semilinear structure and $\mathcal V$ is order regular
with respect to $\leqslant$.  Then any superadditive map $T : \mathcal S \to \mathcal V$ has an additive submap $S \leqslant T$
given by the formula
\begin {equation}
\label {thesubmap}
S (z) = \inf_{z = \sum_j x_j} \sum_j T (x_j), \quad z \in \mathcal S,
\end {equation}
where the infinum is taken over all finite decompositions $x = \sum_j x_j$.
Map $S$ is the greatest additive submap of $T$,
i. e. if $R \leqslant T$ is an additive submap of $T$ then $R \leqslant S$.
If $T$ is positive homogeneous then $S$ is linear.
\end {theorem}
The proof is straightforward.
First, note that for any $z \in \mathcal S$
the interpolation property applied to any two decompositions $z = \sum_j x_j$ and $z = \sum_k y_k$ yields
a decomposition $z = \sum_{jk} z_{jk}$ satisfying $x_j = \sum_k z_{jk}$ and $y_k = \sum_j z_{jk}$ for all $j$ and $k$,
so by subadditivity
$$
\sum_{jk} T (z_{jk}) \leqslant \sum_j T (x_j) \bigwedge \sum_k T (y_k),
$$
which means that
the infinum in \eqref {thesubmap} exists as an infinum of a lower directed set because of the order regularity of $\mathcal V$ with respect to $\leqslant$.

Let us now verify that $S$ is additive.  Suppose that $x, y \in \mathcal S$,
$z = x + y$ and $z = \sum_j z_j$ is an arbitrary finite decomposition in $\mathcal S$.  The interpolation property of $\mathcal S$
in the form of the Riesz decomposition property shows that there exist some $x_j, y_j \in \mathcal S$ such that $z_j = x_j + y_j$ for all $j$.
Therefore
$$
S (x) + S (y) \leqslant \sum_j T (x_j) + \sum_j T (y_j) = \sum_j [T (x_j) + T (y_j)] \leqslant \sum_j T (z_j),
$$
and taking infinum over all such decompositions $z = \sum_j z_j$ shows that
$S (x) + S (y) \leqslant S (z) = S (x + y)$.
On the other hand, for any decompositions $x = \sum_j x_j$ and $y = \sum_k y_k$ into $x_j, y_k \in \mathcal S$ we have
$x + y = \sum_j x_j + \sum_k y_k$, so
$S (x + y) \leqslant \sum_j T (x_j) + \sum_k T (y_k)$.  This yields $S (x + y) \leqslant S (x) + S (y)$ by the order regularity of $\leqslant$
in $\mathcal V$.  Thus we have verified the additivity of $S$.

The fact that $S$ is the greatest linear submap of $T$ is almost immediate.  If $R$ is an additive submap of $T$ then
$$
R (x) = \sum_{j = 1}^N R (x_j) \leqslant \sum_{j = 1}^N T (x_j)
$$
for any $N \in \mathbb N$,
$x = \sum_{j = 1}^N x_j$, $\{x_j\}_{j = 1}^N \subset \mathcal S$, and
it suffices to take the infinum over all such decompositions.  Finally, it is easy to see that
positive homogeneity of $T$ implies positive homogeneity of $S$.  The proof of Theorem \ref {linselg} is complete.

Note that under the assumptions of Theorem \ref {linselg} the additive map $S$ is not necessarily linear if $T$ is merely additive, although
some rather weak additional boundedness assumptions imposed on $T$ imply linearity of~$S$ in this case; see Section \ref {addandlin}.
Since the superadditive maps that arise in problems of analysis are usually positive homogeneous, we will not focus our attention on the relationship
between additivity and linearity in most of the present work.

An interesting question is, to what extent the Riesz decomposition property
and other conditions are actually necessary for the conclusion of Theorem \ref {linselg}?
It was mentioned in \cite {zaslavsky1981en} that there exists an example of cones in $\mathbb R^n$ and maps defined on them
indicating that the Riesz decomposition property
condition on $\mathcal S$ cannot be dropped.
In fact, Theorem~\ref {riesznecc} in Section \ref {acocs} below shows via an application of the Choquet theory and various techniques
that for a large class of semilinear spaces
that covers the typical applications of the theory
the Riesz decomposition property is \emph {necessary} for the conclusion of Theorem \ref {linselg},
yielding at once a wide variety of
such examples and also providing an interesting characterization of Choquet simplexes in terms related to the conclusion of Theorem \ref {linselg}.
On the other hand,
it is also possible to construct a nontrivial example
showing that in the present generality the Riesz decomposition property is not \emph {always} necessary for the
conclusion of Theorem \ref {linselg}.
Let $\mathcal V = \{0, 1\}$ be the two-point semilinear space described in Section \ref {introduction},
and define a cone
\begin {multline*}
\mathcal S = \left\{ (x, y, z) \mid z \geqslant 0, \text {and} \right.
\\
\left.
\text {either $x = 0$, $0 \leqslant y \leqslant z$, or both $0 < x < z$ and $0 < y < z$} \right\}
\end {multline*}
in $\mathbb R^3$ with open square base\footnote {For the definition of the base of a cone see, e. g., in Section \ref {acocs}.
It is slightly easier to visualise this example in the affine reformulation; see Section \ref {ascm}.}
$0 < x < 1$, $0 < y < 1$ plus one edge $x = 0$, $0 \leqslant y \leqslant 1$.
It is easy to see that any superlinear map $T : \mathcal S \to \mathcal V$ is constant on the
topological interior
$$
\tint \mathcal S = \left\{ (x, y, z) \mid z > 0,\, 0 < x < z,\, 0 < y < z \right\},
$$
on the rays
$$
\mathcal S_1 = \{ (x, y, z) \mid z > 0,\, x = y = 0\},
$$
$$
\mathcal S_2 = \{ (x, y, z) \mid z > 0,\, x = 0, y = z\}
$$
and on the open support cone
$\mathcal S_0 = \{ (x, y, z) \mid z > 0,\, x = 0, 0 < y < z\}$; furthermore, there are exactly 6 superlinear maps
$T : \mathcal S \to \mathcal V$ as follows: $T = 0$ everywhere, $T = 1$ on $\tint \mathcal S$ and $T = 0$ elsewhere, and $T = 1$ on
$\mathcal S_0 \cup \tint \mathcal S$ with arbitrary constant values of~$T$ on~$\mathcal S_1$ and~$\mathcal S_2$.
All of these maps are already linear except for the map~$T$ defined by $T = 1$ on $\mathcal S_0 \cup \tint \mathcal S$
and $T = 0$ on $\mathcal S_1 \cup \mathcal S_2$.
Expression \eqref {thesubmap} yields the correct greatest linear submap
$S = 1$ on $\tint \mathcal S$ and $S = 0$ elsewhere, so the conclusion of Theorem \ref {linselg} still holds true for this map $T$.
However, cone $\mathcal S$ does not have the Riesz decomposition property\footnote {This easily follows
from the fact that the closure of the base
of $\mathcal S$ is not a simplex; see, e.~g., \cite {shashkin1973} or \cite {convexityanal} and also Section \ref {acocs} below.}
since, e.~g.,
$\left(\frac 1 2, \frac 1 2, 1\right) \in \left[(0, 0, 0), \left(0, \frac 1 2, 1\right) + \left(\frac 1 2, 0, 1\right)\right]$
but
$$
\left(\frac 1 2, \frac 1 2, 1\right) \notin \left[(0, 0, 0), \left(0, \frac 1 2, 1\right)\right] + \left[(0, 0, 0), \left(\frac 1 2, 0, 1\right)\right].
$$
It is also easy to construct an example where there is a greatest linear submap for every superlinear map, but it is not always given by the
formula~\eqref {thesubmap}.  For example, take $\mathcal V$ as before and define a cone
$\mathcal S = \{ (x, y, z) \mid z > 0, \, x^2 + y^2 \leqslant z \}$ with disk $\{z = 1,\, x^2 + y^2 \leqslant 1\}$ as a base.
In this case linear maps from $\mathcal S$ to $\mathcal V$ are just constants and maps that are $0$ on a single ray intersecting the circle
$\{z = 1, x^2 + y^2 = 1\}$ and $1$ everywhere else, so there always is a greatest submap for any given superlinear map $T : \mathcal S \to \mathcal V$.
However, if $T$ is zero on rays intersecting points $z = x = 1$, $y = \pm 1$ and $1$ everywhere else then the greatest linear submap of $T$ is zero,
but the expression~\eqref {thesubmap} yields a sublinear map $S$ that takes the value $1$ on every point of the circle
$\{z = 1, x^2 + y^2 = 1\}$ except $z = x = 1$, $y = \pm 1$.

As for the space $\mathcal V$, in Proposition \ref {cbrdpnn} of Section \ref {conebases} 
we will see that replacing the Riesz decomposition property assumption by a much stronger
one (as will be evident in Section \ref {sotctc})
of the existence of a cone-basis for $\mathcal S$ makes any requirement of Theorem \ref {linselg} set on the space $\mathcal V$ unnecessary.
On the other hand,
it is fairly easy to construct an example showing that the lower completeness assumption imposed on $\mathcal V$ cannot be dropped in general.
Take, for example, the space $\mathcal S = \left[ \lclass {1} {[0, 1]} \right]_+ \setminus \{0\}$ of nonnegative summable functions on the
unit segment $[0, 1]$ that are not identically $0$ (up to a set of measure $0$), the space $\mathcal V = \mathbb R_+ \setminus \{0\}$ of positive real
numbers and the map $T : \mathcal S \to \mathcal V$ defined by the formula
$T (f) = \left|\{f > 0\}\right|\, \int f$; here $\left|\{f > 0\}\right|$ denotes the Lebesgue measure of the support of $f$.
Positive homogeneity of $T$ is trivial, and superadditivity is also easy to verify: for any
$f, g \in \mathcal S$ we have
\begin {multline*}
T (f + g) = \left|\{f + g > 0\}\right|\, \int f + \left|\{f + g > 0\}\right|\, \int g \geqslant
\\
\left|\{f > 0\}\right|\, \int f + \left|\{g > 0\}\right|\, \int g = T (f) + T (g).
\end {multline*}
However, for any $f \in \mathcal S$ and $N \in \mathbb N$ we can divide the support $\{f > 0\}$ of $f$ into $N$ parts $A_j$ of equal measure
and decompose $f$ into a sum of its parts $f_j = f \chi_{A_j}$ on the sets $A_j$,
so for any linear submap $R$ of $T$ we have
\begin {equation*}
R (f) = \sum_{j = 1}^N R (f_j) \leqslant \sum_{j = 1}^N T (f_j) =
\frac 1 N \left|\{f > 0\}\right| \sum_{j = 1}^N \int f_j = \frac 1 N T (f)
\end {equation*}
implying that $R (f) < \varepsilon$ for any $\varepsilon > 0$ and therefore $R (f) \notin \mathcal V$.  Thus $T$ has no linear submaps.
Apparently, the nature of restictions on $\mathcal V$ necessary for the conclusion of Theorem \ref {linselg} to hold
depends on the~space~$\mathcal S$.

\section {Certain facts about continuity}

\label {cfac}

In this section we introduce and discuss the auxiliary facts of topological nature that we need for the main results of this paper;
see, e.~g., \cite {aubinfrankowska1990} or \cite {rubinov1980} (and also \cite {worth1970}, \cite {huh1974})
for a more detailed exposition of the topological aspects.
Let $\mathcal S$ and $\mathcal W$ be topological spaces.
We say that a map $a : \mathcal S \to 2^{\mathcal W}$ is \emph {upper semicontinous} if for every
closed set $F \subset \mathcal W$ its preimage
\begin {multline*}
a^{-1} (F) = \bigcup_{\weightw \in F} a^{-1} (\weightw) =
\\
\bigcup_{\weightw \in F} \{x \in \mathcal S \mid a (x) \ni \weightw\} =
\left\{ x \in \mathcal S \mid a (x) \cap F \neq \emptyset \right\}
\end {multline*}
under $a$ is also closed.
It is easy to see that $a$ is upper semicontinuous if and only if for every open set $G \subset \mathcal W$
the set $\{ x \in \mathcal S \mid a (x) \subset G \}$ is also open.
This implies that for single-valued maps 
upper semicontinuity is equivalent to the usual continuity.

The following simple proposition, closely following \cite [Proposition 4.4, 4.5] {rubinov1980},
shows that, in the case of regular topology, upper semicontinuity is the same as closedness of the graph for maps with closed values.
Recall that a topological space $\mathcal W$ is called \emph {regular} if every closed set $F \subset \mathcal W$ and every point 
$x \in \mathcal W \setminus F$ can be separated from each other by some open neighbourhoods of $F$ and $x$, i. e.
there exists a couple of open sets $U$ and $V$ in $\mathcal W$ such that $U \ni x$, $V \supset F$ and $U \cap V = \emptyset$.
It is well known that every Hausdorff linear topological space is regular.
\begin {proposition}
\label {grclos}
Suppose that $\mathcal S$ and $\mathcal W$ are topological spaces, $\mathcal W$ is regular and a map
$a : \mathcal S \to 2^{\mathcal W}$ has nonempty closed values.  The following conditions are equivalent.
\begin {enumerate}
\item
$a$ is upper semicontinuous.
\item
The graph $\gr a = \left\{ (x, y) \mid x \in \mathcal S, y \in a (x) \right\}$ of the map $a$ is closed in $\mathcal S \times \mathcal W$.
\end {enumerate}
\end {proposition}
Indeed, suppose that $a$ is upper semicontinuous under the conditions of Proposition \ref {grclos}.
It is sufficient to prove that for every point $(x, z) \notin \gr a$ there exists a neighbourhood of $(x, z)$ that does not intersect $\gr a$.
Since $z \notin a (x)$, by regularity of $\mathcal W$ there exist open sets $V \ni z$ and $G \supset a (x)$ in $\mathcal W$
such that $V \cap G = \emptyset$.
By upper semicontinuity of $a$ the set $U = \{ x \in \mathcal S \mid a (x) \subset G \}$ is open.
Then $U \times V$ is a neighbourhood of $(x, z)$ that does not intersect $\gr a$.

Conversely, let the graph $\gr a$ be closed and $F \subset \mathcal W$ be an arbitrary closed set; we need to show that
$a^{-1} (F) = \left\{ x \in \mathcal S \mid a (x) \cap F \neq \emptyset \right\}$ is also closed.
It is easy to see that $a^{-1} (F)$ is the canonical projection of a closed set $(\mathcal S \times F) \cap \gr a$ onto $\mathcal S$.
Therefore $a^{-1} (F)$ is closed, since
the canonical projection from $\mathcal S \times \mathcal W$ onto $\mathcal S$ is an open map.  The proof of Proposition \ref {grclos} is complete.

The following proposition is used in the proof of Theorems~\ref {everyls}, \ref {riesznecc} and \ref {everylsl1} below.
\begin {proposition}
\label {superctolinc}
Suppose that $\mathcal S$ is a semilinear space equipped with a topology such that for every $x \in \mathcal S$ and open $U \subset \mathcal S$
the set $x + U$ is also open,
$\mathcal W$ is a Hausdorff linear topological space and $\mathcal V$
is the semilinear space of all nonempty compact sets of $\mathcal W$ ordered by inclusion.  Let $T : \mathcal S \to \mathcal V$ be a superlinear map.
If $T$ is upper semicontinuous then the maximal linear submap $S$ of $T$ defined by \eqref {thesubmap} is also upper semicontinuous.
\end {proposition}

First of all, it is easy to verify that the space $\mathcal V$ is order regular so the infinum in \eqref {thesubmap} is well-defined;
see Proposition~\ref {csgs} in Section~\ref {lsolm} below.
Let $x \in \mathcal S$ and $G \subset \mathcal W$ be an open set such that $S (x) \subset G$; it suffices to establish that there exists
a neighbourhood $U \subset \mathcal S$ of $x$ such that
$S (y) \subset G$ for $y \in U$.

Notice that for every lower directed family of compact sets $K_\alpha \subset \mathcal W$ such that $\bigcap_\alpha K_\alpha \subset G$ there is
an index $\alpha$ such that $K_\alpha \subset G$; otherwise $K_\alpha \setminus G$ would have been a centered family of compact sets and therefore
$\bigcap_\alpha K_\alpha \setminus G = \bigcap_\alpha \left(K_\alpha \setminus G\right) \neq \emptyset$.
Applying this observation to \eqref {thesubmap}, we get a decomposition $x = \sum_{j = 1}^n x_j$, $x_j \in \mathcal S$, such that
$\sum_{j = 1}^n T (x_j) \subset G$.  Now it is sufficient to find some open sets $G_j \subset \mathcal W$ satisfying $T (x_j) \subset G_j$
and $\sum_{j = 1}^n G_j \subset G$; once this is done, the sets $U_j = \{ z_j \in \mathcal S \mid T (z_j) \subset G_j \}$
become open neighbourhoods of $x_j$
by the upper semicontinuity of $T$, so the set $U = \sum_{j = 1}^n U_j$
is an open neighbourhood of $x$ such that
$U \subset \{y \in \mathcal S \mid S (y) \subset G\}$.

Let $A = \left\{ \{y_j\}_{j = 1}^n \mid \sum_{j = 1}^n y_j \in G \right\}$.  Since $\mathcal W$ is a linear topological space,
the set $A \subset \mathcal W^n$ is open.  It contains the direct product $P = \prod_{j = 1}^n T (x_j)$ of compact sets which is
a compact set in $A$.  The set $\mathcal B = \left\{ \prod_{j = 1}^n V \mid \text {$V$ is a neighbourhood of $0$ in $\mathcal W$ } \right\}$
is a base of the neighbourhoods of $0$ for the topology of
$\mathcal W^n$, so
$$
X = \{ y + W \mid y \in P,\, W \in \mathcal B,\, y + W \subset A \}
$$
is an open cover of $P$.
Since $P$ is a compact set, there exists a finite subcover
$\left\{y_k + W_k = \left\{y_k^{(j)}\right\}_{j = 1}^n + \prod_{j = 1}^n V_{k, j} \right\}_{k = 1}^M \subset X$ of $P$.
Then $G_j = \bigcap \left\{ \bigcup S_j \mid S_j \subset \left\{y_k^{(j)} + V_{k, j}\right\}_{k = 1}^M,\,\, T (x_j) \subset \bigcup S_j \right\}$
is an intersection of a finite number of open sets, so $G_j$ is an open set in $\mathcal W$ containing $T (x_j)$.
On the other hand, $\prod_{j = 1}^n G_j \subset \bigcup X \subset A$, so
$\sum_{j = 1}^n G_j \subset G$ and thus $\{G_j\}_{j = 1}^n$ is a suitable collection of open sets.
The proof of Proposition \ref {superctolinc} is complete.

In order to work with the issues of continuity of linear selections
it is desirable to have certain natural agreement between semilinear and topological structures.
Let $\mathcal V$ be a semilinear space.  A set $A \subset \mathcal V$ is called \emph {balanced with respect to a point $x \in A$} if
for any $z \in A$ there exists $z' \in A$ such that $x = \frac 1 2 (z + z')$, i. e. $x$ is a midpoint of a segment $[z, z'] \subset A$
for any $z \in A$.  We denote the topological interior of a set $A$ by $\tint A$.
We say that $\mathcal V$ is a \emph {semilinear topological space} if $\mathcal V$ is equipped with a topology satisfying the following conditions.
\begin {enumerate}
\item
The semilinear operations of addition and multiplication by positive constants are continuous.
\item
For every open set $U \subset \mathcal V$ and $x \in \mathcal V$ the set $x + U$ is also open.
\item
$\tint \mathcal V \neq \emptyset$, and 
for every point $x \in \tint \mathcal V$ which is not an additive identity and every open set $U \subset \mathcal V$ containing $x$ there exists an open
set $V \subset U$ which is balanced with respect to $x$, i. e. for every point $x \in \mathcal V$ there exists a base of neighbourhoods of
$x$ that are balanced with respect to $x$.
\end {enumerate}
Naturally, a linear topological space is a semilinear toplogical space.
Throughout the paper,
Properties 2 and 3 of the above definition are directly used only in this section in two critical places apiece.
Note that if $X$ is a normed lattice of measurable functions
then the set $X_+$ of the functions positive a. e. with the topology induced from $X$ is generally \emph {not} a semilinear topological space
as defined above
unless $X = \lclassg {\infty}$ because of property 3.
This is of no relevance for the present paper, even though such semilinear spaces $X_+$ are the focus of it application-wise,
because in the case of cones in a normed lattice continuity of the selections usually follows at once from boundedness and linearity.
However, to give a more interesting example (which is not used in the present work),
there is a natural family of topologies on $X_+$ turning it into a semilinear topological space.
The topologies are defined for every measurable function $h$ satisfying $0 < h \leqslant 1$ by a base of neighbourhoods for $f \in X$
consisting of the sets $\{ g \mid |g - f| < \varepsilon h |f|\}$ for $0 < \varepsilon < 1$,
which for $h = 1$ coincides with the topology induced from $X$ if
$X = \lclassg {\infty}$ and is (much) stronger otherwise.


The following proposition contains an implication from \cite [Theorem~4.1] {rubinov1980}.
\begin {proposition}
\label {funcont}
Suppose that $\mathcal S$ and $\mathcal W$ are semilinear spaces equipped with a topology, and let
$T : \mathcal S \to 2^{\mathcal W} \setminus \{\emptyset\}$ and
$g : \mathcal W \to \mathbb R$ be some mappings.  Define the support function
$g^*_T (x) = \sup_{y \in T (x)} g (y)$ for $x \in \mathcal S$.
If $T$ and $g$ are upper semicontinuous then $g^*_T$ is also upper semicontinuous.
\end {proposition}
Indeed, suppose that $\lambda \in \mathbb R$ is an arbitrary number and maps $T$ and $g$ are upper semicontinuous.
Then
$$
\{x \in \mathcal S \mid g^*_T (x) < \lambda\} =
\bigcup_{0 < \varepsilon < 1} \left\{x \in \mathcal S \mid T (x) \subset \{ y \in \mathcal W \mid g (y) < \lambda - \varepsilon \} \right\}
$$
is an open set, so $g^*_T$ is upper semicontinuous.

Observe that under the assumptions of Proposition~\ref {funcont} if the map $T$ and the functional $g$ are both linear then the support function
$g^*_T$ is also linear. The following proposition allows us to obtain full continuity of the support function $g^*_T$ from
upper semicontinuity of $T$ and $g$ alone.
\begin {proposition}
\label {scfp}
Suppose that $\mathcal S$ is a semilinear topological space and
$f : \mathcal S \to \mathbb R$ is an affine functional on $\mathcal S$.
If $f$ is upper (or lower) semicontinuous then $f$ is continuous.
\end {proposition}
Indeed, suppose that under the conditions of Proposition \ref {scfp} linear functional $f$ is upper semicontinuous
(the case of lower semicontinuous $f$ is naturally reduced to this case by replacing $f$ with $-f$).
This means that for any $x \in \tint \mathcal S$ which is not an additive identity and $\varepsilon > 0$
there exists an open neighbourhood $U \subset \mathcal S$ of $x$ balanced with respect to $x$ such that
\begin {equation}
\label {e153}
f (u) \leqslant f (x) + \varepsilon
\end {equation}
for any $u \in U$.  Since $U$ is balanced with respect to $x$ for any $z \in U$, there exists some
$y \in U$ such that $x = \frac 1 2 (y + z)$, and therefore
$$
f (x) = \frac 1 2 (f (y) + f (z)).
$$
Substituting this into \eqref {e153}
with $u = z$ and $u = y$
yields $|f (z) - f (y)| \leqslant 2 \varepsilon$.
Thus $|f (z) - f (x)| = \frac 1 2 |f (z) - f (y)| \leqslant \varepsilon$ for any $z \in U$.
This implies that $f$ is continuous on a set
$$
[\tint \mathcal S]_+ = \left\{ x \in \tint \mathcal S \mid \text {$x$ is not an additive identity}\right\}.
$$
Note that $\tint \mathcal S$ is dense in $\mathcal S$ because $\mathcal S$ is a convex set.
If $[\tint \mathcal S]_+$ is dense in $\mathcal S$ then it follows that $f$ is continuous on the entire semilinear space $\mathcal S$.
If, on the other hand, $[\tint \mathcal S]_+$ is not dense in $\mathcal S$, then there is an additive identity $0 \in \mathcal S$ and
the set $[\tint \mathcal S]_+ \cup \{0\}$ is dense in $\mathcal S$.  In this case the set
$\{0\}$ is open in $\mathcal S$, so $f$ is automatically continuous at $0$ and therefore in all $\mathcal S$.
The proof of Proposition \ref {scfp} is complete.

\begin {proposition}
\label {swfl}
Suppose that $\mathcal S$, $\mathcal W$ are semilinear spaces,
$g : \mathcal W \to \mathbb R$ and $T : \mathcal S \to \mathcal V \subset \mathcal C (\mathcal W)$
are linear maps,
and $g$ is bounded from above on the sets from $\mathcal V$, i.~e.
$\sup_{v \in V} g (v) \leqslant b_V$ with some $b_V \in \mathbb R$ for any $V \in \mathcal V$.
Then the support function $g^*_T$ is linear.  If, additionally, $\mathcal S$ is a semilinear topological space and $\mathcal W$ is equipped with
a topology such that both $T$ and $g$ are upper semicontinuous then the support function $g^*_T$ is continuous.
\end {proposition}
Linearity of the support function $g^*_T$ is trivial:
\begin {multline*}
g^*_T (\lambda x + \mu y) = \sup_{w \in T (\lambda x + \mu y)} g (w) =
\\
\sup_{u \in T (x), v \in T (y)} g (\lambda u + \mu v) =
\sup_{u \in T (x), v \in T (y)} \left[ \lambda g (u) + \mu g (v) \right] =
\\
\sup_{u \in T (x)} \lambda g (u) + \sup_{v \in T (y)} \mu g (v) = \lambda g^*_T (x) + \mu g^*_T (y)
\end {multline*}
for any $x, y \in \mathcal S$ and $\lambda, \mu > 0$.  The boundedness from above assumed of $g$ on the sets $\mathcal V$
guarantees that $g^*_T (x)$ is finite for all $x \in \mathcal S$, i. e. $g^*_T$ is a linear functional.
Therefore, Propositions \ref {funcont} and \ref {scfp} combined yield Proposition \ref {swfl}.


\section {Linear selections of linear maps}

\label {lsolm}

In this section we state the main results of this paper concerning linear selections of linear or superlinear maps.
A most interesting case to which Theorem \ref {linselg} can be applied
is a semilinear space $\mathcal V \subset \mathcal C (\mathcal W)$ of convex sets in some semilinear space $\mathcal W$
ordered by inclusion.
The semilinear space $\mathcal V$ has to be order regular for Theorem \ref {linselg} to be applicable in this setting.
This implies some restrictions on the collection of convex sets $\mathcal V$.
The most natural restriction is compactness, but in some situations it is too strong.
However, compactness can often be safely replaced by a weaker property.
We say that $\mathcal V \subset \mathcal C (\mathcal W)$ \emph {has compact type} if for every centered family of sets
$A_\alpha \in \mathcal V$ their intersection $\bigcap_\alpha A_\alpha$ also belongs to~$\mathcal V$.
This implies lower completenes of $\mathcal V$, and order regularity is easily verified in this case;
see the proof of Proposition~\ref {csgs} below.  Theorem~\ref {linselg} then provides for a superlinear map $T$
the greatest linear submap $S$ taking values in $\mathcal V$, so the problem of linear selections of superlinear maps gets reduced
to the problem of linear selections of linear maps by Theorem~\ref {linselg}.

In order to work with linear selections of linear maps taking values in $\mathcal V \subset \mathcal C (\mathcal W)$
we need a dual structure $\mathcal D$ of linear functionals on $\mathcal W$ satisfying certain properties allowing us to characterize points of sets
from $\mathcal V$ by values of functionals from $\mathcal D$ at these points and vice versa.
Suppose that $\mathcal W$ is a semilinear space, $\mathcal D \subset \mathcal W'$ is a set of linear
real-valued functionals on $\mathcal W$ and
$\mathcal V \subset \mathcal C (\mathcal W)$.  We denote the affine hyperplanes of a linear functional $f$ on $\mathcal W$ by
$
f^{-1} (c) = \{ x \in \mathcal W \mid f (x) = c \}.
$
A linear functional $f$ on $\mathcal D$ is said to be
\emph {consistent with $\mathcal V$} if for every $A \in \mathcal V$ and $c \in \mathbb R$ nonemptiness of
$f^{-1} (c) \cap A$ implies that $f^{-1} (c) \cap A \in \mathcal V$; $f$ is said to be \emph {defining for $\mathcal V$}
if $f (A)$ is a segment $[a, b]$ for all $A \in \mathcal V$.
We say that the set of functionals $\mathcal D$ is \emph {exhaustive for $\mathcal W$} if values $\{ f (x) \}_{f \in \mathcal D}$ uniquely
identify
every point $x \in \mathcal W$.

The main example of a space $\mathcal V \subset \mathcal C (\mathcal V)$ good enough for the main results
and all but one application presented in this paper
is given in the following simple proposition, which outlines the well-known properties of the space of nonempty compact convex sets
of a locally convex linear topological space.
\begin {proposition}
\label {csgs}
Suppose that
$\mathcal W$ is a locally convex Hausdorff linear topological space and a semilinear space $\mathcal V = \mathcal V (\mathcal W)$
consists of all nonempty compact convex sets in $\mathcal W$.
Then $\mathcal V$ is a subspace of the semilinear space $\mathcal C (\mathcal W)$, i. e. $\mathcal V$ is closed under the semilinear operations of
$\mathcal C (\mathcal W)$.
$\mathcal V$ has compact type and is order regular.
The set $\mathcal W'$ of all continuous linear functionals on $\mathcal W$ is exhaustive for $\mathcal W$, and all $f \in \mathcal W'$
are defining for and consistent with $\mathcal V$.
\end {proposition}
Indeed, if $A, B \in \mathcal V$ then $A + B$ is convex.  If $z_\alpha \in A + B$ is a net then
$z_\alpha = x_\alpha + y_\alpha$ for some nets $x_\alpha \in A$ and $y_\alpha \in B$.  Since $A$ and $B$ are compact, we can replace
the nets by some cofinal ones so that $x_\alpha \to x$ and $y_\alpha \to y$ for some $x \in A$ and $y \in B$.
This means that $z_\alpha \to x + y \in A + B$, so $A + B$ is compact and thus closed.  We have just shown that $A + B \in \mathcal V$ for any
$A, B \in \mathcal V$.
Closedness of $\mathcal V$ under multiplication by positive constants is trivial.
Therefore, $\mathcal V$ is a subspace of $\mathcal C (\mathcal W)$.
$\mathcal V$ has compact type because intersection of a centered family of nonempty compact convex sets is also a nonempty compact convex set.
The inclusion order $\subset$ is lower complete in $\mathcal V$
for the same reason, and for any lower directed set $\mathfrak A \subset \mathcal V$ we have $\inf \mathfrak A = \bigcap \mathfrak A$.
Let us verify the remaining part of the order regularity property for $\mathcal V$.
Suppose that
$\mathfrak A, \mathfrak B \subset \mathcal V$ are lower directed and $Z \in \mathcal V$ satisfies $Z \subset A + B$ for all
$A \in \mathfrak A$ and $B \in \mathfrak B$.  It is evident that for a fixed $B \in \mathcal V$ the set
$\mathfrak A + B = \{ A + B \mid A \in \mathfrak A\}$ is lower directed, and $\inf (\mathfrak A + B) = (\inf \mathfrak A) + B$
since $\inf = \bigcap$.  Taking infinums over $B \in \mathfrak B$ in the same way shows that
$Z \subset \inf \mathfrak A + \inf \mathfrak B$, which proves that $\mathcal V$ is order regular.
The fact that the set $\mathcal W'$ is exhaustive easily follows from a separation theorem.
Finally, if $f \in \mathcal W'$ is a continuous linear functional then $f (A)$ is compact and convex for any $A \in \mathcal V$
as an image of a compact convex set, thus $f (A) = [a, b]$ for some $a, b \in \mathbb R$.  If $A \in \mathcal V$ and $c \in f (A)$ then
$f^{-1} (c) \cap A$ is a nonempty convex closed subset of a compact set $A$, so $f^{-1} (c) \cap A \in \mathcal V$.
The proof of Proposition \ref {csgs} is complete.

We remark in passing that in Proposition~\ref {csgs} compact sets can be replaced by nonempty bounded sets that are closed in measure
if $\mathcal W$ is a Banach ideal lattice of measurable functions on $\mathbb Z$ satisfying the Fatou property,
as will be discussed in detail in Section \ref {scm} below.

Recall that in the present setting Theorem \ref {linselg} characterizes linear submaps $A : \mathcal S \to \mathcal V$
of a superlinear map $T$ under certain assumptions
imposed on $\mathcal S$ and $\mathcal V \subset \mathcal C (Y)$ by the explicit formula \eqref {thesubmap} for the greatest linear submap $S$ of $T$.
The next step in the study of linear selections of $T$ is to describe linear selections $L (A)$ of such a linear submap $A$.
Let $V_A (x) = \{ a (x) \mid a \in L (A) \}$ be the set of all possible values of linear selections of $A$ at some point $x \in \mathcal S$.
In \cite [Theorem 2] {zaslavsky1981en} it was shown (using the Zorn lemma)
for the standard case of $Y$ being locally convex and $\mathcal V$ consisting of all compact
convex sets in $Y$ that the closure of $V_A (x)$ in the weak topology of $Y$ coincides with $A (x)$ and thus the sets $V_A (x)$ are nonempty;
consequently, $L (A)$ is also nonempty.
While this result provides a positive answer to the question about the existence of linear selections for superlinear maps
in this setting,
it does not by itself say everything about the set $V_A (x)$, specifically whether $V_A (x) = A (x)$.
In other words, is it possible to extend every single-valued linear submap $l \leqslant A$ defined on a ray spanned by $x$ to a submap
of $A$ defined
on the entire $\mathcal S$?
Since $V_A (x)$ is convex, it \emph {might} be possible to prove that $V_A (x)$ is closed in the weak topology, which would imply
that $V_A (x) = A (x)$.
However, there is a more direct approach.
Application of a simple technique to the general idea of the proof of \cite [Theorem 2] {zaslavsky1981en}
allows us to establish that $V_A (x) = A (x)$ in a more general setting.
Moreover, the same technique leads to an \emph {explicit} inductive construction of a linear selection of $A$ passing through a given point in $A (x)$
in addition to the quick proof based on the Zorn lemma, which allows us to dispose of the compact type assumption
if the set of defining functionals is finite, which is the case if $\mathcal Y$ is a finite-dimensional topological space.

\begin {theorem}
\label {linselnd}
Suppose that $\mathcal S$ and $\mathcal W$ are semilinear spaces,
$\mathcal V$ is a subspace of the semilinear space $\mathcal C (\mathcal W)$,
$\mathcal D \subset \mathcal W'$ is an exhaustive set of linear functionals that are defining for and consistent with $\mathcal V$,
and either $\mathcal V$ has compact type or $\mathcal D$ is finite.
Then for every linear map $A : \mathcal S \to \mathcal V$, $x \in \mathcal S$ and $y \in A (x)$ there exists a linear selection $a$ of $A$
satisfying $a (x) = y$.  Suppose, in addition, that $\mathcal S$ and $\mathcal W$ are equipped with topologies such that all sets in $\mathcal V$ are closed,
$\mathcal W$ is a regular topological space, $\mathcal S$ is a semilinear topological space,
all functionals in $\mathcal D$ are continuous and $A$ is upper semicontinuous.  Then
there exists a continuous linear selection $a$ of $A$ satisfying $a (x) = y$.
\end {theorem}
In Section \ref {lslm} below we provide a simple proof of Theorem \ref {linselnd} (for the case of $\mathcal V$ having compact type)
based on the Zorn lemma,
which is essentially a refinement of \cite [Theorem 2] {zaslavsky1981en}.
We also give a more explicit construction of the linear selection with required properties (that covers the case of finite $\mathcal D$),
which is based on a certain parametrization of points
in an arbitrary set $A \in \mathcal V$.
This construction is based on transfinite recursion in the general case, so it also involves the Zorn lemma indirectly;
it is truly explicit only when the set of defining functionals $\mathcal D$ for $\mathcal W$ is at most countable.

Note that for general spaces $\mathcal S$ and $\mathcal V$ not every linear set-valued map has a linear selection.
Perhaps, the simplest example is
the two-point semilinear space $\mathcal S = \{0, 1\}$ mentioned before and the positive ray~$\mathcal W = \{\lambda \mid \lambda > 0\}$.
Let $T : \mathcal S \to \mathcal C (\mathcal W)$ be a map defined by $T (f) = \mathcal W$ for all $f \in \mathcal S$.
Map $T$ is linear, but it is easy to see that in this case there are no linear maps $S : \mathcal S \to \mathcal W$ at all.
Section \ref {sotctc} contains an example that shows that the compact type assumption
imposed on $\mathcal V$ cannot be dropped from the statement of Theorem \ref {linselnd}.
It is easy to see, however, that no conditions on $\mathcal V$ of this sort are actually necessary if,
for example, $\mathcal S = \mathbb R_+$.
More generally, by Proposition \ref {conebise} in Section \ref {conebases} below the non-topological conclusion of Theorem \ref {linselnd}
is valid if $\mathcal S$ has a cone-basis and $\mathcal V$ is merely a semilinear subspace of $\mathcal C (\mathcal W)$.
However, as it will be evident in Section \ref {sotctc},
this approach is rather limited when infinite-dimensional spaces are involved.
An interesting particular case of upper semicontinuous maps acting on the space of measures $\mathcal M (K)$
with virtually no restrictions placed on $\mathcal V$ is treated below in Theorem~\ref {csselectm} of Section~\ref {mosom}.
Another interesting particular case of bounded upper semicontinuous maps acting on the space $\left[\lclassg {1}\right]_+$
of nonnegative a.~e. summable functions on a separable measurable space into convex sets
of a Banach ideal space of measurable functions satisfying
the Fatou property that are closed in measure is treated in Section~\ref {ssvm}.

The following result, which is a simple consequence of Theorem \ref {linselnd},
is inspired by (and related to) \cite [Theorem 2] {gorokhovik2008}; see also Theorem \ref {asvmch} in Section \ref {ascm} below.
\begin {theorem}
\label {lsvmch}
Suppose that $\mathcal S$ and $\mathcal W$ are semilinear spaces,
$\mathcal V$ is a subspace of the semilinear space $\mathcal C (\mathcal W)$,
$\mathcal D \subset \mathcal W'$ is an exhaustive set of linear functionals that are defining for and consistent with $\mathcal V$,
and either $\mathcal V$ has compact type or $\mathcal D$ is finite.  Let $T : \mathcal S \to \mathcal V$ be a superlinear set-valued map.
The following conditions are equivalent.
\begin {enumerate}
\item
$T$ is linear.
\item
For any $x \in \mathcal S$ and $y \in T (x)$ there exists a linear selection~$S$ of~$T$ such that $S (x) = y$.
\end {enumerate}
\end {theorem}
Indeed, implication $1 \Rightarrow 2$ follows from Theorem \ref {linselnd}.
Conversely, suppose that condition 2 is satisfied.  We need to show that for arbitrary $a, b \in \mathcal S$ and $\lambda, \mu > 0$
we have $T (\lambda a + \mu b) = \lambda T (a) + \mu T (b)$.  Since $T$ is superlinear, it suffices to establish that 
$T (\lambda a + \mu b) \subset \lambda T (a) + \mu T (b)$.
By condition 2 for $x = \lambda a + \mu b$ and any $y \in T (x)$ there exists a linear selection $S$ of $T$ satisfying $S (x) = y$.
This means that
$$
y = S (x) = \lambda S (a) + \mu S (b) \in \lambda T (a) + \mu T (b)
$$
as claimed.  The proof of Theorem \ref {lsvmch}
is complete.

Theorems \ref {linselg} and \ref {linselnd} together yield the following result.
\begin {theorem}
\label {everyls}
Suppose that $\mathcal S$ and $\mathcal W$ are semilinear spaces, $\mathcal S$ has the Riesz decomposition property,
$\mathcal V$ is a subspace of the semilinear space $\mathcal C (\mathcal W)$ that has compact type and
$\mathcal D \subset \mathcal W'$ is an exhaustive set of linear functionals that are defining for and consistent with $\mathcal V$.
Let $T : \mathcal S \to \mathcal V$ be a superlinear map.
Then $T$ admits a linear selection.
Additionally, fix some points $x \in \mathcal S$ and $y \in T (x)$.
Then the following conditions are equivalent.
\begin {enumerate}
\item
There exists a linear selection $a$ of $T$ satisfying $a (x) = y$.
\item
$y \in \sum_{j = 1}^N T (x_j)$ for any $N \in \mathbb N$ and $\{x_j\}_{j = 1}^N \subset \mathcal S$ satisfying $x = \sum_{j = 1}^N x_j$.
\end {enumerate}
If, in addition, $\mathcal S$ is a semilinear topological space, $\mathcal W$ is a locally convex Hausdorff linear topological space,
$\mathcal V$ consists of all nonempty compact convex sets in $\mathcal W$,
all functionals in $\mathcal D$
are continuous and $T$ is upper semicontinuous, then the selection in condition 1 can be assumed to be continuous.
\end {theorem}
The topological clause of Theorem \ref {everyls} follows from Proposition \ref {superctolinc}, which shows that
under the stated conditions upper semicontinuity of $T$ implies upper semicontinuity of its maximal linear submap determined by the formula
\eqref {thesubmap}.

\section {Tomographical coordinates}

\label {lslm}
In this section we prove Theorem \ref {linselnd} and discuss some related issues.
Suppose that under the conditions of Theorem \ref {linselnd}
we are given a linear map $A : \mathcal S \to \mathcal V$ and some points $x \in \mathcal S$, $y \in A (x)$;
we need to construct a suitable linear selection $a$ of $A$ satisfying $a (x) = y$.

First, every $f \in \mathcal D$ is bounded on any $K \in \mathcal V$, so the \emph {support functions}
$
R_f (K) = \sup_{z \in K} f (z)
$
and
$
L_f (K) = \inf_{z \in K} f (z)
$
are well-defined for all $K \in \mathcal V$
and take finite values.  Surely, $L_f = -R_{-f}$ and $R_f (K) = f^*_K$, $L_f (K) = -(-f)^*_K$ in terms of functions
defined in Proposition \ref {funcont} above.
Let the \emph{$\theta$-section of $K$ by $f$} be defined for $0 \leqslant \theta \leqslant 1$ by
$$
\sect_f (\theta, K) = K \cap f^{-1} \left( (1 - \theta) L_f (K) + \theta R_f (K) \right).
$$
Since $f$ is a defining functional consistent with $\mathcal V$, $\sect_f (\theta, K) \in \mathcal V$
for all $\theta$ and $\bigcup_{0 \leqslant \theta \leqslant 1} \sect_f (\theta, K) = K$; this union is disjoint unless
$L_f (K) = R_f (K)$.  The following proposition shows that taking such sections of a set from $\mathcal V$ is a well-behaved operation.
\begin {proposition}
\label {linind}
Let $\mathcal S$ and $\mathcal W$ be semilinear spaces, $\mathcal V$ be a semilinear subspace of $\mathcal C (\mathcal W)$,
$f$ be a linear functional on $\mathcal V$ defining for and consistent with $\mathcal V$,
and let $T : \mathcal S \to \mathcal V$ be a linear map.  Then for every fixed $0 \leqslant \theta \leqslant 1$ map
$\sect_{\theta, f} T$ defined by
$\sect_{\theta, f} T (x) = \sect_f (\theta, T (x))$, $x \in \mathcal S$, is a linear submap of $T$.
Suppose, additionally, that $\mathcal S$ and $\mathcal W$ are endowed with topologies so that all sets in $\mathcal V$ are closed,
$\mathcal W$ is a regular topological space, $\mathcal S$ is a semilinear topological space, $f$ is continous
and $T$ is upper semicontinuous.  Then the map $\sect_{\theta, f} T$ is also upper semicontinuous.
\end {proposition}
Indeed,
$F (x) = \left( (1 - \theta) L_f (T (x)) + \theta R_f (T (x)) \right)$ defined for $x \in \mathcal S$ is a linear functional by Proposition \ref {swfl}.
Therefore
$$
x \mapsto \sect_{\theta, f} T (x) = \sect_f (\theta, T (x)) = T (x) \cap \{y \in \mathcal W \mid f (y) = F (x) \},
$$
$x \in \mathcal S$,
is a linear map, because it is an intersection of two linear maps.
Now suppose that the topological assumptions of Proposition~\ref {linind} are satisfied.  Then
$F$ is a continuous linear functional by the same Proposition~\ref {swfl}, so the graph of the map
the map $x \mapsto \{y \in \mathcal W \mid f (y) = F (x) \}$ is closed. By Proposition \ref {grclos} the graph of the map $T$ is also closed, so the graph of
the map $\sect_{\theta, f} T$ is closed as an intersection of two closed graphs; in particular, values of $\sect_{\theta, f} T$ are closed.
By the same Proposition \ref {grclos} this means that $\sect_{\theta, f} T$ is upper semicontinuous. The proof of Proposition \ref {linind} is complete.

We are now ready to give a simple proof of Theorem \ref {linselnd} for the case of $\mathcal V$ having compact type;
the other case is covered by a more elaborate construction that will be presented after that.
Suppose that the assumptions of Theorem \ref {linselnd} hold true.  Let $\mathcal L$ be the set of all linear submaps $B$ of $A$
(i. e. $B (z) \subset A (z)$ for all $z \in \mathcal S$) satisfying $y \in B (x)$.
We induce a partial order $\prec$ in $\mathcal L$ by the reverse set inclusion order on the graphs of functions from $\mathcal L$, i. e.
for any $D, E \in \mathcal L$ we have
$D \prec E$ if and only if $\gr E \subset \gr D$. Let $\mathfrak D$ be a linearly ordered set $\mathfrak D = \{D_\alpha\} \subset \mathcal L$.
Then the map $D$ defined by $D (x) = \bigcap_\alpha D_\alpha (x)$ for $x \in \mathcal S$ also belongs to $\mathcal L$ by the compact type of $\mathcal V$,
so $D$ is an upper bound for $\mathfrak D$.  Thus we can apply the Zorn lemma to obtain a maximal element $M \in \mathcal L$.
Let us verify that $M$ is single-valued.  If, on the contrary, $M (z)$ has at least two points $a, b \in M (z)$ for some $z \in \mathcal S$,
then by the assumptions there exists a functional $f \in \mathcal D$ such that $f (a) \neq f (b)$.
Let
$\theta = \frac {f (y) - L_{f} (M (x))} {R_{f} (M (x)) - L_{f} (M (x))}$ if $R_{f} (M (x)) \neq L_{f} (M (x))$ and $\theta = 0$ otherwise.
Then by Proposition \ref {linind} map $N = \sect_{\theta, f} M$ also belongs to $\mathcal L$, and since $L_f (N (z)) = R_f (N (z))$, $N$ is a proper
submap of $M$, a contradiction.  Thus $M$ is single-valued, and the non-topological claims of Theorem \ref {linselnd} are satisfied.
The topological claims are obtained by repeating this argument with $\mathcal L$ being the set of all upper semicontinuous linear submaps of $A$.
The simple proof of Theorem \ref {linselnd} is complete.

We now present a construction which leads to another proof of Theorem \ref {linselnd}.
Suppose that
$K \in \mathcal V$ and $z \in K$ under the assumptions of Theorem \ref {linselnd}.
We want to characterize the position of the point $z$ relative to $K$ in terms of values of functionals
$\mathcal D$ in a way that would allow us to take advantage of Proposition \ref {linind} to build a linear selection of $A$ with graph passing through a
prescribed point.  In order to do this in the general case we need to invoke transfinite recursion twice, which (at least in its full generality)
can be avoided and the construction becomes much more transparent
if $\mathcal D$ is at most countable.  This particular case is also often sufficient for applications, so we are going to treat it first before
presenting the construction for the general case (which does not depend on the particular case of at most countable $\mathcal D$).

Let us assume that $\mathcal D = \{f_j \}_{j \in I}$, where either $I = \{1 \leqslant j \leqslant N\}$ for some $N \in \mathbb N$ or $I = \mathbb N$.
Define the \emph {tomographical coordinates $\Theta = \{\theta_f\}_{f \in \mathcal D} \subset [0, 1]^{\mathcal D}$
of a point $z \in K$ relative to $K \in \mathcal V$} inductively as follows: $K_0 = K$,
$$
\theta_{f_j} = \frac {f_j (z) - L_{f_j} (K_{j - 1})} {R_{f_j} (K_{j - 1}) - L_{f_j} (K_{j - 1})}
$$
and $K_j = \sect_{\theta_j, f_j} K_{j - 1}$ for $j \in \mathbb N$.
Here, as usual, $\frac 0 0 = 0$.

Conversely, since every $f \in \mathcal D$ is defining for and consistent with $\mathcal V$ and $\mathcal V$ has compact type,
for every coordinates $\Theta = \{\theta_f\}_{f \in \mathcal D} \subset [0, 1]^{\mathcal D}$
and $B \in \mathcal V$ we can define a sequence of sets $B_j \in \mathcal V$ by $B_0 = B$ and
$B_j = \sect_{\theta_j, f_j} B_{j - 1}$ for $j \in I$.
This sequence is nonincreasing, so $C_{B, \Theta} = \bigcap_{j \in I} B_j \in \mathcal V$ by the compact type of $\mathcal V$ if $\mathcal D$ is infinite,
otherwise we just have $C_{B, \Theta} = B_N \in \mathcal V$.
Observe that
$$
f_j \left(C_{B, \Theta}\right) = (1 - \theta_j) L_{f_j} (B_{j - 1}) + \theta_j R_{f_j} (B_{j - 1})
$$
for all $j \in I$.
By assumptions $\mathcal D$ is exhaustive for $\mathcal V$,
so the set $C_{B, \Theta}$ consists of a single point $B (\Theta) \in B$ uniquely identified by the coordinates
$\Theta$.  It is easy to see that if $B = K$ then $B_j = K_j$, where $K_j$ are the sets from the construction of the coordinates $\Theta$ above.
Therefore $K (\Theta) = z$, and thus every point in $K$ can be recovered from its tomographical coordinates.

We now prove Theorem \ref {linselnd} in the case of at most countable $\mathcal D$.
Let $\Theta = \{\theta_j\}$ be the tomographical coordinates of $y$ relative to $A (x)$.
For any $B \in \mathcal V$ the construction above yields a sequence of sets $B_j$ such that $B_0 = B$ and
$B_j = \sect_{\theta_j, f_j} B_{j - 1}$ for $j \in \mathbb N$, and $\{B (\Theta)\} = \bigcap_{j \in I} B_j$.
Proposition~\ref {linind} implies that for all $j \in I$ maps
$A_j (z) = [A (z)]_j (\Theta)$, $z \in \mathcal S$, are linear,
so $z \mapsto [A (z)] (\Theta) = \bigcap_{j \in I} [A (z)]_j (\Theta)$
is also a linear map because it is an intersection of linear maps.
Therefore $a : \mathcal S \to \mathcal W$, $a (z) = \Theta (A (z)) \in A (z)$, is a linear mapping which is a linear selection of $A$
satisfying $a (x) = y$.
Now suppose that the topological assumptions of Theorem \ref {linselnd} are satisfied.
Proposition \ref {linind} together with Proposition~\ref {grclos} implies that all functions $A_j$ have closed graphs, so the graph of $a$
is closed, which by Proposition \ref {grclos} means that $a$ is continuous.  The proof of Theorem \ref {linselnd} for the case of at most countable
$\mathcal D$ is complete.

Let us now extend the construction and application of tomographical coordinates to the case of arbitrary cardinality of the set $\mathcal D$.
In order to do this we have to resort to transfinite recursion; the relevant notions and terminology can be found,
for example, in \cite [Section 18] {halmos1960}.
By the Zermelo theorem we can assume that the set $\mathcal D$ is well-ordered with respect to an order $\prec$.  Denote by
$\mathcal I_f = \{ g \in \mathcal D \mid g \prec f \}$ the initial segment of~$\mathcal D$ for~$f \in \mathcal D$.

We are now going to define the coordinates $\Theta = \{\theta_f\}_{f \in \mathcal D} \subset [0, 1]^{\mathcal D}$ of a point
$z \in K$ relative to $K \in \mathcal V$.
Let $\mathcal V_{z, K} = \{ A \in \mathcal V \mid z \in A \subset K\} \subset \mathcal V$,
let $\mathcal F_f$ be
the set of functions $\mathfrak f : \mathcal I_f \to \mathcal V_{z, K} \times [0, 1]$ with its respective components labelled
by $\mathfrak f = (\mathfrak f_{\mathcal V}, \mathfrak f_\theta)$ and
let $\mathcal F = \bigcup_{f \in \mathcal D} \mathcal F_f$.
We now define a \emph {sequence function} $\mathcal N : \mathcal F \to \mathcal V_{z, K} \times[0, 1]$
with its respective components labelled by $\mathcal N = (\mathcal N_{\mathcal V}, \mathcal N_\theta)$
as follows.
Take $\mathfrak f \in \mathcal F_f$, $f \in \mathcal D$, and
define a set $K_{\mathfrak f}$ by $K_{\mathfrak f} = \bigcap_{g \prec f} \mathfrak f_{\mathcal V} (g)$ if $\mathcal I_f \neq \emptyset$ and
by $K_{\mathfrak f} = K$ otherwise.
The compact type of $\mathcal V$ implies that $K_{\mathfrak f} \in \mathcal V_{z, K}$, so
the succeeding coordinate
$$
\theta_{\mathfrak f} = \frac {f (z) - L_{f} (K_{\mathfrak f})} {R_{f} (K_{\mathfrak f}) - L_{f} (K_{\mathfrak f})}
$$
is well-defined.  We complete the definition of the sequence function $\mathcal N$ by setting
$\mathcal N (\mathfrak f) = \left( \sect_{\theta, f} (K_{\mathfrak f}), \theta_{\mathfrak f} \right)$.
By the Transfinite Recursion Theorem there exists a unique function $U : \mathcal D \to \mathcal V_{z, K} \times [0, 1]$
such that $U (f) = \mathcal N (U_f)$ for all $f \in \mathcal D$, where $U_f$ is the restriction of $U$ onto the initial segment
$\mathcal I_f$.
We label the components of $U$ as $U = (U_{\mathcal V}, U_{\theta})$ and define the tomographical coordinates
$\Theta = \Theta (z, K)$ of the point $z \in K$ relative to $K$ by $\theta_f = U_\theta (f)$ for all $f \in \mathcal D$.

Now suppose that we are given some coordinates $\Theta = \{\theta_f\}_{f \in \mathcal D} \in [0, 1]^{\mathcal D}$ and a set $B \in \mathcal V$.
Let $\mathcal G_f$ be the set of functions
$$
\mathfrak g : \mathcal I_f \to \mathcal V_B  = \{ A \in \mathcal V \mid A \subset B \}
$$
for all $f \in \mathcal D$ and let
$\mathcal G = \bigcup_{f \in \mathcal D} \mathcal G_f$.
Define a sequence function $\mathcal M : \mathcal G \to \mathcal V_B$
for $\mathfrak f \in \mathcal G_f$, $f \in \mathcal D$
by $M (\mathfrak g) = \sect_{\theta_f, f} \bigcap_{g \prec f} \mathfrak g (g)$ if $\mathcal I_f \neq \emptyset$ and $M (\mathfrak g) = B$ otherwise.
Then by the Transfinite Recursion Theorem there exists a unique function $V : \mathcal D \to \mathcal V_B$ satisfying
$V (f) = \mathcal M (V_f)$ for all $f \in \mathcal D$, $V_f$ being the restriction of $V$ onto $\mathcal I_f$.
An application of the Transfinite Induction Theorem easily shows that if $f, g \in \mathcal D$ and $f \prec g$ then $V (f) \subset V (g)$,
which means that $\{V (f)\}_{f \in \mathcal D}$ is a centered family of sets in $\mathcal V$, and so
$C_{B, \Theta} = \bigcap_{f \in \mathcal D} V (f)$ belongs to $\mathcal V$ by the compact type of $\mathcal V$.
Let us define $B_f = \bigcap_{g \prec f} V (g)$ for $\mathcal I_f \neq \emptyset$ and $B_f = B$ otherwise.
Surely $f \left(C_{B, \Theta}\right) = (1 - \theta_f) L_{f} (B_{f}) + \theta_f R_{f} (B_{f})$ for all $f \in \mathcal D$.
Since $\mathcal D$ is exhaustive for $\mathcal V$ by the assumptions, the set $C_{B, \Theta}$ consists of a single point $B (\Theta) \in B$
uniquely identified by the coordinates $\Theta$.  If under the previous assumptions $B = K$ and $\Theta = \Theta (z, K)$, then 
another application of the Transfinite Induction Theorem easily shows that $B_f = U_{\mathcal V} (f) \ni z$ for all $f \in \mathcal D$,
so $B (\Theta) = z$, i. e. every point $z$ in a set $K \in \mathcal V$ can be recovered from
values of its tomographical coordinates $\Theta (z, K)$.

We now prove Theorem \ref {linselnd} in the general case.
Let $\Theta = \{\theta_f\}_{f \in \mathcal D}$ be the tomographical coordinates of $y$ relative to $A (x)$.
For any $B \in \mathcal V$ the above construction yields a collection of sets $B_f (\Theta) = B_f$ such that
$B_f = \sect_{\theta_f, f} \bigcap_{g \prec f} B_g$ if $\mathcal I_f \neq \emptyset$ and $B_f = B$ otherwise.
Then proposition \ref {linind} and the Transfinite Induction Theorem together show that $z \mapsto [A (z)]_f (\Theta)$
is a linear map defined for all $z \in \mathcal S$, so the map
$z \mapsto [A (z)] (\Theta) = \bigcap_{f \in \mathbb D} [A (z)]_f (\Theta)$
is also linear as an intersection
of linear maps, and we also have $[A (x)] (\Theta) = y$ by the construction of the tomographical coordinates.
Therefore $a : \mathcal S \to \mathcal W$, $a (z) = \Theta (A (z)) \in A (z)$, is a linear map which is a linear selection of
the map $A$
satisfying $a (x) = y$.
Now suppose that the topological assumptions of Theorem \ref {linselnd} are satisfied.
Note that $a (z) = \bigcap_{f \in \mathcal D} A_f (z)$ for $z \in \mathcal S$,
where $A_f = \sect_{\theta_f, f} \bigcap_{g \prec f} A_g$ if $\mathcal I_f \neq \emptyset$ and $A_f = A$ otherwise.
Proposition \ref {linind} together with Proposition \ref {grclos}
and the Transfinite Induction Theorem show that all functions $A_f$ have closed graphs, so the graph of the map $a$
is also closed, which by Proposition \ref {grclos} means that the map $a$ is continuous.
The proof of Theorem \ref {linselnd} for the general case is complete.

\section {Affine selections of convex maps}

\label {ascm}
In this section we show how the main results for linear selections of superlinear maps
naturally imply corresponding results for affine selections of convex maps (and sometimes vice versa).
These results are of independent interest (see, e.~g., \cite {smajdor1990}, \cite {smajdor1996});
we will need them in Section~\ref {acocs} below to obtain a characterization of Choquet simplexes and a partial converse to
Theorem~\ref {everyls}.

The \emph {suspension} of a convex set $K \subset \mathcal S$ in a semilinear space $\mathcal S$ is the semilinear space
$K_{sus} = \{ (\lambda, \lambda x) \mid \lambda > 0, x \in K\} \subset \mathbb R_+ \times K$
with respect to the semilinear operations inherited from $\mathbb R \times K$.  If $\mathcal S$ is a linear space
then $K_{sus}$ is a cone in the linear space $\mathbb R \times \mathcal S$.
If $K$ is equipped with a topology then there is a natural topology on $K_{sus}$ induced from $\mathbb R \times K$.
If $\mathcal S$ is a semilinear topological space and $\tint K \neq \emptyset$
then the natural topology on $K_{sus}$ has a base
\begin {multline*}
\mathcal B_{(\mu, \mu z)} =
\left\{ (\mu - \varepsilon, \mu + \varepsilon) \times (V \cap K) \mid 0 < \varepsilon < \mu,\right. 
\\
\left.
\text {$V \subset \mathcal S$ is a balanced open neighbourhood of the point $z$}
\right\}
\end {multline*}
of open neighbourhoods for every point $(\mu, z) \in K_{sus}$. The base $\mathcal B_{(\mu, \mu z)}$ has a sub-base of balanced open neighbourhoods for
every point
$$
(\mu, z) \in \tint K_{sus} = (\tint K)_{sus},
$$
so $K_{sus}$ is also a semilinear topological space.

Let $\mathcal S$ and $\mathcal V$ be semilinear spaces, and suppose that $\mathcal V$ is equipped with an order $\leqslant$ consistent with the semilinear
structure of $\mathcal V$.
Suppose that $T : K \to \mathcal V$
is a map defined on a convex set $K \subset \mathcal S$.
The map $T$ is said to be \emph {convex} if
$$
(1 - \theta) T (x) + \theta\, T (y) \leqslant T ((1 - \theta) x + \theta y)
$$
for all $x, y \in K$ and
$0 < \theta < 1$.
$T$ is \emph {affine} if
$$
(1 - \theta) T (x) + \theta\, T (y) = T ((1 - \theta) x + \theta y)
$$
for all $x, y \in K$ and $0 < \theta < 1$.
These properties are related to superlinearity and linearity
in the following way (which is well known but somewhat implicit in the literature; see, e. g., \cite [Chapter 11] {convexityanal}).
For every map $T$ on $K$ we define the \emph {natural extension $T_{sus}$ of $T$ onto $K_{sus}$} by
$T ((\lambda, \lambda x)) = \lambda T (x)$ for $x \in K$.
It is evident that the natural extension $T_{sus}$ is always positive homogeneous and
the map $T$ is easily recovered from its
natural extension $T_{sus}$ by its restriction $T (x) = T_{sus} ((1, x))$ on the set $\{1\} \times K$.
It is easy to see that if $\mathcal S$ is a semilinear topological space, $\mathcal W$ is a semilinear space equipped with a topology
in which its semilinear operations are continuous, and $T$ takes values in some subspace $\mathcal V \subset \mathcal C (\mathcal W)$ of closed sets,
then $T$ is upper semicontinuous if and only if $T_{sus}$ is upper semicontinuous.
\begin {proposition}
\label {coaff}
Suppose that $\mathcal S$ and $\mathcal V$ are semilinear spaces, $K \subset \mathcal S$ is a convex set and
$T : K \to \mathcal V$ is a map.
\begin {enumerate}
\item
$T$ is affine if and only if $T_{sus}$ is linear (equivalently, additve).
\item
$T$ is convex if and only if $T_{sus}$ is superlinear (equivalently, superadditive).
\end {enumerate}
\end {proposition}
Indeed, let $(\lambda, \lambda x)$ and $(\mu, \mu y)$ be arbitrary points in $K_{sus}$.
If $T$ is convex then
\begin {multline}
\label {conveq}
T_{sus} ((\lambda, \lambda x)) + T_{sus} ((\mu, \mu y)) =
\lambda T (x) + \mu T (y) =
\\
(\lambda + \mu) \left( \frac \lambda {\lambda + \mu} T (x) + \frac \mu {\lambda + \mu} T (y) \right) \leqslant
\\
(\lambda + \mu) T \left(\frac \lambda {\lambda + \mu} x + \frac \mu {\lambda + \mu} y\right) =
\\
T_{sus} \left(\left(\lambda + \mu, (\lambda + \mu) \left(\frac \lambda {\lambda + \mu} x + \frac \mu {\lambda + \mu} y\right)\right)\right) =
\\
T_{sus} \left((\lambda, \lambda x) + (\mu, \mu y)\right),
\end {multline}
so $T_{sus}$ is additive and hence linear.  If $T$ is affine then we can replace the inequality sign $\leqslant$ by the equality sign $=$
in \eqref {conveq} to see that
$T_{sus}$ is superadditive and hence superlinear.
On the other hand, let $x, y \in K$ be arbitrary points and $0 < \theta < 1$ be an arbitrary number.
If $T_{sus}$ is superadditive then in the same way
\begin {multline}
\label {supereq}
(1 - \theta) T (x) + \theta T (y) =
\\
T_{sus} (((1 - \theta), (1 - \theta) x)) + T_{sus} ((\theta, \theta y)) \leqslant
\\
T_{sus} (((1 - \theta), (1 - \theta) x) + (\theta, \theta y)) =
\\
T_{sus} ((1, (1 - \theta) x + \theta y)) =
T \left( (1 - \theta) x + \theta y \right),
\end {multline}
so $T$ is convex.  If $T_{sus}$ is additive then  we can replace the inequality sign $\leqslant$ with the equality sign $=$
in \eqref {supereq} to see that $T$ is affine.

Proposition \ref {coaff} allows us to transfer directly the results about linear submaps of superlinear maps and their linear selections
to the case of maps that are convex and affine. 
The following is an analogue and consequence of Theorem \ref {linselg}.
\begin {theorem}
\label {linselga}
Suppose that $K \subset \mathcal S$ is a convex set in a semilinear space $\mathcal S$, $K_{sus}$ has the Riesz decomposition property,
$\mathcal V$ is a semilinear space
equipped with a partial order $\leqslant$ compatible with the semilinear structure, and $\mathcal V$ is order regular
with respect to $\leqslant$.  Then any convex map $T : K \to \mathcal V$ has an affine submap $S \leqslant T$
given by the formula
\begin {equation*}
\label {thesubmapa}
S (x) = \inf_{x = \sum_j \theta_j x_j} \sum_j \theta_j T (x_j), \quad x \in K,
\end {equation*}
where the infinum is taken over all finite convex combinations
$$
x = \sum_j \theta_j x_j, \quad 0 \leqslant \theta_j \leqslant 1, \quad \sum_j \theta_j = 1.
$$
Map $S$ is the greatest affine submap of $T$,
i. e. if $R \leqslant T$ is another affine submap of $T$ then $R \leqslant S$.
\end {theorem}

The following is an analogue and consequence of Theorem \ref {linselnd}.
\begin {theorem}
\label {linselnda}
Suppose that $\mathcal S$ and $\mathcal W$ are semilinear spaces, $K \subset \mathcal S$ is a convex set,
$\mathcal V$ is a subspace of the semilinear space $\mathcal C (\mathcal W)$,
$\mathcal D \subset \mathcal W'$ is an exhaustive set of linear functionals that are defining for and consistent with $\mathcal V$,
and either $\mathcal D$ is finite or $\mathcal V$ has compact type.
Then for every affine map $A : K \to \mathcal V$, $x \in K$ and $y \in A (x)$ there exists an affine selection $a$ of $A$
such that $a (x) = y$.  Suppose, in addition, that $\mathcal S$ and $\mathcal W$ are equipped with topologies so that all sets in $\mathcal V$ are closed,
$\mathcal W$ is a regular topological space, $\mathcal S$ is a semilinear topological space,
all functionals in $\mathcal D$ are continuous and $A$ is upper semicontinuous.  Then
there exists a continuous affine selection $a$ of $A$ satisfying $a (x) = y$.
\end {theorem}

The following result is an analogue and consequence of Theorem \ref {lsvmch}.
\begin {theorem}
\label {asvmch}
Suppose that $\mathcal S$ and $\mathcal W$ are semilinear spaces, $K \subset \mathcal S$ is a convex set,
$\mathcal V$ is a subspace of the semilinear space $\mathcal C (\mathcal W)$,
$\mathcal D \subset \mathcal W'$ is an exhaustive set of linear functionals that are defining for and consistent with $\mathcal V$,
and either $\mathcal D$ is finite or $\mathcal V$ has compact type.  Let $T : K \to \mathcal V$ be a convex map.
The following conditions are equivalent.
\begin {enumerate}
\item
$T$ is affine.
\item
For any $x \in \mathcal S$ and $y \in T (x)$ there exists an affine selection $S$ of $T$ satisfying $S (x) = y$.
\end {enumerate}
\end {theorem}
Theorem \ref {asvmch} extends \cite [Theorem 2] {gorokhovik2008} that was obtained for the finite-dimensional case by a different method
involving parametrization of the set of all affine selections via its extreme or exposed points (cf. Theorems \ref {riesznecc} and \ref {csselect}
in Section \ref {acocs} below).
Our method seems to be stronger since it does not depend on the existence of convex or exposed points for values $\mathcal V$ of the map in question.

The following is an analogue and consequence of Theorem \ref {everyls}.
\begin {theorem}
\label {everylsa}
Suppose that $\mathcal S$ and $\mathcal W$ are semilinear spaces, $K \subset \mathcal S$ is a convex set,
$K_{sus}$ has the Riesz decomposition property,
$\mathcal V$ is a subspace of the semilinear space $\mathcal C (\mathcal W)$ that has compact type and
$\mathcal D \subset \mathcal W'$ is an exhaustive set of linear functionals that are defining for and consistent with $\mathcal V$.
Let $T : K \to \mathcal V$ be a convex map.  Then $T$ admits an affine selection.
Additionally, fix some points $x \in K$ and $y \in T (x)$.
Then the following conditions are equivalent.
\begin {enumerate}
\item
There exists an affine selection $a$ of $T$ such that $a (x) = y$.
\item
$y \in \sum_{j = 1}^N \theta_j T (x_j)$ for any $N \in \mathbb N$,
$\{\theta_j\}_{j = 1}^N \subset [0, 1]$ and a sequence $\{x_j\}_{j = 1}^N \subset K$ such that $x = \sum_{j = 1}^N \theta_j x_j$ and $\sum_{j = 1}^N \theta_j = 1$.
\end {enumerate}
If, additionally, $\mathcal S$ is a semilinear topological space,  $\mathcal W$ is a locally convex Hausdorff linear topological space,
$\mathcal V$ consists of all nonempty compact convex sets in $\mathcal W$,
all functionals in $\mathcal D$
are continuous and $T$ is upper semicontinuous, then the selection in condition 1 can be taken continuous.
\end {theorem}

\section {Characterization of Choquet simplexes}

\label {acocs}

In this section we will demonstrate that under suitable restrictions on the semilinear spaces
the existence of linear selections for arbitrary superlinear maps
implies that the semilinear space the maps are acting on is a lattice and thus has the Riesz decomposition property,
thus establishing a partial converse to Theorem~\ref {everyls}.
This result relies heavily on the Choquet theory; specifically, with the help of the Choquet-Meyer theorem
we show that under some appropriate restrictions
on the spaces every convex map $T : K \to \mathcal V \subset \mathcal C (\mathcal W)$ admits an affine selection
if and only if $K$ is a Choquet simplex.

First, we need to introduce a few relevant notions.
The suspension construction introduced in Section~\ref {ascm} above
relates a convex semilinear space to every convex set of a semilinear space in a particular way.
It is well known that many (but not all) semilinear spaces can be represented as suspensions of convex sets.
Let $\mathcal S$ be a semilinear space.
We say that a convex set $B \subset \mathcal S$ is a \emph {base of $\mathcal S$} if for every $y \in \mathcal S$ which is not an additive identity
there exist unique
$\lambda > 0$ and $x \in B$ such that $y = \lambda x$.
It is easy to see that 
for every convex set $K \subset \mathcal V$ in a semilinear space $\mathcal V$
the set $\{1\} \times K$ is a base for the suspension $K_{sus}$.

It is said that
convex sets $B \subset \mathcal S$ and $K \subset \mathcal V$ in semilinear spaces $\mathcal S$ and $\mathcal V$ are \emph {affinely isomorphic}
if there exists an affine bijection $T : B \to K$ onto $K$; the inverse map $T^{-1}$ is then also an affine bijection $T : K \to B$ onto $B$.
For example, the set $\{1\} \times K \subset K_{sus}$ is affinely isomorphic to $K$.
It is easy to see that (see \cite [Proposition 11.2] {convexityanal})
an affine isomorphism $R$
between bases $B \subset \mathcal S$ and $K \subset \mathcal V$ of respective semilinear spaces $\mathcal S$ and $\mathcal V$
can be extended to a linear isomorphism $T$ between $\mathcal S$ and $\mathcal V$ in the following way: for $y \in \mathcal S$
we take $\lambda > 0$ and $x \in B$ such that $y = \lambda x$ and set $T (y) = \lambda R (x)$; moreover, if semilinear spaces $\mathcal S$
and $\mathcal V$ are equipped with topologies such that their respective semilinear operations are
continuous and $R$ is continuous then $T$ is also continuous;
thus if $R$ is a homeomorphism then $T$ is an affine homeomorphism between $\mathcal S$ and $\mathcal V$.
This implies that if a semilinear space $\mathcal V$ has a base $K$ then $\mathcal V$ is affinely isomorphic to $K_{sus}$, and if $\mathcal V$ is a semilinear
topological space then this isomorphism is also a homeomorhism.  In other words, any semilinear space $\mathcal V$ with a base $K$
is just the suspension $K_{sus}$ of the base $K$ up to an affine isomorphism.
Not every cone has a cone-base, but it is a fairly common property.
Let us mention a few examples.  The cone of finite nonnegative measures $\mathcal M_+ (\Omega)$ on a measurable space $\Omega$ has the set of
probability measures $\mathcal M_{+, 1}$ as a base.
Likewise, we can restrict measures in this example to those that are absolutely continuous with respect
to a fixed measure $\mu$; this leads to the fact that the cone of all
nonnegative $\mu$-summable functions $\left[ \lclass {1} {\Omega, \mu} \right]_+$ has a base,
for example the set of measurable functions $f \geqslant 0$ such that
$\int f d\mu = 1$.  On the other hand, if a cone $C$ is not generating (i. e. if $C \cap -C$ has more than one point) then $C$ does not have a base.
The cone $K = \{ (x, y) \mid (y = 0 \wedge x \geqslant 0) \vee y > 0 \}$
corresponding to the lexicographic ordering of the plane $\mathbb R^2$ is generating, but it is easy to see that it does not have a
base.  For more information on cones that admit bases see, e.~g., \cite [1.9] {jameson1970}.

We now introduce the last set of relevant notions before stating the main result of this section; see, e.~g.,
\cite [Chapter 9--11] {convexityanal} for a detailed introduction into Choquet theory.
Let $K$ be a compact convex set in a locally convex Hausdorff linear topological space $\mathcal V$.  Then for every probability Baire measure $\mu$ on $K$
by \cite [Theorem 9.1] {convexityanal} there exists a unique point $x \in K$ such that for all continuous affine functionals $\ell$ on $K$
one has $\ell (x) = \int_K \ell (y) d\mu (y)$.  The point $x$ is called a \emph {barycenter} or \emph {the center of mass} of the measure $\mu$. 
If $\mathcal V$ is finite-dimensional then this is equivalent to
$x = \int_K y d\mu (y)$.
Let $\mathcal E (K)$ be the set of extreme points of $K$.
The famous Choquet theorem \cite [Theorem 10.7] {convexityanal} states that if $K$ is a \emph {metrizable} compact convex set in a locally convex
Hausdorff linear topological space then the set of extreme points $\mathcal E (K)$ is Baire measurable and
for any $x \in K$ there exists a probability Baire measure $\mu$ on $K$ such that $\mu (\mathcal E (K)) = 1$
and $x$ is the barycenter of $\mu$.
A convex set $K \subset \mathcal S$ in a semilinear space $\mathcal S$ is called an
\emph {algebraic simplex} or just \emph {simplex}
if $K_{sus}$ is a lattice.
If, additionally, $\mathcal S$ is a locally convex Hausdorff linear topological space and $K$ is a compact set then
$K$ is called a \emph {Choquet simplex}.
The well-known corollary to the Choquet-Meyer theorem (that we will introduce in a moment)
\cite [Theorem 11.12] {convexityanal} states the following: if $K$ is a metrizable compact set in a locally convex Hausdorff
linear topological space then $K$ is a Choquet simplex if and only if for any $x \in K$ there exists a \emph {unique} probability
measure $\mu_x$ on $K$ with barycenter $x$ such that $\mu_x (\mathcal E (K)) = 1$.

But what if $K$ is not metrizable?  This case is also fairly well understood.
Without the assumption of metrizability the
set of extreme points $\mathcal E (K)$ is not necessarily Baire measurable, so we need more subtle terms to express
the condition that a measure $\mu$ is located at $\mathcal E (K)$.
We denote by $C (B)$ the space of all continuous (real-valued) functions on a set $B \subset K$,
by $\mathcal A (K) \subset C (K)$ the space of all continuous affine functions on $K$ and by $\mathscr C (K) \subset C (K)$
the set of all continuous convex functions on $K$.
For a function $f \in \mathscr C (K)$ its \emph {concave envelope} $\hat f$ is defined by
$$
\hat f (x) = \inf \left\{ g (x) \mid -g \in \mathscr C (K), g \geqslant f \right\}, x \in K.
$$
The set $\mathscr C (K)$ can be replaced by the set of continuous affine functions $\mathcal A (K)$ in this definition.
\begin {theoremchoquetmeyer}[{\cite [Theorem 11.5]{convexityanal}\footnote {
We state the two last equivalent conditions for the case of metrizable $K$ only to avoid discussing maximal measures in the Choquet order
since we do need it in the present work.
}}]
Let $K$ be a compact convex set in a locally convex Hausdorff linear topological space.
Recall that $\mathcal M_{+, 1} (K)$ is the space of Baire probability measures on $K$.
The following conditions are equivalent.
\begin {enumerate}
\item
$K$ is a simplex.
\item
For each $f \in \mathscr C (K)$ its concave envelope $\hat f$ is affine.
\end {enumerate}
Additionally, if $K$ is metrizable then there are two more equivalent conditions.
\begin {enumerate}
\item[(3)]
For any $f \in \mathscr C (K)$ and $\mu \in \mathcal M_{+, 1} (K)$ having $x \in K$ as the barycenter it is true that $\int_K f d\mu = \hat f (x)$.
\item[(4)]
For each $x \in K$ there exists a unique Baire probability measure $\mu \in \mathcal M_{+, 1} (K)$ having $x$ as its barycenter.
\end {enumerate}
\end {theoremchoquetmeyer}
We are now ready to state the main result.

\begin {theorem}
\label {latnecc}
Suppose that $\mathcal S$ is a semilinear topological space with a base $K$ which is a compact convex set in a locally convex
Hausdorff linear topological space.  Let $\mathcal W$ be any nontrivial locally convex Hausdorff linear topological space and
$\mathcal V (\mathcal W) \subset \mathcal C (\mathcal W)$ be the space of all nonempty compact convex sets in $\mathcal W$.
The following conditions are equivalent.
\begin {enumerate}
\item
$K$ is a Choquet simplex.
\item
$\mathcal S$ is a lattice.
\item
$\mathcal S$ has the Riesz decomposition property.
\item
For any superlinear
map $T : \mathcal S \to \mathcal V (\mathcal W)$ there exists a linear selection $S$ of $T$.
\item
For any convex
map $T : K \to \mathcal V (\mathcal W)$ there exists an affine selection $S$ of $T$.
\item
Every superlinear map $T : \mathcal S \to \mathcal V (\mathbb R)$ has a linear selection $S$.
\item
Every convex map $T : K \to \mathcal V (\mathbb R)$ has an affine selection $S$.
\end {enumerate}
\end {theorem}
Note that equivalence of conditions 1--3 of Theorem \ref {latnecc} (and also condition 8 of Theorem \ref {riesznecc} below)
is rather well known; see, e. g., \cite [\S 3] {shashkin1973};
however, we establish it using nothing more than the Choquet-Meyer theorem and the results of this paper.
All sequential transitions $1 \Rightarrow 2 \Rightarrow 3 \Rightarrow 4 \Leftrightarrow 5 \Rightarrow 6 \Leftrightarrow 7$ are either trivial or have already
been covered: $1 \Rightarrow 2$ is the definition of a Choquet simplex which is that $K_{sus}$ is a lattice together with the fact that
$K_{sus}$ is affinely isomorphic to $\mathcal S$ because $K$ is a base of $\mathcal S$,
$2 \Rightarrow 3$ is a simple property mentioned in the introduction,
$3 \Rightarrow 4$ is a consequence of Theorem \ref {everyls}, 
$4 \Leftrightarrow 5$ and $6 \Leftrightarrow 7$ follow from Proposition \ref {coaff}, and implication $5 \Rightarrow 7$ (or $4 \Rightarrow 6$) is trivial
as we can immerse $\mathcal V (\mathbb R)$ into $\mathcal V (\mathcal W)$ by mapping $\mathbb R$ onto a line in $\mathcal W$.
The remaining implication $7 \Rightarrow 1$ is a simple consequence of the Choquet-Meyer theorem.  Indeed, suppose that $f \in \mathscr C (K)$;
we need to show that its concave envelope $\hat f$ is affine.
Define a set-valued map $T : K \to \mathcal V (\mathbb R)$ by $T (x) = \{y \mid f (x) \leqslant y \leqslant \hat f (x) \}$, $x \in K$.  It is easy to see that
$T$ is a convex map because it is
an intersection of two convex maps $\{y \geqslant f (x)\}$ and $\{y \leqslant \hat f (x)\}$.
Condition 7 implies that there exists an affine selection $S$ of $T$.  Then $S \geqslant \hat f$ by the definition of concave envelope and
$S \leqslant \hat f$ because $S$ is a selection of $T$, thus $S = \hat f$ and $\hat f$ is affine.  The proof of Theorem \ref {latnecc} is complete.

If $K$ is a Choquet simplex and $x \in K$, let us denote by $\mu_x$ the unique Baire probability measure on $K$ such that $\mu_x (\mathcal E (K)) = 1$
and $\mu_x$ has barycenter $x$.  When is the map $x \mapsto \mu_x$ continuous?
The answer is provided by the following particular case of a well known result.
\begin {theorem}
[{\cite [Theorem 11.13] {convexityanal}\footnote {Again, we avoid discussing maximal measures in the Choquet order by making use of the fact that
the set of maximal measures
coincides with the family of measures $\{\mu_x\}_{x \in K}$ if the compact convex set $K$ is metrizable.}}]
\label {bairesimp}
Let $K$ be a metrizable Choquet simplex in a locally convex Hausdorff linear topological space.  The following conditions are equivalent.
\begin {enumerate}
\item
The set $\mathcal E (K)$ of extreme points of $K$ is closed.
\item
The map $x \mapsto \mu_x$ is continuous in the weak* topology\linebreak{} $\sigma (\mathcal M (K), C (K))$ of $\mathcal M (K)$.
\item
Family $\{\mu_x\}_{x \in K}$ is closed in the weak* topology of $\mathcal M (K)$.
\item
The convex envelope $\hat f$ is continuous for all $f \in \mathscr C (A)$.
\end {enumerate}
\end {theorem}
If any (and hence all) of the conditions of Theorem \ref {bairesimp} is satisfied, $K$ is called a \emph {Bauer simplex}.
In this case we can refine the characterization provided by Theorem \ref {latnecc} 
as follows.
\begin {theorem}
\label {riesznecc}
Suppose that $\mathcal S$ is a semilinear topological space with a base $K$ which is a compact metrizable convex set in a locally convex
Hausdorff linear topological space and the set $\mathcal E (K)$ of extreme points of $K$ is closed.
Let $\mathcal W$ be any nontrivial locally convex Hausdorff linear topological space and
$\mathcal V (\mathcal W) \subset \mathcal C (\mathcal W)$ be the space of all nonempty compact convex sets in $\mathcal W$.
Then the following conditions are equivalent.
\begin {enumerate}
\item
$K$ is a Choquet (and thus Bauer) simplex.
\item
$\mathcal S$ is a lattice.
\item
$\mathcal S$ has the Riesz decomposition property.
\item
For any upper semicontinuous superlinear
map $T : \mathcal S \to \mathcal V (\mathcal W)$ there exists a continuous linear selection $S$ of $T$.
\item
For any upper semicontinuous convex
map $T : K \to \mathcal V (\mathcal W)$ there exists a continuous affine selection $S$ of $T$.
\item
Any upper semicontinuous superlinear map $T : \mathcal S \to \mathcal V (\mathbb R)$ has a continuous linear selection $S$.
\item
Any upper semicontinuous convex map $T : K \to \mathcal V (\mathbb R)$ has a continuous affine selection $S$.
\item
$\mathcal A (K) = C (\mathcal E (K))$ in the sense that the restriction of all functions from $\mathcal A (K)$
onto the set $\mathcal E (K)$ of extreme points of $K$ coincides with
$C (\mathcal E (K))$.
\end {enumerate}
\end {theorem}
All sequential implications $1 \Rightarrow 2 \Rightarrow 3 \Rightarrow 4 \Leftrightarrow 5 \Rightarrow 6 \Leftrightarrow 7$ are
established as in Theorem \ref {latnecc}, making use of the topological clauses of Theorem~\ref {everyls}.
Implication $8 \Rightarrow 1$ follows from the particular case of the Choquet-Meyer theorem (see also \cite [Theorem 11.5] {convexityanal}).
If condition 8 is satisfied
and there are two Baire measures $\mu_1$ and $\mu_2$ concentrated on $\mathcal E (K)$ having $x \in K$ as the barycenter then their difference
$\mu = \mu_1 - \mu_2$ satisfies $\int_K \ell (y) d\mu (y) = 0$ for all $\ell \in \mathcal A (K) = C (\mathcal E (K))$ which is a predual
for the space of Borel measures $\mathcal M (\mathcal E (K))$ on $\mathcal E (K)$
(note that the $\sigma$-algebras of Borel and Baire sets on $K$ coincide because $K$ is metrizable)
having finite variation, and therefore $\mu = 0$, which proves the uniqueness condition in the Choquet-Meyer theorem.
The remaining implication $7 \Rightarrow 8$ is not quite as trivial.  Let us prove it.

First observe that a map $M : K \to \mathcal M (\mathcal E (K))$ defined by
\begin {multline*}
M (x) = \left\{ \mu \mid \text {$\mu$ is a Baire measure on $K$,} \right.
\\
\left.\text {$\mu (\mathcal E (K)) = 1$ and $x$ is the barycenter of $\mu$}\right\}
\end {multline*}
is convex.  Indeed, for any $0 < \theta < 1$, $x, y \in K$, $\mu \in M (x)$, $\nu \in M (y)$ and a continuous affine functional $\ell$ on $K$
we have
\begin {multline*}
\ell ((1 - \theta) x + \theta y) = (1 - \theta) \ell (x) + \theta \ell (y) =
\\
(1 - \theta) \int_K \ell (z) d\mu (z) + \theta \int_K \ell (z) d\nu (z) =
\\
\int_K \ell (z) d\left[(1 - \theta) \mu + \theta \nu\right] (z),
\end {multline*}
so $(1 - \theta) \mu + \theta \nu \in M ((1 - \theta) x + \theta y)$ and thus
$$
(1 - \theta) M (x) + \theta M (y) \subset M ((1 - \theta) x + \theta y).
$$

Now fix $f \in C (\mathcal E (K))$; we need to find a function $S$ in $\mathcal A (K)$ that coincides with $f$ on $\mathcal E (K)$.
Define a map $T : K \to \mathcal V (\mathbb R)$ by
$$
T (x) = (M \circ f) (x) = \left\{ \int_{\mathcal E (K)} f (z) d\mu (z) \mid \mu \in M (x) \right\}
$$ 
for $x \in K$, which is just composition of $f$ and $M$.  It is easy to derive from convexity of $M$ that $T$ is a convex map.
Let us verify that the graph of $T$ is closed.  Suppose that $\{x_n\} \subset K$ is an arbitrary sequence
converging to some $x \in K$ and a sequence $y_n \in T (x_n)$ converges to some $y \in \mathbb R$; we need to show that $y \in T (x)$.
Note that $y_n = \int_{\mathcal E (K)} f (z) d\mu_n (z)$ for some measures $\mu_n \in M (x_n)$.
Observe that since $\mathcal E (K)$ is closed in $K$, $\mathcal E (K)$ is a compact set, and thus $C (\mathcal E (K))$ is a Banach space and its dual is
$\mathcal M (\mathcal E (K))$.
The sequence $\{\mu_n\}$ is bounded in
$\mathcal M (\mathcal E (K))$, so we can pass to a subsequence and assume that $\mu_n$ converges in the weak$*$ topology
to some $\mu \in \mathcal M (\mathcal E (K))$.  As for all $\ell \in \mathcal A (K)$ (and, in particular, for $\ell = 1$) we have
$\ell (x_n) = \int_{\mathcal E (K)} \ell (z) d\mu_n (z)$, passage to the limit shows that $\mu \in M (x)$.
Surely $y_n$ converges to $y = \int_{\mathcal E (K)} f (z) d\mu (z) \in \Phi (x)$, so the graph of $T$ is indeed closed and thus
$T$ is upper semicontinuous by Proposition \ref {grclos}.

If condition 7 of Theorem \ref {riesznecc} is satisfied, then $T$ has a continuous affine selection $S$.
Note that $M (x) = \{\delta_x\}$ for $x \in \mathcal E (K)$ where $\delta_x$ is a unit point mass concentrated at $x$;
this follows at once from the Bauer theorem (see, e. g., \cite [Theorem 9.3] {convexityanal}).  Therefore $T (x) = \{f (x)\}$ and $S (x) = f (x)$
for $x \in \mathcal E (K)$, i. e. $S$ coincides with $f$ on the set of extreme points $\mathcal E (K)$.
We have just verified implication $7 \Rightarrow 8$.  The proof of Theorem~\ref {riesznecc} is complete.

We mention that in the finite-dimensional case, i. e. when $\mathcal S \subset \mathbb R^{n + 1}$ or $K \subset \mathbb R^n$,
Theorems~\ref {latnecc} and \ref {riesznecc} can be obtained with much simpler means by making use of cone-bases\footnote {Plural of ``cone-basis''.}
and simple arithmetic
instead of tomographical coordinates and the Riesz decomposition property;
see the end of Section~\ref {conebases}.

\section {Parametrization of the set of selections}

\label {mosom}
In this section we explore how the results of the previous Section~\ref {acocs}
can be used in certain cases to naturally parametrize the set of all affine selections of a given convex map by their values
on the extreme points of the simplex, which leads to some extensions of Theorems~\ref {everylsa}~and~\ref {everyls}.
These extensions, however, are applicable essentially only to the case of upper semicontinuous maps acting on the semilinear space of Borel measures
$\mathcal S = \mathcal M_{+}(K)$ with the weak$*$ topology on a metrizable compact set $K$.

Condition 8 of Theorem~\ref {riesznecc} in Section~\ref{acocs} suggests a natural parametri\-zation for the set of continuous affine selections of an
upper semicontinuous convex set-valued map into $\mathcal V (\mathbb R)$ by its values on $\mathcal E (K)$.
More generally, an easy application of Choquet theory yields a parametrization for all
continuous affine selections of an upper semicontinuous convex map acting into a Banach space $X$ by its values on the set $\mathcal E (K)$.
In order to do this we employ the Bochner integral instead of the scalar one; see, e. g., \cite [\S4, \S5, Chapter V] {ioshidaen}.
Note that if $A \subset K$ is a Baire measurable subset of $K$ and $f : A \to X$
is a continuous function then $f (A)$ is separable in $X$ because $A$ is separable
(the closure of $A$ in $K$ is a metrizable compact set that admits $\varepsilon$-nets for any $\varepsilon > 0$); $f$ is weakly measurable because
$f$ is weakly continuous (i. e. with respect to the weak topology $\sigma (X, X^*)$ in $X$),
and thus by the Pettis theorem $f$ is strongly measurable.  Therefore the Bochner integral $\int_A f d\mu$ is well-defined
for any $\sigma$-finite Baire measure $\mu$ on $K$ if $\int_A \|f\|_X d\mu < \infty$; in particular, this is the case if $\mu (K) = 1$ and
$A$ is closed and thus $f$ is bounded on~$A$.
We denote by $C (A \to X)$ the set of $X$-valued continuous functions on~$A$.
\begin {theorem}
\label {csselect}
Suppose that $K$ is a metrizable compact set in a locally convex Hausdorff linear topological space, $K$ is a Bauer simplex
(i. e. the set $\mathcal E (K)$ of extreme points of $K$ is closed in $K$),
$\mu_x$ is the measure on $\mathcal E (K)$ with barycenter $x$ for $x \in K$,
$X$ is a Banach space,
$T : K \to 2^X \setminus \{\emptyset\}$
is an upper semicontinuous convex set-valued map and values of $T$ are closed in $X$.
Then for every $f \in C (\mathcal E (K) \to X)$ satisfying $f (x) \in T (x)$ for all $x \in \mathcal E (K)$
function $S_f (x) = \int_{\mathcal E (K)} f (y) d\mu_x (y)$, $x \in K$,
is a continuous affine selection $S$ of $T$.  Conversely, any continuous affine selection
$S$ of $T$ has the representation $S = S_f$ where $f$ is the restriction of $S$ onto $\mathcal E (K)$.
\end {theorem}

Indeed, by the proof of Theorem~\ref {riesznecc} the map $x \mapsto \mu_x$ is affine, and therefore $S_f$ is also affine.
Let us verify that $S_f$ is continuous.  Suppose that $x_n \in K$ is an arbitrary sequence converging to some $x \in K$ in $K$;
we need to establish that $S_f (x_n) \to S_f (x)$.  Let $\varepsilon_k > 0$ be a sequence of numbers converging to $0$.
Since $f \left(\mathcal E (K)\right)$ is a compact set, for every~$k$ there exists a finite cover $\mathcal B_k = \{B_j^k\}_{j = 1}^{m_k}$
of $f \left(\mathcal E (K)\right)$ by balls $B_j^k$ with respective centers $z_j^k$ having radius $\varepsilon_k$.
Let $\{\varphi_j^k\}_{j = 1}^{m_k}$ be a continuous partition of unity on $f \left(\mathcal E (K)\right)$ subordinate to $\mathcal B_k$.
We introduce approximating functions $f_k (y) = \sum_{j = 1}^{m_k} \varphi_j^k (f (y)) z_j^k$ for $y \in \mathcal E (K)$; let us verify that
the sequence $f_k$ converges to $f$ uniformly in $\mathcal E (K)$.
Suppose that $y \in K$. Then $f (y) \in B_j^k$ for some $j$,
condition $\varphi_l^k (f (y)) > 0$ implies that balls $B_j^k$ and $B_l^k$ intersect, and thus $\|z_l^k - z_j^k\|_X \leqslant 2 \varepsilon_k$
because both of these balls have radius $\varepsilon_k$.  Therefore
\begin {multline}
\label {fkdiff}
\|f_k (y) - f (y)\|_X =
\left\|\sum_{\{l \mid \varphi_l^k (f (y)) > 0\}} \varphi_l^k (f (y)) z_l^k - z_j^k + z_j^k - f (y)\right\|_X \leqslant
\\
\left\|\sum_{\{l \mid \varphi_l^k (f (y)) > 0\}} \varphi_l^k (f (y)) z_l^k - z_j^k\right\|_X + \left\|z_j^k - f (y)\right\|_X \leqslant
\\
2 \varepsilon_k + \varepsilon_k = 3 \varepsilon_k.
$$
\end {multline}
Since $K$ is a Bauer simplex, the map $x \mapsto \mu_x$ is continuous with respect to the weak* topology of $\mathcal M (\mathcal E (K))$.
This means that for any $k$ there exists some $N_k$ such that for all $n \geqslant N_k$ we have
\begin {multline}
\label {mukdiff}
\left| \int_{\mathcal E (K)} \varphi_j^k (f (y)) d\mu_{x_n} (y) -  \int_{\mathcal E (K)} \varphi_j^k (f (y)) d\mu_{x} (y)\right| \leqslant
\\
\frac {\varepsilon_k} {m_k \left( 1 + \bigvee_{j = 1}^{m_k} \|z_j^k\|_X \right)}
\end {multline}
for all $j$, $1 \leqslant j \leqslant m_k$.
Combining \eqref {mukdiff} and \eqref {fkdiff} together and taking advantage of the properties of the partition of unity, we easily arrive~at
\begin {multline*}
\|S_f (x_n) - S_f (x)\|_X \leqslant
\\
\|S_f (x_n) - S_{f_k} (x_n)\|_X + \|S_{f_k} (x_n) - S_{f_k} (x)\|_X + \|S_{f_k} (x) - S_f (x)\|_X \leqslant
\\
\left\|\int_{\mathcal E (K)} [f (y) - f_k (y)] d\mu_{x_n} (y)\right\|_X +
\\
\left\| \int_{\mathcal E (K)} \left[ \sum_{j = 1}^{m_k} \varphi_j^k (f (y)) z_j^k \right] d\mu_{x_n} (y) -
\int_{\mathcal E (K)} \left[ \sum_{j = 1}^{m_k} \varphi_j^k (f (y)) z_j^k \right] d\mu_{x} (y) \right\|_X +
\end {multline*}
\begin {multline*}
\left\| \int_{\mathcal E (K)} [f_k (y) - f (y)] d\mu_{x} (y) \right\|_X \leqslant
\\
\int_{\mathcal E (K)} \|f (y) - f_k (y)\|_X d\mu_{x_n} (y) +
\\
\sum_{j = 1}^{m_k} \left| \int_{\mathcal E (K)} \varphi_j^k (f (y)) d\mu_{x_n} (y) -  \int_{\mathcal E (K)} \varphi_j^k (f (y)) d\mu_{x} (y)\right| \|z_j^k\|_X +
\\
\int_{\mathcal E (K)} \|f_k (y) - f (y)\|_X d\mu_{x} (y) \leqslant
3 \varepsilon_k +  \varepsilon_k + 3 \varepsilon_k = 7 \varepsilon_k
\end {multline*}
for all $n \geqslant N_k$, which means that $S_f (x_n)$ converges to $S_f (x)$ as claimed.
It is easy to verify that $S_f$ is a selection of $T$.  Indeed, for any convex combination
$x = \sum_{j = 1}^N \theta_j x_j$, $0 \leqslant \theta_j \leqslant 1$, $\sum_{j = 1}^N \theta_j = 1$ of extreme points $x_j \in \mathcal E (K)$
we have $\mu_x = \sum_{j = 1}^N \theta_j \delta_{x_j}$ where $\delta_{x_j}$ is a point mass concentrated at $x_j$.
Therefore
$$
S_f (x) = \sum_{j = 1}^N \theta_j f (x_j) \in \sum_{j = 1}^N \theta_j T (x_j) \subset T \left(\sum_{j = 1}^N \theta_j x_j\right) = T (x)
$$
by the convexity of $T$, so $S_f \in T$ on the convex hull $\co \mathcal E (K)$ of $\mathcal E (K)$.
Since $T$ is upper semicontinuous the graph of $T$ is closed by Proposition \ref {grclos}, and therefore it contains the closure of the
restriction of the graph of $S_f$ onto $\co \mathcal E (K)$, which coincides with the entire graph of $S_f$ because $\co \mathcal E (K)$ is dense
in $K$ by the Krein-Milman theorem and $S_f$ is continuous; this proves that $S_f$ is a selection of $T$.
Finally, if two continuous affine maps $a, b : K \to X$ coincide on $\mathcal E (K)$ then they coincide on the convex hull $\co \mathcal E (K)$ of
$\mathcal E (K)$ which is dense in $K$, so $a = b$ everywhere on $K$ by continuity. Thus every continuous affine
selection $S$ of $T$ is uniquely determined by its values on $\mathcal E (K)$ (alternatively, we can invoke condition 8 of Theorem \ref {riesznecc}).
The proof of Theorem \ref {csselect} is complete.

Note that Theorem~\ref {csselect} extends Theorem~\ref {everylsa}: if $K$ is a Bauer simplex and $T$ is assumed to be an upper semicontinuous convex
map acting into arbitrary (not necessarily convex or compact) closed subsets of a Banach space $X$,
this is sufficient for the existence of affine selections of $T$ and even their complete characterization.
Taking Proposition~\ref {coaff} into account, this result also extends Theorem~\ref {everyls} in a similar way.

\begin {theorem}
\label {csselectsl}
Suppose that $\mathcal S$ is a semilinear space with a base $K$ which
is a metrizable compact set in a locally convex Hausdorff linear topological space,
$K$ is also a Bauer simplex,
$\mu_x$ is the measure on $\mathcal E (K)$ with barycenter $x$ for $x \in K$,
$X$ is a Banach space, $T : \mathcal S \to 2^X \setminus \{\emptyset\}$
is an upper semicontinuous superlinear set-valued map and values of $T$ are closed in $X$.
Then for every $f \in C (\mathcal E (K) \to X)$ satisfying $f (x) \in T (x)$ for all $x \in \mathcal E (K)$
function $S_f (\lambda x) = \lambda \int_{\mathcal E (K)} f (y) d\mu_x (y)$, $x \in K$,
$\lambda > 0$, is a continuous linear selection $S$ of $T$ satisfying $S (x) = f (x)$ for $x \in \mathcal E (K)$.
Conversely, any continuous linear selection
$S$ of $T$ has the form $S = S_f$ where $f$ is the restriction of $S$ onto~$\mathcal E (K)$.
\end {theorem}
Theorem~\ref {csselectsl} follows at once from Proposition~\ref {coaff} and Theorem~\ref {csselect} applied to the restriction of $T$ onto $K$.
These results, however, extend a very special case of the general Theorems~\ref {everylsa} and \ref {everyls}, since the semilinear space
$\mathcal S$ under the conditions of
Theorem~\ref {csselectsl} is isomorphic to $\mathcal M_+(\mathcal E (K))$ with the weak* topology
and the compact set $K$ under the conditions of Theorem~\ref {csselect}
is isomorphic to $\mathcal M_{+, 1} (\mathcal E (K))$ with the weak* topology.  These observations lead to the following generalization of
Theorem~\ref {csselectsl}.
\begin {theorem}
\label {csselectm}
Suppose that $K$ is a metrizable compact topological space,
$X$ is a Banach space, $T : \mathcal M_{+} (K) \to 2^X \setminus \{\emptyset\}$
is a superlinear set-valued map, values of $T$ are closed in $X$ and $T$ is upper semicontinuous
with respect to the weak$*$ topology of $\mathcal M_+ (K)$.
Then for every $f \in C (K \to X)$ satisfying $f (x) \in T (\delta_x)$ for all $x \in K$
function $S_f (\mu) = \int_{K} f (y) d\mu (y)$, $\mu \in \mathcal M_+ (K)$
is a unique continuous linear selection of $T$ satisfying $S(\delta_x) = f (x)$ for $x \in K$.  Conversely, any continuous linear selection
$S$ of $T$ has the form $S = S_f$ where $f (x) = S (\delta_x)$ for~$x \in K$.
\end {theorem}
Indeed, let $\mathcal K = \mathcal M_{+, 1} (K)$ with the weak$*$ topology induced from $\mathcal M_+ (K)$.
Thus $\mathcal K$ is a metrizable compact set since $K$ and consequently $C (K)$ are separable.
The space $\mathcal M_{+} (K)$ is a lattice: for every $\mu_1, \mu_2 \in M_{+}(K)$ both of these measures are absolutely continuous with
respect to $\nu = \mu_1 + \mu_2$, thus $d\mu_1 = \sigma_1d\nu$, $d\mu_2 = \sigma_2d\nu$ with some $\sigma_1, \sigma_2 \in \lclass {1} {\nu}$
and $d (\mu_1\vee \mu_2) = (\sigma_1 \vee \sigma_2) d\nu$ is the least upper bound for $\mu_1$ and $\mu_2$.  Therefore $\mathcal K$ is a simplex.
It is easy to see that the map $x \mapsto \delta_x$ establishes a homeomorphism $\rho$ between $K$ and $\mathcal E (\mathcal K)$,
thus $\mathcal K$ is a Bauer simplex.  Moreover, the measure on $\mathcal E (\mathcal K)$ representing any $\mu \in \mathcal K$
coincides with $\rho \circ \mu$.
Therefore an application of Theorem~\ref {csselectsl} to the compact set $\mathcal K$ yields Theorem~\ref {csselectm}.
Of course, Theorem~\ref {csselectm} can also be verified directly without any use or mention of the Choquet theory by repeating the proof
of Theorem~\ref {csselect} with necessary changes.

\section {Sets closed in measure}

\label {scm}
In this section we establish an analogue of Proposition \ref {csgs} for sets that are convex, bounded and closed in measure
in a Banach ideal space of measurable functions, which makes it possible to apply some of the main results of this paper in this
somewhat obscure but highly useful setting.  For more detail on basic definitions and properties see, e.~g.,~\cite {kantorovichold}.

Let $(\Omega, \mu)$ be a $\sigma$-finite measurable space, and $\sclass {\Omega} {\mu}$ be the set of all measurable functions on $(\Omega, \mu)$
with its usual lattice structure.
It is said that $X \subset \sclass {\Omega} {\mu}$ is an
\emph {ideal\footnote {The word ``ideal'' in this term
is often omitted but understood in the literature working with general lattices of measurable functions;
we keep it here to avoid confusion.}
lattice of measurable functions on the measurable space $(\Omega, \mu)$}
if $X$ is an ideal in $\sclass {\Omega} {\mu}$, i. e. if conditions $f \in X$, $g \in \sclass {\Omega} {\mu}$ and $|g| \leqslant |f|$ a. e.
imply that $g \in X$, and the support $\supp X = \bigcup_{f \in X} \supp f$ of $X$ coincides with the entire set $\Omega$ up to a set of $\mu$-measure $0$.
It is said that $X$ is a \emph {normed ideal lattice} of measurable functions on $(\Omega, \mu)$
if $X$ is an ideal lattice equipped with a norm $\|\cdot\|_X$ such that for
any $f \in X$ and $g \in \sclass {\Omega} {\mu}$ satisfying $|g| \leqslant f$ one has $\|g\|_X \leqslant C \|f\|_X$ with some constant $C_X$ independent
of $f$ and $g$.  We say that $X$ is a \emph {Banach ideal lattice} of measurable functions on $(\Omega, \mu)$ if $X$ is a normed ideal lattice
that is a Banach space, i. e. if $X$ complete with respect to the norm $\|\cdot\|_X$.

For a Banach ideal lattice $X$ of measurable functions on $(\Omega, \mu)$
any order continuous functional $f$ on $X$ (order continuity means that given a sequence $x_n \in X$
such that $\sup_n |x_n| \in X$ and $x_n \to 0$ a. e. it follows that $f (x_n) \to 0$) has an integral representation
$f (x) = \int x y_f$ for some function $y_f$ that is identified with $f$.
The set of all such functionals $X'$ is a Banach ideal lattice with
the norm defined by $\|f\|_{X'} = \sup_{g \in X, \|g\|_X = 1} \int |f g|$.
The lattice $X'$ is called the \emph {order dual} of the lattice $X$.  If $X$ is a Banach ideal lattice then every functional from $X'$
is continuous on $X$.

A Banach ideal lattice of measurable functions on $(\Sigma, \mu)$ is said to satisfy \emph {the Fatou property} if
for any
$f_n, f \in X$ such that $\|f_n\|_X \leqslant 1$ and the sequence $f_n$ 
converges to $f$ a. e. it is also true that $f \in X$ and $\|f\|_X \leqslant 1$.
We can always assume that $C_X = 1$ if $X$ has the Fatou property.
The Fatou property of a lattice $X$ is equivalent to $\mu$-closedness of the unit ball
$B_X$ of the lattice $X$
(here and elsewhere $\mu$-convergence denotes 
convergence in measure in all measurable sets $E \subset \Omega$ satisfying $\mu (E) < \infty$).
If $X$ is a Banach ideal lattice then the Fatou property of $X$ is equivalent to order reflexivity of $X$, i. e. to the relation $X'' = X$.
If $X$ is a Banach ideal lattice having the Fatou property then the lattice $X'$ is a norming set of functionals for $X$,
i. e. $\|f\|_X = \sup_{g \in X', \|g\|_{X'} = 1} \int f g$ for all $f \in X$.

It is easy to see that $\mu$-closedness in a Banach ideal lattice $X$ satisfying the Fatou property is not weaker than norm-closedness,
and it is often strictly stronger.  For example, the set of continuous functions~$f$ satisfying $|f| \leqslant 1$
is norm-closed in the lattice $X = \lclass {\infty} {[0, 1]}$ but its closure in measure coincides with the unit ball of $X$.
Sets that are convex, bounded and closed in measure in a Banach ideal lattice $X$ having the Fatou property behave in some ways like compact
sets.  It is easy to see that the collection of such sets $\mathcal V (X)$ contains all compact convex sets of $X$ (in the strong topology of $X$)
since every sequence $f_n \in K$
in a compact set $K \subset X$ converging in measure to some $f \in X$ has a subsequence converging in $X$ and thus in measure to some $g \in K$
and therefore $f = g \in K$.  Also observe that
$\mathcal V (X)$ is contained in the collection of nonempty convex sets of $X$ that are compact in the weak topology
of $X$ since sets from $\mathcal V (X)$ are closed in $X$; therefore $\mathcal V (X)$ is an intermediate space between the spaces of
convex sets that are closed in the strong and the weak topologies of $X$.
The following theorem shows that $\mathcal V (X)$ has compact type.
\begin {theorem}[see {\cite[Theorem 3, Chapter X, \S5]{kantorovichold}}]
\label {fatoucompact}
Let $X$ be a Banach ideal lattice on $(\Sigma, \mu)$ having the Fatou property
and let $\{V_{\xi}\}_{\xi \in \Xi}$ be a centered family of sets in $X$ that are convex, bounded and $\mu$-closed.
Then $\bigcap_{\xi \in \Xi} V_{\xi}$ is not empty.
\end {theorem}

The following proposition that we are going to use repeatedly
is, essentially, rather well known (cf., e. g., \cite [1.2] {kisliakov2002en}, \cite [Proposition 15] {rutsky2011en}
or \cite {bukhvalovlozanovsky1977}).
\begin {proposition}
\label {fatoucs}
Suppose that $X$ is a Banach lattice of measurable functions on $(\Omega, \mu)$ having the Fatou property and
$F_n \subset X$ is a decreasing sequence of nonempty sets in $X$ that are bounded, convex and
$\mu$-closed.  Then for any sequence $f_n \in F_n$
there exists a sequence of finite convex combinations
$$
g_n = \sum_{k \geqslant n} \alpha_k^{(n)} f_k, \quad \alpha_k^{(n)} \geqslant 0, \quad \sum_{k \geqslant n} \alpha_k^{(n)} = 1
$$
such that $g_n$
converges a. e. to some $f \in F = \bigcap_n F_n$.
\end {proposition}

In order to prove Proposition~\ref {fatoucs} we introduce a nonincreasing sequence of sets
$$
A_n = \clos \co_{k \geqslant n} \{ f_k \},
$$
where $\clos$ denotes the closure in measure
and $\co$ denotes the convex hull.
By the assumptions of the proposition $A_n \subset F_n$
for all $n$, so Proposition \ref {fatoucompact} applied to the family $\{A_n\}$ indicates that
$\bigcap_n A_n$ is not empty.  Take any $f \in \bigcap_n A_n$.
Then $f \in F$ and for every $n$ there exists a sequence $g_k^n \in \co_{k \geqslant n} \{f_k \}$
such that $g_k^n$ converges to $f$ in measure on all sets of finite measure as
$k \to \infty$.  Since convergence in measure is metrizable, one can easily find a sequence of indices $k_n$
such that $g_{k_n}^n$ converges to $f$ in measure on all sets of finite measure.
Therefore some subsequence of this sequence converges to $f$ almost everywhere, which concludes the proof of Proposition \ref {fatoucs}.

\begin {proposition}
\label {csgsl}
Suppose that $X$ is a Banach ideal lattice of measurable functions on a $\sigma$-finite measurable space $(\Omega, \mu)$
and a semilinear space $\mathcal V = \mathcal V (X)$
consists of all nonempty bounded convex sets in $X$ that are $\mu$-closed.
Then $\mathcal V$ is a subspace of the semilinear space $\mathcal C (X)$, i. e. $\mathcal V$ is closed under the semilinear operations of
$\mathcal C (X)$.
$\mathcal V$ has compact type and is order regular.
The set $X'$ of all order continuous linear functionals on $X$ is exhaustive for $X$, and all $f \in X'$
are defining for~$\mathcal V$.
Suppose also that $\Omega = \mathbb Z$.  Then the set of point evaluation functionals $\mathcal D = \{f_j\}_{j \in \mathbb Z}$,
$f_j (\{x_k\}) = x_j$, is exhaustive for $X$ and consistent with~$\mathcal V$.
\end {proposition}

The proof is quite similar to the proof of Proposition \ref {csgs}.
First, suppose that $A, B \in \mathcal V$; we need to show that $A + B$ is closed in measure and thus $A + B \in \mathcal V$.
Let $z \in \clos (A + B)$.  Then there exists a sequence $z_n \in A + B$, $z_n = x_n + y_n$,
$x_n \in A$, $y_n \in B$, such that $z_n \to z$ in measure.  By Proposition \ref {fatoucs} there exists a sequence of finite convex combinations
$a_n = \sum_{m \geqslant n} \alpha_{n, m} x_n$ converging to some $x \in A$ in $X$.  Applying Proposition \ref {fatoucs} again to the sequence
$b_n = \sum_{m \geqslant n} \alpha_{n, m} y_n$ allows us to conclude that there exists a sequence of convex combinations
$d_n = \sum_{m \geqslant n} \beta_{n, m} b_m$ converging in $X$ to some $y \in B$.
Let $c_n = \sum_{m \geqslant n} \beta_{n, m} a_m$; then $c_n \in A$, $d_n \in B$,
$c_n \to x$ and $d_n \to y$ in $X$, and therefore in measure.
At the same time $c_n + d_n$ is a sequence of convex combinations of $z_n$, so $c_n + d_n \to z$.  This means that
$z = x + y$, so $A + B$ is indeed closed in measure and $\mathcal V$ is a subspace of $\mathcal C (X)$.
Compact type of $\mathcal V$ follows from Theorem \ref {fatoucompact}.
We have $\inf \mathfrak A = \bigcap \mathfrak A$ for any lower directed family of sets $\mathfrak A \subset \mathcal V$,
and so the rest of order regularity is established exactly as in the proof of Proposition~\ref {csgs}.
Finally, the set $X'$ is exhaustive for $X$ because $X'$ is norming for $X$.  Let $f \in X'$ be an order continuous functional in $X$ and $A \in \mathcal V$.
Since $A$ is bounded and $f$ is continuous, $f (A)$ is also bounded and convex.
$A$ is closed in $X$ (because it is closed in measure), so by the Open Mapping Theorem $f (A)$ is closed and
thus $f (A) = [a, b]$ for some $a, b \in \mathbb R$, so $f$ is defining for $\mathcal V$.
Let us now verify the claims about the case $\Omega = \mathbb Z$.
It is obvious that $\mathcal D$ is exhaustive for $X$.
Since convergence in measure implies convergence everywhere in $\Omega$, functionals $f_j$ are continuous for all $j \in \mathbb Z$
with respect to convergence in measure and thus for every $A \in \mathcal V$ and $c \in \mathbb R$ nonemptiness of $f^{-1} (c) \cap A$
implies that $f^{-1} (c) \cap A$ is closed in measure and nonempty, so $f^{-1} (c) \cap A \in \mathcal V$ and $f_j$ is consistent with $\mathcal V$.
The proof of Proposition~\ref {csgsl} is complete.

Proposition \ref {csgsl} makes it possible to apply the main results of this paper to certain situations lacking compactness, for example,
to bounded sets in $X = \lsclass {1}$ (and other non-dual lattices) that are closed in measure.
Note that in general lattices order continuous functionals may not be consistent with $\mathcal V$ defined in Proposition \ref {csgsl};
for example, the mean $f (x) = \int_\Omega x (\omega) d\omega$ is an order continuous functional in $X = \lclass {1} {[0, 1]}$, but it is easy to see that
$$
0 \in \left(\clos f^{-1} (1) \cap B_{\lclass {1} {[0, 1]}}\right) \setminus \left(f^{-1} (1) \cap B_{\lclass {1} {[0, 1]}}\right),
$$
where $B_{\lclass {1} {[0, 1]}}$
denotes the closed unit ball of $X$.  In fact, it is easy to see that
no nontrivial order continuous functional on $X = \lclass {1} {[0, 1]}$ is consistent with
$\mathcal V = \mathcal V (X)$  described in Proposition \ref {csgsl}.
However, another technique of obtaining linear selections allows us to extend Theorem~\ref {everyls} to the case of maps
acting on the cone of nonnegative functions in $\lclassg {1}$ into $\mathcal V (X)$ under additional assumptions;
see Theorem~\ref {everylsl1} in Section~\ref {ssvm}.
It is unclear whether the conclusion of Theorem~\ref {everyls} holds true for this particular
space $\mathcal V (X)$ under general conditions or for general Banach ideal lattices.

\section {Semilinear spaces with cone-bases}

\label {conebases}
In \cite {smajdor1990} and \cite {smajdor1996} some results on additive and affine selections were obtained for maps defined on a cone
with a cone-basis.  In this section we will show how the main results of this paper become very clear and instantaneous under the assumption that
the semilinear space the maps are acting on has a cone-basis.
However, this assumption is stronger than the interpolation property, and in Section \ref {sotctc} below it will become clear that it is often violated
outside the case of cones in a finite-dimensional space.

A set $B \subset K$ in a semilinear space $K$ is said to be a \emph {cone-basis} of~$K$ if for every $x \in K$ there exist
some unique nonnegative coordinates
$\alpha_b (x) \geqslant 0$ that are zero for all $b \in B$ except for a finite number and $x = \sum_{b \in B} \alpha_b (x) b$.
\begin {proposition}
\label {coneprop}
Any cone with a cone-basis has the interpolation property.
\end {proposition}
Indeed, suppose that
$z \in K$ and we are given some
$$
\{x_j\}_{j = 1}^m, \{y_k\}_{k = 1}^n \subset K
$$
satisfying
$z = \sum_{j = 1}^m x_j = \sum_{k = 1}^n y_k$.  Then
\begin {gather*}
z = \sum_{b \in B} \alpha_b (z) b,\quad
x_j = \sum_{b \in B} \alpha_b (x_j) b, \quad 
y_k = \sum_{b \in B} \alpha_b (y_k) b.
\end {gather*}
Let $z_{jk} = \sum_{b \in B} \frac {\alpha_b (x_j) \alpha_b (y_k)} {\alpha_b (z)} b$; here we use the usual convention $\frac 0 0 = 0$
(note that $\alpha_b (z) = 0$ implies $\alpha_b (x_j) = 0$ and $\alpha_b (y_j) = 0$).
Surely,
$x_j = \sum_{k = 1}^n z_{jk}$ and $y_k = \sum_{j = 1}^m z_{jk}$ for all $j$ and $k$, so the interpolation property of $K$ is established.

On the other hand, not every cone with the interpolation property has a cone-basis.
In fact, Theorem~\ref {linfcb} in Section~\ref {sotctc} below
indicates that in no infinite-dimensional normed ideal lattice of measurable functions
the cone of nonnegative functions has a cone-basis.
The following proposition gives a very easy particular case of Theorem \ref {linfcb}.
\begin {proposition}
\label {l1notcb}
The cone $K = \left[ \lclass {1} {[0, 1]}\right]_+$ of nonnegative summable functions on the segment $[0, 1]$ has
the interpolation property, but it has no cone-basis.
\end {proposition}
The cone $K$ has the interpolation property as the cone of positive
functions in a lattice.  Let us show that $K$ does not have a cone-basis.
Indeed, suppose that $B \subset K$ is a cone-basis of $K$.  We can assume by renorming that
$\|b\|_{\lclassg {1}} = 1$ for all $b \in B$.  Then every function $z \in B$ has to be an extreme point of the unit ball of $\lclassg {1}$: if
$z = (1 - \theta) x + \theta y$ for some $0 < \theta < 1$ and $x, y \in \lclassg {1}$ such that $\|x\|_{\lclassg {1}} \leqslant 1$ and 
$\|y\|_{\lclassg {1}} \leqslant 1$, then we can make $x$ and $y$ positive (and thus belonging to $K$)
by replacing them with $x \vee 0$ and $y \vee 0$ respectively, and the definition of a cone-basis
implies that $\alpha_b (z) = (1 - \theta) \alpha_b (x) + \theta \alpha_b (y)$ for all $b \in B$,
which means that $x = y = b$.  But it is well known that the unit ball of $\lclassg {1}$ has no extreme points, so $K$ does not have a cone-basis.

This example shows that the existence of a cone-basis for cones in infinite-dimensional spaces is not that common.
In fact, in Theorem \ref {linfcb} below we show that the only ideal normed lattices in which the cone of nonnegative functions has
a cone-basis are finite-dimensional lattices.
This contrasts starkly with the fact that every linear space has a Hamel (i.~e. algebraic) basis.
However, the existence of a cone-basis 
allows one to parametrize the set of all linear selections of a superlinear set-valued map with almost no effort.
The following proposition follows immediately from the definitions.
\begin {proposition}
\label {cbrdpnn}
Suppose that $\mathcal S$ is a semilinear space having a cone-basis $B$ and
$T : \mathcal S \to \mathcal V$ is a superlinear map acting into a semilinear space $\mathcal V$ equipped with a transitive relation $\leqslant$
compatible with the semilinear structure of $\mathcal V$.
Then for every $f : B \to \mathcal V$ satisfying $f \leqslant T$ on the set $B$
the map defined by $S_f (x) = \sum_b \alpha_b (x) f (b)$ is a linear submap
$S_f \leqslant T$ of $T$ such that $S_f \leqslant f$ on $B$.  In particular,
if $f$ is the restriction of $T$ on $B$ then $S_T = S_f$ is the greatest linear submap of $T$.
\end {proposition}

The next proposition strengthens Theorem \ref {linselnd} in the case when the space $\mathcal S$ has a cone-basis.
\begin {proposition}
\label {conebise}
Suppose that $\mathcal S$ and $\mathcal W$ are semilinear spaces,
$\mathcal V$ is a subspace of the semilinear space $\mathcal C (\mathcal W)$ and $\mathcal S$ has a cone-basis $B$.
Let $A : \mathcal S \to \mathcal V$ be a linear map.
Then for any $x \in \mathcal S$ and $y \in A (x)$ there exists a linear selection $a : \mathcal S \to \mathcal W$ of $A$ satisfying
$a (x) = y$.
\end {proposition}
Indeed, let $x = \sum_{b \in B} \alpha_b (x) b$, and $B_x = \{b \in B \mid \alpha_b (x) \neq 0\}$.
Then $A (x) = \sum_{b \in B_x} \alpha_b (x) A (b)$ by linearity of $A$, so
$y = \sum_{b \in B_x} \alpha_b (x) y_b$ with some $y_b \in A (b)$ for $b \in B_x$.
Select arbitrary $y_b \in A (b)$ for all $b \in B \setminus B_x$ (possibly making use of the Axiom of Choice)
and define the mapping $a$ by $a (z) = \sum_{b \in B} \alpha_b (z) y_b$ for any $z = \sum_{b \in B} \alpha_b (z) b$.
Then
$$
a (z) \in \sum_{b \in B} \alpha_b (z) A (b) = A \left( \sum_{b \in B} \alpha_b (z) y_b \right) = A (z)
$$
and $a (x) = y$, which concludes the proof of Proposition \ref {conebise}.

Finally,
it is easy to obtain from the previous two propositions
a version of Theorem \ref {everyls} for the case when space $\mathcal S$ has a cone-basis.
\begin {proposition}
\label {everylscb}
Suppose that $\mathcal S$ and $\mathcal W$ are semilinear spaces, $\mathcal S$ has a cone-basis $B$ and
$\mathcal V$ is a subspace of the semilinear space $\mathcal C (\mathcal W)$.  Let $T : \mathcal S \to \mathcal V$ be a superlinear map.
Then the following conditions are equivalent.
\begin {enumerate}
\item
There exists a linear selection $a$ of $T$ satisfying $a (x) = y$.
\item
$y \in \sum_{b \in B} \alpha_b (x) T (b)$.
\end {enumerate}
\end {proposition}

Note that in the finite-dimensional case by the Krein-Milman theorem a compact convex set $K \subset \mathbb R^n$ is a Choquet (or Baire) simplex
if and only if the set $\mathcal E (K)$ of the extreme points of $K$ is a cone-basis in $K_{sus}$.
Thus it is possible to obtain Theorems \ref {latnecc} and \ref {riesznecc} in the case $K \subset \mathbb R^n$
(or, equivalently, $\mathcal S \subset \mathbb R^{n + 1}$)
from the simple results described in this section.


\section {Sharpness of the compact type condition}

\label {sotctc}

In this section we present an example that shows that the requirement
for $\mathcal V$ to have compact type cannot be dropped from the conditions of Theorem~\ref {linselnd}.
We also prove, using essentially the same idea,
that the cone of positive functions in an ideal normed lattice of measurable functions on a measurable space $(\Omega, \mu)$ 
has no cone-basis unless measure $\mu$ is discrete.

Let $\mathcal S = \left[\lsclass {\infty}\right]_+$ be the set of bounded sequences of nonnegative real numbers,
$c_0$ be the set of sequences that converge to 0 with the usual norm $\|\{a_n\}\|_{c_0} = \sup_n |a_n|$,
$\mathcal W = \left[c_0\right]_+$ be the cone of sequences of nonnegative numbers in $c_0$
and $\mathcal V \subset \mathcal C (\mathcal W)$ be the set of all nonempty bounded convex closed sets in $\mathcal W$.
It is easy to see that $\mathcal V$ does not have compact type: take, for example, a diverging sequence $a$ bounded by 1
and a decreasing sequence of sets
$$
A_n = \left\{ \{b_j\} \in \mathcal W \mid \text {$|b_j| \leqslant 1$ for all $j$}, \text {$b_j = a_j$ for $1 \leqslant j \leqslant n$} \right\}
$$
that all lie in $\mathcal V$ but their intersection $\bigcap_n A_n = \{a\}$ is not a subset of $\mathcal W$.
However, $\mathcal V$ is a subspace of $\mathcal C (\mathcal W)$ that has
an exhaustive set of linear functionals $f_j (\{a_k\}) = a_j$ which are defining for and consistent with $\mathcal V$.
Let us define a map $A : \mathcal S \to \mathcal V$ by $A (\{a_j\}) = \mathcal W \cap \prod_j [0, a_j]$.
Linearity of $T$ is easily verified, but, as we will now show, for any sequence 
$x = \{x_j\} \in c_0$ such that $x_j > 0$ for all $j \in \mathbb N$ there is no linear selection $a$ of $A$ satisfying $a (x) = x$.
Suppose, on the contrary, that there is such a selection $a$.
For any $I \subset \mathbb N$ we denote by $\chi_I$ the indicator function
for $I$ on $\mathbb N$.
As $a (\chi_I x) \leqslant \chi_I x$ for any $I \subset \mathbb N$, it follows from additivity of $a$ that $a (\chi_I x) = \chi_I x$
for any $I \subset \mathbb N$.
Let $e_j = \chi_{\{k \mid k = j\}} \in c_o$ be the usual coordinate vectors.
Then $a (e_j) = e_j$, so for any finite $I \subset \mathbb N$ and any $z \in \mathcal S$ we have $a (\chi_I z) = \chi_I z$, and
$a (z) = a (\chi_I z) + a (\chi_{\mathbb N \setminus I} z) \geqslant \chi_I z$.  Therefore $a (z) = z$ for all $z \in \mathcal S$,
a contradiction since, for example, $z = \chi_{\mathbb N} \notin \mathcal W$.
Thus we have constructed an example of semilinear spaces $\mathcal S$, $\mathcal W$, $\mathcal V$ and a map $A$ that satisfy all conditions
of Theorem \ref {linselnd} except for the compact type of $\mathcal V$ and for which the conclusion of Theorem \ref {linselnd} does not hold.

Let $\nu$ be the counting measure on $\mathbb N$.
The construction just described actually shows that for no Banach ideal lattice $X$ on the measurable space $(\mathbb N, \nu)$ the cone
$\mathcal S = X_+$ of nonnegative functions has a cone-basis.
Let $\mathcal W$ be any ideal lattice of measurable functions which is a proper sublattice of $X$.
We can take, for example,
$$
\mathcal W = \left\{ \{a_n\} \in X \mid \left\{n a_n \right\} \in X \right\},
$$
which is a Banach ideal lattice with respect to the norm
$$
\left\|\{b_n\}\right\|_{\mathcal W} = \left\|\{n a_n\}\right\|_X.
$$
It is easy to see that $\mathcal W$ is a proper sublattice of $X$: for example, sequence $n^{-\frac 1 2} \frac {e_n} {\|e_n\|_X}$ converges to
$0$ in $X$ but diverges in $\mathcal W$.
Suppose that $\mathcal S = X_+$ has a cone-basis.  Then by Proposition \ref {conebise} map $A$ would have had a linear selection $a$ satisfying
$a (x) = x$ for some $x \in \mathcal W_+$.  There exists some $x \in \mathcal W_+$ such that $\supp x = \mathbb N$; for example,
$x = \left\{\frac {2^{-n}} {\|e_n\|_{\mathcal W}}\right\}$.  Then the argument from the example shows that $a (x) = x$ for all $x \in \mathcal S$ and thus
$X_+ = \mathcal W_+$, a contradiction.
Therefore for no Banach ideal lattice $X$ on $(\mathbb N, \nu)$ the cone
$X_+$ has a cone-basis.  This result can be extended to all nondiscrete measurable spaces.


\begin {theorem}
\label {linfcb}
Suppose that $X$ is a normed ideal lattice of measurable functions on a $\sigma$-finite measurable space $(\Omega, \mu)$.
The following conditions are equivalent.
\begin {enumerate}
\item
The cone of nonnegative functions $X_+$ in $X$ has a cone-basis.
\item
Measure $\mu$ is discrete, i. e. $\mu$ is a finite sum of point masses $\mu = \sum_{j = 1}^N c_j \delta_{\omega_j}$, $\omega_j \in \Omega$.
\item
$X$ is finite-dimensional.
\end {enumerate}
\end {theorem}
Implication $2 \Rightarrow 1$ is straightforward: if condition 2 is satisfied then $W = \{ \chi_{\omega_j} \}_{j = 1}^N$ is a cone-basis for $X_+$.
Equivalence of conditions 2 and 3 is also simple.  If condition 2 is satisfied then $W$ is a basis for $X$,
so implication $2 \Rightarrow 3$ holds true. On the other hand, if measure $\mu$
is nondiscrete then there exists a countable family $\{E_j\}_{j \in \mathbb N}$ of disjoint measurable sets in $\Omega$ having positive finite measure,
and with it an infinite family of functions from $X$ supported in $E_j$ that is linearly independent, so if condition 2 is not satisfied then condition
3 is also not satisfied, which proves implication $3 \Rightarrow 2$.
Let us establish the remaining implication $1 \Rightarrow 2$.  Suppose that it does not hold for some measure $\mu$, i. e. $\mu$ is nondiscrete but
$X_+$ has a cone-basis; we will now devise a construction from these two assumptions that would lead to a contradiction.
Since measure $\mu$ is nondiscrete and $\sigma$-finite there exists a countable partition $\mathcal E = \{E_n\}_{n \in \mathbb N}$ of $\Omega$
(i. e. sets $E_n$ are pairwise disjoint and $\bigcup_n E_n = \Omega$) such that $\chi_{E_n} \in X$
(see, e. g., \cite [Chapter IV, \S 3, Corollary 2] {kantorovichen} or \cite {kantorovichold}).
If $\weightw$ is a measurable function which is positive a. e. the weighted normed lattice $X (\weightw)$ is defined by
$X (\weightw) = \{ f \weightw \mid f \in X\}$ with the weighted norm $\|g\|_{X (\weightw)} = \|g \weightw^{-1}\|_X$.
Let $\weightw = \sum_{n \in \mathbb N} 2^{-n} \chi_{E_n}$, $\mathcal W = \left[X (\weightw)\right]_+$ and $\mathcal S = X_+$.
It is easy to see that $X (\weightw) \subset X$ (continuous inclusion) because $\weightw \leqslant 1$, and therefore $\mathcal W \subset \mathcal S$.
Define $\mathcal V \subset \mathcal C (\mathcal W)$ by
$\mathcal V = \left\{ \{f \in \mathcal W \mid 0 \leqslant f \leqslant a\} \mid a \in \mathcal W \right\}$,
and define a map $A : \mathcal S \to \mathcal V$ by $A (g) = \{f \in \mathcal W \mid 0 \leqslant f \leqslant g\}$.
Let 
$g = \sum_{n \in \mathbb N} 2^{-n} \frac {\chi_{E_n}} {\|\chi_{E_n}\|_{X (\weightw)}}$.
Since by assumptions $X_+ = \mathcal S$ has a cone-basis, by Proposition \ref {conebise} there exists a linear selection $a$ of $A$ satisfying
$a (g) = g$.  Observe that for any $f \in X_+$ we have $a (f) \leqslant f$.  Therefore, for any measurable $E \subset \Omega$
we have $a (g \chi_E) \leqslant g \chi_E$ and $a (g \chi_{\Omega \setminus E}) \leqslant g \chi_{\Omega \setminus E}$.  But
$g = g \chi_E + g \chi_{\Omega \setminus E}$ and $a (g) = g$, so $a (g \chi_E) = g \chi_E$; in particular, $a (\chi_{E_n}) = \chi_{E_n}$
for any~$n$.
Let $h = \sum_{n \in \mathbb N} 2^{-\frac n 2} \frac {\chi_{E_n}} {\|\chi_{E_n}\|_{X}}$.  It is easy to see that $h \in \mathcal S$.
By the assumptions, $a (h) \in \mathcal W$, but
\begin {multline*}
\|a (h)\|_{X (\weightw)} = \|a (h \chi_{E_n}) + a (h \chi_{\Omega \setminus E_n})\|_{X (\weightw)} \geqslant \|a (h \chi_{E_n})\|_{X (\weightw)} =
\\
\left\| \frac {2^{-\frac n 2}} {\|\chi_{E_n}\|_X} a (\chi_{E_n})\right\|_{X (\weightw)} =
\left\| \frac {2^{-\frac n 2}} {\|\chi_{E_n}\|_X} \chi_{E_n}\right\|_{X (\weightw)} =
\\
\left\| \frac {2^{\frac n 2}} {\|\chi_{E_n}\|_X} \chi_{E_n}\right\|_{X} = 2^{\frac n 2}
\end {multline*}
for any $n \in \mathbb N$, so $a (h) \notin \mathcal W$.  This contradiction concludes the proof of Theorem \ref {linfcb}.

\section {Linear and superlinear set-valued maps of $\left[\lclassg {1}\right]_+$}

\label {ssvm}
In Section~\ref {mosom} above we have seen that for upper semicontinuous superlinear maps acting on the semilinear space of Borel measures
$\mathcal M_{+} (K)$ on a compact set $K$ the conditions necessary for the existence of linear selections are significantly weaker compared to the
general case treated in the main results of this paper, and indeed we have a natural parametrization of the set of all linear selections in this case.
The semilinear space that often arises in applications, however, is
the cone $\mathcal S = \left[\lclassg {1}\right]_+$ of measurable functions on some measurable space $(\Omega, \mu)$ with a $\sigma$-finite measure $\mu$
that are positive a.~e. and summable.  In this section we will see how in certain cases the main results of this survey can be extended in this setting;
the most useful, perhaps, is a simple construction proving the existence of a linear selection for every bounded linear map
$T : \left[\lclassg {1}\right]_+ \to \mathcal C (X)$ acting into a Banach space $X$ if the measurable space $\Omega$ is separable.
Additional assumptions provide a result similar to Theorem~\ref {linselnd} for the case of maps taking values in sets closed in measure.

Let $\mathcal S$ be a monoid semilinear space. 
We say that
$$
\mathcal B = \{b_{n, k}\}_{{n \geqslant 0,} \atop {0 \leqslant k < 2^n}} \subset \mathcal S
$$
is a \emph {nesting basis} of $\mathcal S$ if it satisfies the following conditions.
\begin {enumerate}
\item
$b_{n, k} = b_{n + 1, 2 k} + b_{n + 1, 2 k + 1}$ for all $n$ and $k$.
\item
$b_{n_1, k_1} = b_{n_2, k_2}$ implies that either $n_1 = n_2$ and $k_1 = k_2$
\linebreak{}
or $b_{n_1, k_1} = b_{n_2, k_2} = 0$.
\item
The conic hull of $\mathcal B$ coincides with $\mathcal S$.
\end {enumerate}
Let us give an example right away.  Suppose that $\mathfrak A = (\Omega, \Sigma, \mu)$ is a separable measurable space.
This means that there exists a sequence of sets $\{A_j\}_{j \geqslant 1} \subset \Sigma$ such that the minimal $\sigma$-algebra
$\Sigma_0 = \sigma \left( \{A_j\}_{j \geqslant 1} \right)$ containing all of them is $\mu$-dense in $\mathfrak A$.
By renorming measure $\mu$ if necessary we can assume that $\mu (\Omega) = 1$.
We will now construct a nesting basis for the space $\mathcal S$ of simple functions in $\left[ \lclass {1} {\Omega, \Sigma_0, \mu} \right]_+$.
Let $b_{0, 0} = \chi_{\Omega}$
and $\Sigma_n = \sigma \left( \{A_j\}_{1 \leqslant j < n}\right)$; this $\sigma$-algebra consists from disjoint unions of at most $2^n$ atoms
that we are going to denote by~
$B_{n, k}$, $0 \leqslant k < 2^n$, possibly with empty sets in certain places.  We can also assume that
$B_{n + 1, 2 k} \cup B_{n + 1, 2 k + 1} \subset B_{n, k}$ for all $n$ and $k$.  Then functions $b_{n, k} = \chi_{B_{n, k}}$ form a nesting basis.
For example, for the unit segment $\Omega = [0, 1]$ with a Lebesgue measure we can take the usual dyadic partition
$B_{n, k} = \left[2^{-n} k, 2^{-n} (k + 1)\right]$.
\begin {proposition}
\label {nbassel}
Suppose that $\mathcal S$ is a monoid semilinear space that has a nesting basis $\mathcal B = \{b_{n, k}\}_{{n \geqslant 0,} \atop {0 \leqslant k < 2^n}}$,
$\mathcal W$ is a semilinear space, let $T : \mathcal S \to 2^{\mathcal W} \setminus \{\emptyset\}$ be a linear map and $y_0 \in T (b_{0, 0})$.
Then $T$ has a linear selection satisfying $S (b_{0, 0}) = y_0$.
\end {proposition}
Indeed, we can define a partial additive selection $S_0 : \mathcal B \to \mathcal W$ of $T$
inductively by selecting $S_0 \left(b_{0, 0}\right) = y_0 \in T \left(b_{0, 0}\right)$ and
\begin {multline*}
\left(S_0 \left(b_{n + 1, 2 k} \right), S_0 \left(b_{n + 1, 2 k + 1}\right)\right) \in
\\
\left\{ (g, h) \mid
g \in T \left(b_{n + 1, 2 k}\right),
h \in T \left(b_{n + 1, 2 k + 1}\right),
g + h = S_0 \left(b_{n, k}\right)
\right\};
\end {multline*}
this set is nonempty because $T \left(b_{n + 1, 2 k}\right) + T \left(b_{n + 1, 2 k + 1}\right) = T \left(b_{n, k}\right)$ due to the linearity of $T$.
We then define a suitable linear selection $S$~of~$T$~by
$$
S \left(\sum_{j = 1}^N \alpha_j b_{n_j, k_j} \right) = \sum_{j = 1}^N \alpha_j S_0 \left( b_{n_j, k_j}\right)
$$
for all finite sequences of $n_j$, $k_j$ and $\alpha_j \geqslant 0$.

Let $X$ and $Y$ be normed spaces, and $\mathcal S \subset X$  We say that a map $T : \mathcal S \to 2^Y$ is \emph {bounded} if
$\sup_{y \in T (x)} \|y\|_Y \leqslant C \|x\|_X$ for all $x \in \mathcal S$ with some constant $C$ independent of $x$.
Boundedness sometimes allows us to extend selections from dense sets.
\begin {theorem}
\label {blmls}
Suppose that $X$ is a Banach space, $(\Omega, \mu)$ is a separable measurable space
and $T : \left[\lclassg {1} {(\Omega, \mu)}\right]_+ \to 2^X$ is a bounded upper semicontinuous linear map.
Then $T$ has a continuous linear selection.
\end {theorem}
Indeed, we use the construction before Proposition~\ref {nbassel} to obtain a $\sigma$-algebra $\left(\Omega, \Sigma_0, \mu\right)$
that is $\mu$-dense in $\left(\Omega, \Sigma, \mu\right)$ and a nesting basis $\mathcal B$ for
the semilinear space $\mathcal S_0$ of simple functions on $\left[\lclass {1} {\Omega, \Sigma_0, \mu}\right]_+$.  By Proposition~\ref {nbassel}
the restriction of $T$ onto $\mathcal S_0$ has a linear selection $S_0 : \mathcal S_0 \to X$ which is bounded.
By Proposition~\ref {gccont} the map $S_0$ is naturally extended to a bounded linear selection $S_1$ of $T$
on the lattice $\lclass {1} {\Omega, \Sigma_0, \mu}$.  Then $S_1$ can be extended by continuity to a linear bounded and hence continuous map
$S : \lclass {1} {\Omega, \mu} \to X$ which we restrict onto the cone $\mathcal S = \left[\lclass {1} {\Omega, \mu}\right]_+$;
the graph of $S$ is the closure of the graph of $S_1$ in $\mathcal S \times X$.
By Proposition~\ref {grclos} the graph of $T$ is closed, so the graph of $S$ is contained within the graph of $T$, which means that
$S$ is a bounded linear selection of $T$.

In order to obtain an analogue of Theorem~\ref {linselnd} we will need the following proposition, which shows that at least in the present setting
upper semicontinuous linear maps also posess a sort of lower semicontinuity.
\begin {proposition}
\label {lmils}
Suppose that $(\Omega, \mu)$ is a measurable space,
$Y$ is a Banach space
and $T : \left[\lclass {1} {\Omega, \mu}\right]_+ \to 2^{Y}$ is a bounded linear map.
Then for evey $f \in \left[\lclass {1} {\Omega, \mu}\right]_+$, $g \in T (f)$ and a sequence $f_n \in \left[\lclass {1} {\Omega, \mu}\right]_+$
converging to $f$ in $\lclass {1} {\Omega, \mu}$
there exists a sequence $\left[\lclass {1} {\Omega, \mu}\right]_+$ of convex combinations of $\{f_j\}_{j \geqslant n}$ and $g_n \in T (h_n)$ such that
$h_n \to f$ in $\lclass {1} {\Omega, \mu}$ and $g_n \to g$ in $Y$.
\end {proposition}
Indeed, let $\varepsilon_n \to 0$.  It suffices to find a a convex combination $h_n$ of $\{f_j\}_{j \geqslant n}$ and $g_n \in T (h_n)$
satisfying $\|g_n - g\|_Y \leqslant \varepsilon_n$.
Suppose that no such $h_n$ and $g_n$ exist; this means that the distance between the closed convex hull $G_n$ of $\{T (f_j) \mid j \geqslant n\}$
and $g$ is at least $\varepsilon_n$, so $g \notin G_n$.  By a separation theorem there exists some bounded functional $\varphi \in Y^*$ on~$Y$
satisfying
\begin {equation}
\label {supnc}
\sup \{\varphi (\mathrm g) \mid \mathrm g \in G_n \} < \varphi (g).
\end {equation}
By Proposition~\ref {swfl} support function $\varphi_T^*$ defined by $\varphi_T^* (x) = \sup_{y \in T (x)} \varphi (y)$
is bounded and linear.
Proposition~\ref {gccont} allows us to extend $\varphi_T^*$ boundedly onto $\lclass {1} {\Omega, \mu}$,
and therefore $\varphi_T^*$ is continuous.
This contradicts~\eqref {supnc} since the right-hand part is bounded from above by $\varphi_T^* (f)$
but at the same time the left-hand part is at least $\liminf_n \varphi_T^* (f_n) = \varphi_T^* (f)$.  The proof of Proposition~\ref {lmils} is complete.

We also need the following refinement of Proposition~\ref {nbassel}.
\begin {proposition}
\label {nbasselr}
Suppose that $\mathcal S$ is a monoid semilinear space that has a nesting basis $\mathcal B = \{b_{n, k}\}_{{n \geqslant 0,} \atop {0 \leqslant k < 2^n}}$,
$\mathcal W$ is a semilinear space and $T : \mathcal S \to 2^{\mathcal W} \setminus \{\emptyset\}$ is a linear map.
Then for every $x \in \mathcal S$ and $y \in T (x)$ there exists a linear selection $S$ of $T$ satisfying $S (x) = y$.
\end {proposition}
Indeed, we can write
\begin {equation}
\label {xnbe}
x = \sum_{j = 1}^N \alpha_j b_{n_j, k_j}
\end {equation}
with $\alpha_j > 0$ and an additional assumption that
\begin {equation}
\label {currel}
b_{n_j, k_j} \prec b_{n_l, k_l} \text {  for no $j$ and $l$;}
\end {equation}
here $\prec$ is the intrinsic order of the semilinear space $\mathcal S$ introduced in Section~\ref {introduction}.
Otherwise $b_{n_l, k_l}$ is expressed as a sum of $b_{n_j, k_j}$ and a finite number of elements from $\mathcal B$ that we can substitute
into \eqref {xnbe}; it is easy to see that every such substitution decreases the total number of relations \eqref {currel} and so
after a finite number of such substitutions condition \eqref {currel} is satisfied.
Similarly, we also express $x_r = b_{0, 0} - \sum_{j = 1}^N b_{n_j, k_j}$ as
$x_r = \sum_{j = N + 1}^M b_{n_j, k_j}$.
Because of the linearity of $T$ there exist some $y_{n_j, k_j} \in T \left(b_{n_j, k_j}\right)$ satisfying
$y = \sum_{j = 1}^N \alpha_j y_{n_j, k_j}$.  We then apply Propositon~\ref {nbassel} to subspaces
$\mathcal S_{n_j, k_j} = \left\{\lambda z \in \mathcal S \mid z \preceq b_{n_j, k_j}, \lambda > 0 \right\}$
with their respective bases
$$
\mathcal B_{n_j, k_j} = \left\{b \in \mathcal B \mid b \preceq b_{n_j, k_j} \right\}
$$
to obtain for all $1 \leqslant j \leqslant M$
linear selections $S_{n_j, k_j}$ of $T$ restricted onto $\mathcal S_{n_j, k_j}$ satisfying $S_{n_j, k_j} (b_{n_j, k_j}) = y_{n_j, k_j}$ for $1 \leqslant j \leqslant N$.
It is easy to see that the conic hull of $\bigcup_{j = 1}^M \mathcal S_{n_j, k_j}$ is the entire $\mathcal S$.
Now we can define a linear selection $S$ of $T$ by
$S (z) = \sum_{j = 1}^M S_{n_j, k_j} (z)$, $z \in \mathcal S$, which by the construction satisfies $S (x) = y$.

\begin {theorem}
\label {linselndl1}
Suppose that $X$ is a Banach ideal lattice of measurable functions on a $\sigma$-finite measurable space $(\Omega_Y, \nu)$, $X$ satisfies the
Fatou property,
$(\Omega, \mu)$ is a separable measurable space,
the semilinear space $\mathcal V$ consists of all nonempty bounded convex sets in $X$ that are $\nu$-closed
and $T : \left[\lclassg {1} {(\Omega, \mu)}\right]_+ \to \mathcal V$ is a bounded upper semicontinuous linear map.
Then for any $f \in \left[\lclassg {1} {(\Omega, \mu)}\right]_+$ and $g \in T (f)$ there exists a continuous linear selection $S$ of $T$ satisfying
$S (f) = g$.
\end {theorem}

Indeed, suppose that under the conditions of Theorem~\ref {linselndl1} we are given some
$f \in \left[\lclassg {1} {(\Omega, \mu)}\right]_+$ and $g \in T (f)$.
Following the proof of Theorem~\ref {blmls}, there exists a $\sigma$-algebra $\left(\Omega, \Sigma_0, \mu\right)$
that is $\mu$-dense in $\left(\Omega, \Sigma, \mu\right)$ and a nesting basis $\mathcal B$ for
the semilinear space $\mathcal S_0$ of simple functions on $\left[\lclass {1} {\Omega, \Sigma_0, \mu}\right]_+$.
Then there exists a sequence $f_n \in \mathcal S_0$ converging to $f$ in $\lclassg {1}$.
By Proposition~\ref {lmils}
there exists a sequence $h_n \in \left[\lclass {1} {\Omega, \mu}\right]_+$ of convex combinations of $\{f_j\}_{j \geqslant n}$ and $g_n \in T (h_n)$ such that
$h_n \to f$ in $\lclass {1} {\Omega, \mu}$ and $g_n \to g$ in $X$.
Then by Proposition~\ref {nbasselr} there exist linear selections $S_n$ of $T$ restricted onto $\mathcal S_n$ satisfying $S_n h_n = g_n$,
which by Proposition~\ref {gccont} can be extended by continuity to bounded maps on the entire $\lclass {1} {\Omega, \mu}$, which are uniformly bounded
because they are selections of a bounded map $T$.
Then it is evident that the family $\{S_n\}$ is uniformly continuous.
We define a map $R : \left[\lclassg {1} {(\Omega, \mu)}\right]_+ \to \mathcal V$ by
$R_n (x) = \clos \co \{S_j (x) \mid j \geqslant n\}$, $R (x) = \bigcap_n R_n (x)$ for $x \in \left[\lclassg {1} {(\Omega, \mu)}\right]_+$,
where $\clos$ denotes $\nu$-closure.
Then $R$ is a submap of $T$ having nonempty values
belonging to $\mathcal V$ because of the compact type of $\mathcal V$ guaranteed by Proposition~\ref {csgsl}.
We may assume that $\nu \left(\Omega_Y\right) < \infty$ by replacing the measure $\nu$ by an equivalent one.
It is well-known that convergence in measure is metrizable; let
\begin {equation}
\label {deltametric}
\delta (x_1, x_2) = \int_{\Omega_Y} \frac {|x_1 - x_2|} {1 + |x_1 - x_2|} d\nu, \quad x_1, x_2 \in X
\end {equation}
be a tranlation-invariant metric for this convergence in $X$.  It is also well-known that convergence in $X$ implies convergence in measure for normed
ideal lattices $X$; see, e.~g., \cite [Chapter IV, \S 3, Theorem 1] {kantorovichen} or~\cite {kantorovichold}.
Therefore the inclusion $X \to X_\delta$ is continuous and thus bounded, where $X_\delta$ denotes the space $X$ with metric $\delta$;
together with translation invariance of $\delta$ this means that there exists a constant $C_X$ such that
$\delta (x_1, x_2) \leqslant C_X \|x_1 - x_2\|_X$ for all $x_1, x_2 \in X$.
We are now going to verify that $R (f) = \{g\}$.  Take arbitrary $\varepsilon > 0$.
There exists some $N \in \mathbb N$ such that $\|f_n - f\|_{\lclassg {1}} \leqslant \varepsilon$ and $\|g_n -  g\|_{X} \leqslant \varepsilon$
for all $n \geqslant N$.
Then $\|S_n f - g_n\|_{X} = \|S_n f - S_n f_n\|_{X} \leqslant C \varepsilon$ for all $n \geqslant N$
with some constant $C$ independent of $\varepsilon$ by the uniform boundedness of $\{S_n\}$, and therefore
$\|S_n f - g\|_{X} \leqslant \|S_n f - g_n\|_{X} + \|g_n - g\|_{X} \leqslant (C + 1) \varepsilon$ for all~$n \geqslant N$.
This implies that $\delta (S_n f, g) \leqslant C_X (C + 1) \varepsilon$ for $n \geqslant N$
and therefore by the convexity of \eqref {deltametric} we also have
$\sup_{y \in R_n (f)} \delta (y, g) \leqslant C_X (C + 1) \varepsilon$ for all~$n \geqslant N$.
Since $\varepsilon$ is arbitrary, this means that $R (f) = \{g\}$ as claimed.
It is easy to see that the graphs of $R_n$ are closed in
$\left[\lclass {1} {\Omega, \mu}\right]_+ \times X_\delta$, thus
the graphs of $R_n$ are closed in a stronger topology $\left[\lclass {1} {\Omega, \mu}\right]_+ \times X$.
Therefore the graph of $R$ is closed in $\left[\lclass {1} {\Omega, \mu}\right]_+ \times Y$, and by Proposition~\ref {grclos}
map $R$ is upper semicontinuous.  Finally, we apply Theorem~\ref {blmls} to obtain a continuous linear selection $S$ of $R$.
Then $S (f) = g$ and $S$ is a selection of $T$.  The proof of Theorem~\ref {linselndl1} is complete.

Finally, we present a version of Theorem~\ref {everyls} for the case treated in this section.
\begin {theorem}
\label {everylsl1}
Suppose that $X$ is a Banach ideal lattice of measurable functions on a $\sigma$-finite measurable space $(\Omega_Y, \nu)$, $X$ satisfies the
Fatou property,
$(\Omega, \mu)$ is a separable measurable space,
the semilinear space~$\mathcal V$ consists of all nonempty bounded convex sets in $X$ that are $\nu$-closed
and $T : \left[\lclassg {1} {(\Omega, \mu)}\right]_+ \to \mathcal V$ is a bounded upper semicontinuous superlinear map.
Then $T$ admits a continuous linear selection.
Suppose, additionally, that $x \in \left[\lclassg {1} {(\Omega, \mu)}\right]_+$ and $y \in T (x)$.
Then the following conditions are equivalent.
\begin {enumerate}
\item
There exists a continuous linear selection $a$ of $T$ such that $a (x) = y$.
\item
$y \in \sum_{j = 1}^N T (x_j)$ for any $N \in \mathbb N$ and $\{x_j\}_{j = 1}^N \subset \left[\lclassg {1} {(\Omega, \mu)}\right]_+$
satisfying $x = \sum_{j = 1}^N x_j$.
\end {enumerate}
\end {theorem}
Indeed, by Theorem~\ref {linselg} (applicable due to Proposition~\ref {csgsl} and
the space $\left[\lclassg {1} {(\Omega, \mu)}\right]_+$ being a
lattice and thus having the Riesz Decomposition Property) map $T$ has a greatest
linear selection
$$
S : \left[\lclassg {1} {(\Omega, \mu)}\right]_+ \to \mathcal V
$$
defined by the formula \eqref {thesubmap} which is upper semicontinuous by Proposition~\ref {superctolinc}.
The claims of Theorem~\ref {everylsl1} then follow from Theorem~\ref {linselndl1} applied to the map $S$.

We mention without going into details
that the case of linear maps acting on $\left[\lclassg {1}\right]_+$ treated in this section is closely related to
the theory of set-valued measures and set-valued integrals
that has been fairly extensively developed in the literature; see, e.~g., \cite {cascalesetal2012} and references contained therein.
A linear set-valued map $T$ on $\left[\lclassg {1}\right]_+$ acting into a linear space~$X$ naturally corresponds to
a set-valued measure $M (A) = T (\chi_A)$,
and at least under some additional assumptions linear selections of $T$ can be recovered from selections of $M$
that can in some cases be described as measurable selections $f$ of a Radon-Nikod\'ym derivative $F$ of $M$,
which provides a parametrization of the set of suitable linear selections $S$ of $T$
of the form $S_f (g) = \int g\, df$ in the spirit of Section~\ref {mosom} and also provides additional sufficient conditions for the existence
of linear selections of~$T$.  However, this approach requires that $X$ has the Radon-Nikod\'ym property, which is notoriously absent
in the case $X = \lclassg {1}$ that usually appears in the applications, so we will not pursue it in the present work.

\section {Additivity and linearity}

\label {addandlin}

In this short section we give some remarks about rather close relationship between additivity and linearity for maps between semilinear spaces.
Similar results are discussed in \cite {smajdor1990} (see also \cite {smajdor1996}).
This might be useful in applications if a set-valued map under study is merely superadditive but we want to know if its additive selections
are superlinear.

It is easy to see that if $A : \mathcal S \to \mathcal V$ is an additive map acting between two semilinear spaces
$\mathcal S$ and $\mathcal V$ then
$A (k x) = k A (x)$ for any positive integer $k$, hence $A (\lambda x) = \lambda A (x)$ for all $x \in \mathcal S$ and positive rational
$\lambda \in \mathbb Q_+$.  This property is called \emph {$\mathbb Q_+$-homogeneity}.  Since $\mathbb Q_+$ is dense in $\mathbb R_+$
this property is often as good as the usual homogeneity.
However, it is well known that there exists an additive function $T : \mathbb R \to \mathbb R$ which is not linear
(see, e. g., \cite [\S 2.26] {cia}, disregard the slight confusion of terms).
Equipping $\mathcal V = \mathbb R$ with the trivial partial order ($x \leqslant y$ if and only if $x = y$) shows that in Theorem \ref {linselg}
the map $S = T$ might not be linear despite being $\mathbb Q_+$-linear.
However, introduction of certain rather weak continuity assumptions remedies the situation.
\begin {proposition}
\label {lscont}
Suppose that $\mathcal S$ and $\mathcal V$ are semilinear spaces and $A : \mathcal S \to \mathcal V$ is an additive map.
Suppose also that $\mathcal V$ is a monoid semilinear space equipped with a Hausdorff
topology such that the semilinear operations of $\mathcal V$ are
continuous, and
for every $x \in \mathcal S$ the map $\nu \mapsto A (\nu x)$ is bounded at some $\alpha \geqslant 0$, i. e. for every neighbourhood $U$ of $0$
in $\mathcal V$ there exist some $m > 0$ and $\varepsilon > 0$ such that $A (\nu x) \in m U$ for $|\nu - \alpha| < \varepsilon$.
Then $A$ is linear.
\end {proposition}
Suppose that $\lambda > 0$ and $x \in \mathcal S$; we need to prove that $A (\lambda x) = \lambda A (x)$.
We already know that $A$ is $\mathbb Q_+$-homogeneous.
Let
$$
a_\lambda = \{ \mu > 0 \mid \lambda + \mu \in \mathbb Q_+ \}.
$$
Thus we have
\begin {equation}
\label {ext11}
A (\lambda x) + A (\mu x) = A ((\lambda + \mu) x) = (\lambda + \mu) A (x) = \lambda A (x) + \mu A (x)
\end {equation}
for any $\mu \in a_\lambda$, and we need to pass to the limit $\mu \to 0$ in this relation.
The right hand part of \eqref {ext11} converges to $\lambda A (x)$ because $\mathcal V$ is monoid.  By the assumptions of the proposition
$\nu \mapsto A (x \nu)$ is bounded at some $\alpha \geqslant 0$.  Suppose that $U$ is an arbitrary neighbourhood of 0; we need to show that
$A (\mu x) \in U$ for all sufficiently small $\mu \in a_\lambda$.
By the assumptions there exist some $m > 0$ and $\varepsilon > 0$ satisfying $A (\nu x) \in m U$ for $|\nu - \alpha| < \varepsilon$.
For every $\mu > 0$ there exists some $\beta \in \mathbb Q_+$ so that $|\beta \mu - \alpha| < \varepsilon$,
and therefore $A (\beta \mu x) \in m U$, so $A (\mu x) \in \frac m \beta U$
by the $\mathbb Q_+$-homogeneity of $A$.  It is easy to see that we can make $\frac m \beta \leqslant 1$ for all sufficiently small $\mu$
(for $\alpha > 0$ it is always the case, and in the case $\alpha = 0$ we can take $\beta$ arbitrarily large by making $\mu$ arbitrarily small),
which concludes the proof of Proposition \ref {lscont}.

\begin {proposition}
\label {linselgc}
Suppose that under the assumptions of Theorem~\emph{\ref {linselg}}\linebreak{} $\mathcal V$ is a monoid semilinear space equipped with a Hausdorff
topology such that the semilinear operations in $\mathcal V$ are
continuous and
for every $x \in \mathcal S$ map $\nu \mapsto T (\nu x)$ is bounded at some $\alpha \geqslant 0$.  Suppose also that the order
$\leqslant$ in $\mathcal V$ is continuous in the following sense: if $a_n \leqslant b_n$ for some $a_n, b_n \in \mathcal V$ and
$b_n$ converges to 0 then $a_n$ also converges to 0.
Then any additive submap $A \leqslant T$ of $T$ is linear.  In particular, the map $S$ defined by \eqref {thesubmap} is linear.
\end {proposition}
Indeed, the proof of Proposition \ref {lscont} above shows that it is sufficient to establish that $A (\mu x) \to 0$ as
$a_\lambda \ni \mu_j \to 0$
for any $x \in \mathcal S$ and $\lambda > 0$.  Since $A \leqslant T$, by the continuity of the order we can replace $A$ by $T$.
As before, let $U$ be an arbitrary neighbourhood of 0; we need to show that
$T (\mu x) \in U$ for all sufficiently small $\mu \in a_\lambda$.
By the hypothesis there exist $m > 0$ and $\varepsilon > 0$ such that $T (\nu x) \in m U$ for $|\nu - \alpha| < \varepsilon$.
Since $a_\lambda$ is dense in $\mathbb R_+$,
for every $\delta_j > 0$ and $M_j > 0$ there exists some $\mu_j \in a_\lambda$ and an integer $n > M_j$ satisfying $|n \mu - \alpha| \leqslant \varepsilon$,
and therefore $T (n \mu x) \in m U$.
Note that $n T (\mu x) \leqslant T (n \mu x)$ by the superadditivity of $T$,
so $T (\mu x) \in \frac m n U$.  Since $n$ can be chosen arbitrarily large we can assume that $\frac m n \leqslant 1$,
and therefore $T (\mu x) \in U$ as claimed.
The proof of Proposition~\ref {linselgc} is complete.

\section {Application to inversion of non-injective operators}

\label {ationbo}

The classical Bounded Inverse Theorem (also known as Inverse Mapping Theorem)
states that any bounded linear operator $T : X \to Y$ acting between a couple of Banach spaces
has a bounded inverse $T^{-1}$ if $T$ is bounded and bijective (i. e. $T$ is \emph {onto}, that is $T (X) = Y$, and $T$ is injective,
that is $T x = 0$ implies $x = 0$).  Its applications include situations when one has a problem
in the form of a linear equation $T x = y$ with respect to $x$ that has unique solutions and question arises whether the solutions of this equation
depend continuously on
$y$.

Sometimes, however, solutions are not unique, and thus the Bounded Inverse Theorem is not directly applicable.
In this case, i. e. under the same assumptions except injectivity of $T$, Open Mapping Theorem (also known as the Banach-Schauder theorem)
implies that there still is a constant $C$
such that for any $y \in Y$ there exist some $x \in X$ satisfying $T x = y$ and $\|x\|_X \leqslant C \|y\|_Y$.
Can we find a bounded right inverse for~$T$, i. e. a bounded operator $T^{-1}$ solving $T x = y$ with respect to $y$, that is
$T T^{-1} = I_Y$ where $I_Y$ is the identity operator in $Y$?
Incidentally, in the dual terms, i. e. by passing to adjoints, this becomes a question about the existence of a bounded left inverse for $T^*$
if $T^*$ is injective.
In certain cases results of this sort are readily available and widely used.  Let us examine some of them to look at typical difficulties that may arise.
Since $T$ is coninuous, the kernel $\mathcal N (T) = \{x \in X \mid T x = 0 \}$ of $T$ is closed in $X$, so the quotient $X \slash \mathcal N (T)$
is also a Banach space, and thus we can apply Bounded Inverse Theorem to obtain a bounded inverse $T_*^{-1}$ for the quotient mapping
$T_* : X \slash \mathcal N (T) \to Y$ defined by $T_* ([x]) = T (x)$ for $[x] = \{x + \mathcal N (T)\} \in X \slash \mathcal N (T)$
(for more on quotient mappings see, e. g., \cite {dieudonnev21970}).
If $\mathcal N (T)$ is a complemented subspace of $X$, composition of the projection $P$ onto the complement of $\mathcal N (T)$ in $X$ with $T_*^{-1}$
gives a right inverse $T^{-1}$ for $T$.
In particular, if $X$ is a Hilbert space, we can take the orthogonal projection $P$, which would make
the solution of $T x = y$ given by the resulting operator $T^{-1} = P T_*^{-1}$ have minimal norm in $X$ and thus unique; this technique
is known as the least squares method.  But what if the kernel $\mathcal N (T)$ is not complemented in $X$?
If $Y = \lclassg {1}$, it is easy to see that the set-valued map arising from application of the Open Mapping Theorem is actually superlinear on
the cone of nonnegative functions $\left[\lclassg {1}\right]_+$,
so under certain conditions it is possible to apply Theorem~\ref {everyls} or Theorem~\ref {everylsl1} to this map.

\begin {theorem}
\label {rit}
Suppose that $X$ and $Y$ are Banach spaces, $Y$ has a generating cone $\mathcal S \subset Y$ satisfying the
Riesz decomposition property and the norm $\|\cdot\|_Y$ is linear on $\mathcal S$.
Let $T : X \to Y$ be a continuous linear surjection.  
Suppose also that at least one of the following conditions are satisfied.
\begin {enumerate}
\item
$X$ is a dual space, i. e. $X = Z^*$ for some Banach space $Z$.
\item
$X$ is a Banach ideal lattice of measurable functions on measurable space $\mathbb Z$
satisfying the Fatou property and
the kernel $\mathcal N (T)$ of $T$ is closed in measure.
\item
$X$ is a Banach ideal lattice of measurable functions on a $\sigma$-finite measurable space
satisfying the Fatou property, $Y = \lclassg {1}$ on a separable measurable space,
$\mathcal S = \left[ \lclassg {1} \right]_+$ and
the kernel $\mathcal N (T)$ of~$T$ is closed in measure.
\item
$Y = \lsclass {1}$ and $\mathcal S = \left[ \lsclass {1} \right]_+$.
\end {enumerate}
Then $T$ admits a bounded right inverse $T^{-1}$.
In particular, $T^{-1}$ embeds $Y$ into $X$ as a subspace, and if $Y$ is not reflexive then $X$ is not reflexive.
Also, let $z$ be a point in $X$ satisfying $T z \in \mathcal S$.
The following conditions are equivalent.
\begin {enumerate}
\item
$T$ admits a bounded right inverse $T^{-1}$ that preserves the point $z$, i. e. $T^{-1} T z = z$.
\item
For any finite decomposition $y = \sum_{j = 1}^N y_j$ of $y = T z$ into some $y_j \in \mathcal S$
there exists a decomposition $z = \sum_{j = 1}^N z_j$ of $z$ satisfying $\|z_j\|_X \leqslant c \|y_j\|_{Y}$ and $T z_j = y_j$,
with some constant $c$ that does not depend on the decomposition $\{y_j\}$.
\end {enumerate}
\end {theorem}
Indeed, let us define a set-valued map $\Phi : \mathcal S \to 2^X$ by
$$
\Phi (y) = \left\{ x \in X \mid T x = y, \|x\|_X \leqslant C \|y\|_{Y} \right\}
$$
for $y \in \mathcal S$.  It is easy to see that $\Phi$ is superlinear: positive homogeneity is trivial, and for every $f, g \in Y$ conditions
$u \in \Phi (f)$, $v \in \Phi (g)$ imply that
$\|u + v\|_X \leqslant \|u\|_X + \|v\|_X \leqslant C \left( \|f\|_Y + \|g\|_Y \right) = C \|f + g\|_Y$ by linearity of the norm $\|\cdot\|_Y$.
By the Open Mapping Theorem $\Phi$ has nonempty values if $C$ is chosen sufficiently large.
Values of $\Phi$ are closed in $X$ because $T$ is continuous.
Thus, if the first set of assumptions is satisfied then values of $\Phi$ are compact in the weak* topology $\sigma (X, Z)$,
so $\Phi$ acts into a space $\mathcal V (\mathcal W)$ defined in Proposition \ref {csgs} where we set
$\mathcal W = X$ with the weak* topology $\sigma (X, Z)$.
If the second or third set of assumptions is satisfied then values of $\Phi$ are bounded and closed in measure,
so $\Phi$ takes values in the space $\mathcal V (X)$ defined in Proposition \ref {csgsl}.
Under assumptions 1 and 2 map $\Phi$ satisfies conditions of Theorem~\ref {everyls} by either Proposition~\ref {csgs} or Proposition~\ref {csgsl},
so there exists a linear selection $\varphi$ of $\Phi$.
It is easy to see that map $\Phi$ has closed graph and thus it is upper semicontinuous by Proposition~\ref {grclos}, so
under the third set of assumptions $\Phi$ satisfies the conditions of Theorem~\ref {everylsl1}.
If the fourth set of assumptions is satisfied then we can construct a suitable linear selection in an obvious way: let $e_j$ be the coordinate 
vectors in $\lsclass {1}$, choose $z_j \in \Phi (e_j)$, and let
\begin {equation}
\label {varpph}
\varphi (\{b_j e_j\}) = \sum_j b_j z_j, \quad b_j \geqslant 0.
\end {equation}
Thus under any of the conditions 1--4 map $\Phi$ has a linear selection~$\varphi$.
Since the cone $\mathcal S$ is generating for the space $Y$, we can define a map $T^{-1} : Y \to X$
as the linear extension of $\varphi$ by Proposition \ref {gccont}.  It is easy to see that $T^{-1}$ is a right inverse of $T$
with norm at most $C$.  The rest of the statement \ref {rit}
is an interpretation of the conditions of Theorem~\ref {everyls} or Theorem~\ref {everylsl1} for the map $\Phi$.
Suppose that $T$ admits a bounded right inverse $T^{-1}$ with norm $c$.
Then for any finite decomposition $y = \sum_{j = 1}^N y_j$ of $y = T z \in \mathcal S$ into some elements of $\mathcal S$ elements
$z_j = T^{-1} y_j$ form a decomposition $z = \sum_{j = 1}^N z_j$ of $z$ satisfying
$\|z_j\|_X \leqslant \left\|T^{-1}\right\|_{Y \to X} \|y_j\|_{Y} \leqslant c \|y_j\|_{Y}$.
Conversely, if under any of assumptions 1--3
condition 2 is satisfied and we take constant $C$ in map $\Phi$ to be at least $c$, then by Theorem~\ref {everyls} or Theorem~\ref {everylsl1}
there exists a linear selection $\varphi$ of $\Phi$ satisfying $\varphi (T z) = z$, and therefore $T^{-1} T z = z$ for the map $T^{-1}$.
Assumptions 4 involve a bit of manual labor.  Let $y = \sum_{j = 1}^\infty a_j e_j$.
We apply condition 2 to a sequence of decompositions $y = \sum_{j = 1}^N a_j e_j + \sum_{j > N} a_j e_j$ to obtain $z_{N,j} \in X$ satisfying
$T z_{N,j} = a_j e_j$ for $1 \leqslant j \leqslant N$, $T z_{N,N + 1} = \sum_{j > N} a_j e_j$ and $\sum_{j = 1}^{N + 1} z_{N, j} = z$.  It is evident that we can take
$z_j = z_{j, j} = z_{N, j}$ for all $j$ and thus $y_j = T z_{j, j}$ and $\sum_j z_{j, j} = z$, so the selection \eqref {varpph} satisfies the necessary properties.
The proof of Theorem \ref {rit} is complete.

We mention that there may not be a continuous right inverse for operator $T$ if conditions of Theorem \ref {rit} are not satisfied.
\begin {comment}
For example, let $\{h_n\}$ be the Haar system in $\lclassg {1} = \lclass {1} {[0, 1]}$ normed by $\|h_n\|_{\lclassg {1}} = 1$, and define a map
$T : \lsclass {1} \to \lclassg {1}$ by $T \{a_n\} = \sum_n a_n h_n$.  Then $T$ is bounded with norm $1$.  Since the Haar system is a Schauder
basis for $\lclassg {1}$, $T$ is surjective.  However, $T$ does not admit a bounded right inverse, because $\lclass {1} {[0, 1]}$ is not embeddable
into $\lsclass {1}$ (the space $\lclass {1} {[0, 1]}$ contains an infinite-dimensional Hilbert space $\lclass {2} {[0, 1]}$,
but the space $\lsclass {1}$ does not contain any infinite-dimensional Hilbert spaces because it has the Schur property).
\end {comment}
An interesting example of this is a theorem of Peetre \cite {peetre1979} stating that the trace operator $\mathrm {Tr}$ defined on the Sobolev space
$\mathrm L_1^{(1)} \left( \mathbb R^{n + 1}_+ \right)$ acting from the upper half-space $\mathbb R^{n + 1}_+$ to $\lclass {1} {\mathbb R^n}$,
despite being onto (see, e. g., \cite [Lemma 19] {pelczynskiwojciechowski2003}),
admits no bounded right inverse; see \cite [Theorem 20] {pelczynskiwojciechowski2003}.



\section {Application to approximation of linear spaces}

\label {atass}

In this section we mention how selection theorems for superlinear maps could be applied to certain questions related to approximation.
There is a huge amount of literature dedicated to various problems of approximation of functions.
Although I was not successful in finding any concrete
useful applications for what follows in this section, the general principle seems to be interesting enough to include in this survey.

Let $X$ be a normed vector space and $Y$ its subspace.  The value
\begin {equation}
\label {appxey}
E_Y (f) = \inf_{P \in Y} \|P - f\|_X
\end {equation}
determines how well a given function $f \in X$ can be approximated by a function from $Y$.
There usually are two main questions associated with approximation: how the value $E_Y (f)$ can be estimated for various subspaces $Z$ of $X$,
and how the best (or, failing that, merely just good enough) approximating function can be obtained.  An answer to the first question usually has the form
\begin {equation}
\label {appx}
E_Y (f) \leqslant C \|f\|_Z
\end {equation}
with some constant $C$ and a seminorm $\|\cdot\|_Z$ (or something like a seminorm, e. g. a modulus of continuity).
Estimates of the form \eqref {appx} are usually called
\emph {Jackson-type inequalities}.
A classical example is the space of continuous functions $X = C (\mathbb T)$ on the unit circle $\mathbb T$ with the usual supremum norm,
the space $Y_n$ of trigonometric
polynomials on $\mathbb T$ of degree at most $n - 1$
and the space $Z = C^{(r)} (\mathbb T)$ of $r$ times continuously differentiable functions
with seminorm $\|f\|_{C^{(r)}} = \sup_{z \in \mathbb T} |f^{(r)} (z)|$.  In this case \eqref {appx} is known as the classical Jackson inequality
with $C = \frac {C_r} {n^r}$, where $C_r$ (which is known as the Akhiezer–Krein–Favard constant) depends only on $r$;
see, e. g., \cite {achieser1956}.

The nature of the approximating transformation $f \mapsto P$ used to obtain inequality \eqref {appx} is rarely explicitly mentioned in the literature.
However, on closer examination,
the approximating function $P$ is nearly always expressed as a linear transformation of the function $f$ being approximated,
even in cases where one might expect some degree of nonlinearity (for example, if weighted spaces are involved; see, e. g., \cite {lubinsky2007}).
Thus a question arises: can the approximation satisfying an estimate \eqref {appx}
be carried out by a linear transformation, i. e. does \eqref {appx} imply that there exists a linear operator $T : X \to Y$ satisfying
\begin {equation}
\label {appxT}
\|f - T f\|_X \leqslant C \|f\|_Z
\end {equation}
for all $f \in X$?
Some answers are readily available.  Similarly to what was discussed in Section \ref {ationbo},
if $X$ is a Banach space and $Y$ is a closed subspace of $X$, then $E_Y (f)$ is just the norm of $f$ in the quotient space
$X \slash Y$, and the question is, essentially, whether we can find a bounded projection from $(X \slash Y) \cap Z$ onto $X$.
If $X = Z$ is a Hilbert space and $Y$ is a closed subspace of $X$ then the
orthogonal projection from $X$ onto~$Y$ solves this problem.
One might consider a question stronger than \eqref {appxT} about the best approximation,
i.~e. if we can find an operator $T : Z \to Y$ providing the best (or ``almost the best'' in some sense)
possible approximation in $X$ by $Y$ for every function $f \in Z$.
There is a set-valued map
$\pi_Y : X \to \mathcal C (Y)$ defined by
\begin {equation}
\label {piY}
\pi_Y (f) =  \left\{Q \mid \|Q - f\|_X \leqslant (1 + \varepsilon) \inf_{P \in Y} \|P - f\|_X \right\}
\end {equation}
that gives approximating functions $P$ for any $\varepsilon > 0$.
As there usually is a function $P$ realizing the infinum in \eqref {appxey} (for example, this is always the case if $Y$ is finite-dimensional),
we can take $\varepsilon = 0$ in \eqref {piY}.
In this case $\pi_Y$ is called the \emph {metric projection}.
If $\pi_Y$
admits a linear selection then a linear mapping
$T$ satisfying \ref {appxT} also exists, but apart from the Hilbert space case this happens only rarely.
For example, in \cite [Theorem 3.5] {deutsch1983}
it was shown for the case $X = \lclassg {1}$ on a measurable space $(\Omega, \mu)$
that if $Y$ is a finite-dimensional nontrivial subspace of $X$ and $\pi_Y$ admits a linear selection then measure
$\mu$ has at least one atom.
The situation with \eqref {appxT} is much better as it is easy to see that in some cases the corresponding set-valued map is superlinear,
so we can apply results of this survey to obtain a linear selection.
For simplicity we will state the following theorem for finite-dimensional spaces $Y$ only.

\begin {theorem}
\label {appgen}
Suppose that $X$ is a normed space, $Y, Z \subset X$ are linear sets, $Z$ is equipped with a quasiseminorm with respect to which
estimate \eqref {appx} holds for all $f \in Z$, $\mathcal S \subset Z$ is a generating cone for $Z$ that satisfies the Riesz decomposition property
and the quasiseminorm $\|\cdot\|_Z$ is $1$-concave on $\mathcal S$
(i. e. $\|f\|_Z + \|g\|_Z \leqslant \|f + g\|_Z$ for any $f, g \in \mathcal S$).
Suppose also that $Y$ is finite-dimensional.
Then there exists a linear operator $T : Z \to Y$ satisfying \eqref {appxT}.
\end {theorem}

Indeed, since $Y$ is finite-dimensional, the infinum in \eqref {appx} is always attained.
Define a set-valued map $\Phi : \mathcal S \to 2^Y \setminus \{\emptyset\}$ by
$$
\Phi (f) = \left\{ P \in Y \mid \|f - P\|_X \leqslant C \|f\|_Z \right\}
$$
for $f \in \mathcal S$.
It is evident that $\Phi$ is positive homogeneous and from $1$-concavity of the norm of $Z$ on $\mathcal S$ it follows that
for any $f, g \in \mathcal S$, $P \in \Phi (f)$ and $Q \in \Phi (g)$ we have
\begin {multline*}
\|f + g - (P + Q)\|_X \leqslant \|f - P\|_X + \|g - Q\|_X \leqslant
\\
C \left(\|f\|_Z + \|g\|_Z\right) \leqslant C \|f + g\|_Z,
\end {multline*}
i. e. $P + Q \in \Phi (f + g)$, and therefore $\Phi$ is superlinear.  As $Y$ is finite-dimensional, it is easy to see that
values of $\Phi$ are compact in $Y$, so by Theorem \ref {everyls} there exists a linear selection $\varphi$ of $\Phi$.
Since $\mathcal S$ is generating for $Z$, we can extend $\varphi$ to a linear operator $T$ on $Z$ by Proposition \ref {gccont}.
Thus $T$ satisfies \eqref {appxT}, which concludes the proof of Theorem \ref {appgen}.

Let us give an example how Theorem \ref {appgen} can be applied to what is sometimes refered to as the Jackson-Favard estimate
for suitable functions $f$ on the segment
$[-1, 1]$, $1 \leqslant p \leqslant \infty$ and $n \geqslant r$ (see, e. g., \cite [Chapter 8, (7.1)] {devorelorentz1993}):
\begin {equation}
\label {jackfav}
\inf_{P \in \mathbb P_n} \|P - f\|_{\lclassg {p}} \leqslant C_r n^{-r} \left\| \varphi^r f^{(r)} \right\|_{\lclassg {p}},
\end {equation}
where $Y = \mathbb P_n$ denotes the space of all polynomials of degree at most $n$,
and $\varphi (x) = \sqrt {1 - x^2}$ for $x \in [-1, 1]$.
Note that in this case Theorem \ref {appgen} does
not give any new information (apart from the equality of optimal constants in \eqref {jackfav} and \eqref {jackfavt} below)
since a suitable linear operator $T$ is already implicit in the proof of the inequality.
In order to apply Theorem \ref {appgen}, we have to restrict ourselves to the case $p = 1$.  Then we take $X = \lclassg {1}$,
$Z$ to be the Sobolev space of measurable functions $f$ that are $r$ times differentiable almost everywhere,
$f^{(k)} (0) = 0$ for $0 \leqslant k < r$, and their norm
$\|f\|_Z = \left\| \varphi^r f^{(r)} \right\|_{\lclassg {1}}$ is finite, $\mathcal S = \{f \in Z \mid f^{(r)} \geqslant 0 \text { a. e.}\}$
is isomorphic to $\left[ \lclassg {1} \right]_+$ for a weighted space $\lclassg {1}$, so it has the Riesz decomposition property.
Then application of Theorem \ref {appgen} shows that there exists a linear operator $T$ acting into $\mathbb P_n$ that satisfies
\begin {equation}
\label {jackfavt}
\|T f - f\|_{\lclassg {1}} \leqslant C_r n^{-r} \left\| \varphi^r f^{(r)} \right\|_{\lclassg {1}}
\end {equation}
for all suitable $f$.

\section {Application to analytic functions on polydisk}

\label {aafpd}

In this section we will show how Theorem~\ref {everyls} or Theorem~\ref {blmls}
can be applied to a construction of analytic functions in polydisk $\mathbb D^n$,
$n \geqslant 2$, having prescribed nonnegative radial boundary values almost everywhere.
The result of this application detailing a linear construction with suitable control
will be used below in Section \ref {atcp} to obtain a weak version of the corona
theorem on polydisk.
For basic facts regarding the theory of analytic functions on polydisk see, e.~g., \cite {rudin1969en}.

Let $RP (\mathbb T^n) \subset M (\mathbb T^n)$ be the set of Borel measures on $\mathbb T^n$ having finite variation
such that convolution of the Poisson kernel with any $\mu \in RP (\mathbb T^n)$ is a real part of an analytic function in $\mathbb D^n$.
It is well known (see, e. g., \cite {rudin1969en}) that $\mu$ belongs to $RP (\mathbb T^n)$ if and only if the spectrum $\supp \hat \mu$
of the measure $\mu$ belongs to the set
$$
\widehat {RP_n} =
\mathbb Z^n_+ \cup (-\mathbb Z^n_+) = \left\{ j = (j_1, \ldots, j_n) \mid \text {$\sigma j \geqslant 0$ for some $\sigma \in \{-1, 1\}$}\right\}.
$$
A major difference between the cases of disk ($n = 1$) and polydisk ($n \geqslant 2$)
is this restriction on the spectrum of the real part of an analytic
function (in the sense of boundary values).  This, among other things, presents considerable challenges when one tries to construct analytic functions
in polydisk with controls on the boundary values; a well-known example of this is that there is no factorization in $\hclass {p} {\mathbb D^n}$,
$n \geqslant 2$.  However, there are some weaker results.
By \cite [Theorem 2.4.2] {rudin1969en} for any lower semicontinuous
$f \in \left[ \lclass {1} {\mathbb T^n} \right]_+$ we can find a nonnegative singular measure
$\sigma_f \in M (\mathbb T^n)$ satisfying $f - \sigma_f \in RP (\mathbb T^n)$ and $\|\sigma_f\|_{M (\mathbb T^n)} = \|f\|_{\lclass {1} {\mathbb T^n}}$.
In other words, for any such $f$ we can find an analytic function $F$ in $\mathbb D^n$ such that $\Re F (0) = 0$
and its radial boundary values $F^*$ satisfy $\Re F^* = f$ almost everywhere on $\mathbb T^n$.
We extend this result to show that $\sigma_f$ can be chosen to depend linearly on $f$, which will be very useful in Section \ref {atcp} below.
First, we need a more precise statement of \cite [Theorem 2.4.2] {rudin1969en}.

Let
$G_m$ be the subgroup 
$$
G_m = \left\{ \left(e^{i \theta_1}, \ldots, e^{i \theta_n}\right) \in \mathbb T^n \mid \exp \left( i 2^m \sum_{j = 1}^n \theta_j \right) = 1 \right\}
$$
of the group $\mathbb T^n$
with respect to pointwise multiplication for any integer $m \geqslant 0$.
It is easy to see that $G_m$ is continuous and compact with respect to the topology induced from $\mathbb T^n$.
Let $\mu_m$ be the Haar measure of the group $G_m$.
Observe that $\mu_m$ is, up to a normalization factor of order $2^{-m}$, the
$(n - 1)$-dimensional Lebesgue measure on the diagonal sections
$\left\{ e^{i\theta} \mid \sum_{j = 1}^n \theta_j = 2 \pi \frac {k} {2^m} \right\}$ of the torus $\mathbb T^n$, $k \in \mathbb Z$.
Let $\mu = \sum_{k \geqslant 0} \mu_k$; it is easy to see that $\mu$ is a $\sigma$-finite singular Borel measure on $\mathbb T^n$.
Let $G = \bigcup_{m \geqslant 0} G_m$ with Borel $\sigma$-algebra induced from $\mathbb T^n$.
For convenience, let $G_{-1} = \emptyset$.
Observe that $\chi_{G_s} \mu_m = 2^{s - m} \chi_{G_s} \mu_s$ for all $m \geqslant s$.
It is easy to see, using this scaling condition, that
$$
\chi_{G_s \setminus G_{s - 1}} \mu = \sum_{k \geqslant s} \chi_{G_s \setminus G_{s - 1}} \mu_k =
\\
\sum_{k \geqslant s} \chi_{G_s \setminus G_{s - 1}} 2^{s - k} \mu_s = 2 \chi_{G_s \setminus G_{s - 1}} \mu_s
$$
for all $s \geqslant 0$.  Therefore, again using that $\mu_m = 2^{s - m} \mu_s$ on $G_s$ for~$s \leqslant m$,
\begin {multline*}
\mu_m = \sum_{s = 0}^m \chi_{G_s \setminus G_{s - 1}} \mu_m = \sum_{s = 0}^m \chi_{G_s \setminus G_{s - 1}} 2^{s - m} \mu_s =
\\
\sum_{s = 0}^m \chi_{G_s \setminus G_{s - 1}} 2^{s - m - 1} \mu = \weightw_m d\mu
\end {multline*}
for all $m \geqslant 0$ with 
$\weightw_m = \sum_{s = 0}^m \chi_{G_s \setminus G_{s - 1}} 2^{s - m - 1}$.

\begin {proposition}
\label {rudincorrp}
Let $f$ be a lower semicontinuous function from $\lclass {1} {\mathbb T^n}$ which is positive everywhere on $\mathbb T^n$.
Then there exists a $\mu$-measurable nonnegative function $\rho$ satisfying the following conditions: $\rho \leqslant \frac 1 2 f$ $\mu$-almost everywhere,
$\rho \in \lclass {1} {d\mu}$, $f\, dm - \rho\, d\mu \in RP (\mathbb T^n)$
and
$$
\|\rho d\mu\|_{M (\mathbb T^n)} = \|\rho\|_{\lclass {1} {d\mu}} = \|f\|_{\lclass {1} {\mathbb T^n}}.
$$
\end {proposition}

The proof follows the argument of \cite [Theorem 2.4.2] {rudin1969en}.
As $f$ is positive and lower semicontinuous, it is easy to construct a decomposition $f = \sum_{m \geqslant 1} f_m$
into nonnegative trigonometric polynomials $f_m$ of degree at most $2^m - 1$, $m \in \mathbb N$.
Let $\rho = \sum_{m \geqslant 1} \weightw_m f_m$.
Thus $\rho\, d\mu = \sum_{m \geqslant 1} f_m d\mu_m = \sum_{m \geqslant 1} \sigma_m$
with $\sigma_m = f_m d\mu_m$.
Representing
$$
f_m (\zeta) = \sum_{|l| < 2^m} \hat f_m (l) \zeta^l
$$
with the standard multi-index notation, we see that
\begin {multline}
\label {spectcoinc}
\hat \sigma_m (k) =
\int_{\mathbb T^n} e^{-i \theta \cdot k} \sum_{|l|< 2^m} f_m (\theta) d\mu_m (\theta) =
\\
\int_{\mathbb T^n} e^{-i \theta \cdot k} \sum_{|l|< 2^m} \hat f_m (l) e^{i \theta \cdot l} d\mu_m (\theta) =
\\
\sum_{|l| < 2^m} \hat f_m (l) \int_{\mathbb T^n} e^{i (l - k) \cdot \theta} d\mu_m (\theta) = \hat f_m (k)
\end {multline}
for $k \in \left[\mathbb Z^n \setminus \left(\mathbb Z^n_+ \cup (-\mathbb Z^n_+)\right)\right] \cup \{ 0 \}$,
since by the orthogonality of the characters of the group $G_m$ (or a direct computation since we know what the measure $\mu_n$ is)
the value of $\int_{\mathbb T^n} e^{i (l - k) \cdot \theta} d\mu_m (\theta)$
is $1$ for $k - l = 2^m s\,\mathbf 1$ with some $s \in \mathbb Z$ where $\mathbf 1 = (1, \ldots, 1)$, and $0$ otherwise.
Furthermore,
\begin {multline*}
\| \rho \|_{\lclass {1} {\mu}} = \sum_{m \geqslant 1} \| \weightw_m f_m \|_{\lclass {1} {\mu}} =
\sum_{m \geqslant 1} \int f_m d\mu_m =
\\
\sum_m \hat \sigma_m (0) = \sum_m \hat f_m (0) = \sum_m \int f_m = \|f\|_{\lclassg {1}}
\end {multline*}
and
$$
\rho = \sum_{m \geqslant 1} \weightw_m f_m \leqslant \frac 1 2 \sum_{m \geqslant 1} f_m = \frac 1 2 f
$$
$\mu$-almost everywhere because $\weightw \leqslant \frac 1 2$ $\mu$-almost everywhere.
Finally, \eqref {spectcoinc} means that
$$
\hat f (k) - \hat \sigma (k) = \sum_m \left( \hat f_m (k) - \hat \sigma_m (k)\right) = 0
$$
for
$k \notin \mathbb Z^n_+ \cup (-\mathbb Z^n_+)$, so $f\, dm - \sigma \in RP (\mathbb T^n)$ as claimed.
The proof of Proposition~\ref {rudincorrp} is complete.

Let $C_{0+}$ be the cone of continuous functions on $\mathbb T^n$ that are positive everywhere on $\mathbb T^n$.
It is easy to see that $C_{0+}$ is a lattice and therefore it has the Riesz decomposition property.
We are now ready to prove the main result of this section.
\begin {theorem}
\label {ymapping}
For every $n \in \mathbb N$, $n \geqslant 2$, there exists a bounded with norm $1$ linear mapping
$Y : \lclass {1} {\mathbb T^n, dm} \to \lclass {1} {G, d\mu}$
satisfying
$$
f\, dm - (Y f) d\mu \in RP (\mathbb T^n)
$$
for all $f \in \lclass {1} {\mathbb T^n, dm}$.
Moreover, operator $Y$ is positive, i.~e. $f \geqslant 0$ almost everywhere implies
that $Yf \geqslant 0$ $\mu$-almost everywhere.
If $f \in C_{0+}$ then $Y f \leqslant \frac 1 2 f$ on $\mathbb T^n$.
\end {theorem}

Define a set-valued map $\Phi : C_{0+} \to 2^{\lclass {1} {\mu}}$ by
\begin {multline*}
\label {yzero2}
\Phi (h) = \left\{
\rho \in \lclass {1} {\mu}  \mid \rho \geqslant 0 \text { a. e.}, \,
h\, dm - \rho\, d\mu \in RP (\mathbb T^n), \phantom\int
\right.
\\
\left.
\rho \leqslant \frac 1 2 h \,\, \text { a. e.}, \,
\|\rho\|_{\lclass {1} {\mu}} \leqslant \|h\|_{\lclassg {1}}
\right\}
\end {multline*}
for $h \in C_{0+}$.
By Proposition \ref {rudincorrp} map $\Phi$ takes nonempty values.
Let us verify that values of $\Phi$ are compact in the weak topology of $\lclass {1} {\mu}$.
Observe that for every $h \in C_{0+}$ the set $\Phi (h)$ is uniformly integrable because it is dominated by a bounded function $\frac 1 2 h$,
which by the Dunford-Pettis theorem means that $\Phi (h)$ is relatively compact in the weak topology of $\lclass {1} {\mu}$. Thus it suffices to show that
$\Phi (h)$ is closed in this topology for every $h \in C_{0+}$.
Fix $h \in C_{0+}$ and let $\rho_\alpha \in \Phi (h)$ be a net converging weakly to
some $\rho \in \lclass {1} {\mu}$.  Since by the weak convergence for every measurable set $E \subset G$ we have
$\int_E \rho = \lim_\alpha \int_E \rho_\alpha \leqslant \int_E \frac 1 2 h$, it follows at once that $\rho \leqslant \frac 1 2 h$.
Since functions $\zeta \to e^{i n \zeta}$ are bounded and continuous, the Fourier coefficients of the measure
$h\, dm - \rho_\alpha\, d\mu$ converge to the Fourier coefficients of the measure $h\, dm - \rho_\alpha\, d\mu$.
Therefore $h\, dm - \rho\, d\mu \in RP (\mathbb T^n)$.  Finally,
$\|\rho\|_{\lclass {1} {\mu}} \leqslant \liminf_\alpha \|\rho_\alpha\|_{\lclass {1} {\mu}} \leqslant \|h\|_{\lclassg {1}}$ by the Fatou theorem.
We have verified that $\Phi (h)$ is closed and thus compact in the weak topology of $\lclass {1} {\mu}$.
It is easy to see that $\Phi (h)$ is convex for every $h \in C_{0+}$.
In short, $\Phi$ takes values in the set
$\mathcal V (\lclass {1} {\mu})$ of all nonempty compact convex sets with respect to the weak topology of $\lclass {1} {\mu}$ as
defined in Proposition \ref {csgs}.
It is easy to see that $\Phi$ is a superlinear map, so by Theorem \ref {everyls} combined with
Proposition \ref {csgs} there exists a linear selection $Y$ of $\Phi$ that satisfies the claims of Theorem \ref {ymapping} on
the set $C_{0+}$.  

Let us now extend $Y$ onto the cone $C_+ (\mathbb T^n)$ of all nonnegative continuous functions on $\mathbb T^n$.
Set $Y f = Y (f + \delta) - Y (\delta)$ for any $\delta > 0$ and all $f \in C_+ (\mathbb T^n)$.  This is an extension since for $f \in C_{0+}$
we have $Y (f + \delta) - Y (\delta) = Y (f) + Y (\delta) - Y (\delta) = Y (f)$ due to the additivity of~$Y$.
It is easy to see that this definition does not depend on
the value $\delta > 0$: if $\delta_2 > \delta_1$ then
\begin {multline*}
Y (f + \delta_2) - Y (\delta_2) = Y (f + \delta_1) + Y (\delta_2 - \delta_1) - Y (\delta_2) =
\\
Y (f + \delta_1) + Y (\delta_2 - \delta_1) - (Y (\delta_2 - \delta_1) + Y (\delta_1)) =
\\
Y (f + \delta_1) - Y (\delta_1).
\end {multline*}
Map $Y$ extended this way remains linear: positive homogeneity is conserved because
$$
Y (\lambda f) = Y (\lambda f + \lambda \delta) - Y (\lambda \delta) = \lambda \left[ Y (f + \delta) - Y \delta \right] = \lambda Y (f)
$$
for all $f \in C_+ (\mathbb T^n)$ and $\lambda > 0$, and additivity is conserved because
\begin {multline*}
Y (f + g) = Y (f + g + 2 \delta) - Y (2 \delta) =
\\
[Y (f + \delta) - Y (\delta)] + [Y (g + \delta) - Y (\delta)] = Y (f) + Y (g)
\end {multline*}
for all
$f, g \in C_+ (\mathbb T^n)$.  This extension is still bounded with norm~$1$ because
\begin {multline*}
\|Y (f) \|_{\lclass {1} {\mu}} \leqslant \inf_{\delta > 0} \|Y (f + \delta) - Y (\delta) \|_{\lclass {1} {\mu}} \leqslant
\\
\inf_{\delta > 0} ( \|Y (f + \delta)\|_{\lclass {1} {\mu}} + \|Y (\delta)\|_{\lclass {1} {\mu}}) \leqslant 
\inf_{\delta > 0} ( \|f + \delta\|_{\lclassg {1}} + \|\delta\|_{\lclassg {1}}) \leqslant 
\\
\inf_{\delta > 0} ( \|f\|_{\lclassg {1}} + 2 \|\delta\|_{\lclassg {1}}) = \|f\|_{\lclassg {1}}
\end {multline*}
for all $f \in C_+ (\mathbb T^n)$.  
Since $C_+ (\mathbb T^n)$ is a generating cone for $C (\mathbb T^n)$, which is a lattice with norm $\|\cdot\|_{\lclassg {1}}$,
we can apply Proposition \ref {gccont} to extend $Y$ to a positive operator
$Y : C (\mathbb T^n) \to \lclass {1} {\mu}$ with norm $1$.
Then, because $C (\mathbb T^n)$ is dense in $\lclassg {1}$ we can extend $Y$
to $\lclassg {1}$ by continuity with preservation of norm.
Condition 
$f\, dm - (Y f) d\mu \in RP (\mathbb T^n)$ for all $f \in C (\mathbb T^n)$ is preserved throughout the extension because
the Fourier coefficients $\widehat {\psi d\mu} (k) = \int_{\mathbb T^n} e^{-ik\zeta} \psi d\mu$, $k \in \mathbb Z^n$ are linear and continuous in
$\lclass {1} {G, d\mu}$.
The proof of Theorem~\ref {ymapping} is complete.

\section {Application to the corona problem}

\label {atcp}

In this section we will show how Theorem \ref {ymapping} can be applied to partially extend the
result obtained in \cite {kisliakovrutsky2011en} (see also \cite {thesis} for a simpler version) to the case of polydisk.
See, e. g., \cite {treilwick2005}, \cite {douglassarkar}
and references contained therein for detailed background on the corona problem from the analytic point of view.

One frequently used reformulation of the corona problem in a standard setting is stated roughly as follows.
Suppose that $G \subset \mathbb C^m$ is a domain in $\mathbb C^m$, $m \geqslant 1$, and we are given
a finite collection of bounded analytic functions $f_j \in \hclass {\infty} {G}$ satisfying
$\sup_{z \in G} \sum_j |f_j (z)| > 0$; can we find some bounded analytic functions $g_j \in \hclass {\infty} {G}$
solving the equation $\sum_j f_j g_j = 1$?  As of the time of this writing,
for $m = 1$ and $G = \mathbb D$ this question has been settled for almost five decades, but \emph {no} domains $G \subset \mathbb C^m$, $m \geqslant 2$
are known for which the corona problem has positive solution.
On the other hand, there are examples of smooth pseudoconvex (but not strictly pseudoconvex)
domains for in which the corona problem does not have a solution in general; see \cite {fornaesssibony1993} and \cite {sibony1987}.
However, there are a fairly large number of positive results for fairly regular domains, e. g. in polydisk $\mathbb D^m$ or strictly pseudoconvex domains
like the unit ball $\mathbb B_m$ of $\mathbb C^n$,
if we replace the requirement of uniform boundedness of $g_j$ with
some weaker conditions that still provide quite a lot of control over the behavior of $g_j$, for example $g_j \in \BMO$ (in the sense of boundary values),
or even $e^{\left[(|g_j|/\lambda)^{\frac 1 {2 m - 1}}\right]} \in \lclassg {1}$ for some $\lambda > 0$ (also in the sense of boundary values);
see, e. g., \cite {trent200x}.

Another related question is related to the \emph {Toeplitz Corona Theorem}, which roughly means that we replace $1$ in the right part of the
equation $\sum_j f_j g_j = 1$ with an arbitrary function $h \in X$, where $X$ is some normed space of analytic functions, and try to solve
$\sum_j f_j g_j = h$ with respect to $g_j \in X$ with a suitable estimate on the norms of $g_j$ in $X$.  This question also appears to be
more yielding
with a significant number of results for various domains and norms,
most notably $\hclassg {p}$ for $p < \infty$ and, in particular, $\hclassg {2}$, which corresponds to the original Toeplitz Corona Theorem.
The case $X = \hclass {p} {\mathbb D^n}$ was first solved in \cite {lin1994}, \cite {lin1986}, and there are a number of results
with improved estimates and for other domains in $\mathbb C^n$; see, e.~g., \cite {treilwick2005}, \cite {trentwick2009},
\cite {anderssoncarlsson1996}, \cite {anderssoncarlsson2000} and references contained
therein.
Let us introduce a latest result for the polydisk obtained in \cite {treilwick2005}.
We use the notation $\lsclassr {2} {N} = \left\{ a \in \lsclass {2} \mid \supp a \subset \{ k \mid 1 \leqslant k \leqslant N\}\right\}$
for $N \in \mathbb N$.
\begin {theorem} {\cite {treilwick2005}}
\label {treilwickt}
Let $1 \leqslant p < \infty$.
Suppose that $m > 1$, $k \in \mathbb N$, $N \in \mathbb N$ and
$
F \in \hclass {\infty} {\mathbb D^m; \lsclassr {2} {N} \to \lsclassr {2} {k}}
$
satisfies
$\delta^2 I \leqslant F F^* \leqslant I$ and $0 < \delta^2 < \frac 1 e$.
Let $f \in \hclass {p} {\mathbb D^m, \lsclassr {2} {k}}$.
Then there exists some function $g \in \hclass {p} {\mathbb D^m, \lsclassr {2} {N}}$
satisfying $F g = f$ and
$\|g\|_{\lclass {p} {\lsclassr {2} {N}}} \leqslant C \|f\|_{\hclass {\infty} {\lsclassr {2} {k}}}$
with a constant $C$ independent of $f$ and $N$.
\end {theorem}
The classical Toeplitz corona theorem corresponds to $p = 2$ and $m = 1$.
This approach seems promising because there is some hope in extracting from it some useful information for the case $p = \infty$.
It is well known that for $m = 1$ cases $p = 2$ and $p = \infty$ are equivalent with equivalence of constants;
see, e.~g., \cite [Appendix 3] {nikolsky}.
Recently in \cite {kisliakovrutsky2011en} this equivalence was extended to a wide variety of Hardy-type spaces, in particular to all
$\hclassg {p}$, $1 \leqslant p \leqslant \infty$ for $m = 1$.  Unfortunately, the method used in \cite {kisliakovrutsky2011en}
depends heavily on factorization, which makes it unsuitable for $m > 1$.
However, if we replace factorization with a surrogate provided by Theorem \ref {ymapping} in Section \ref {aafpd},
it is possible to obtain from Theorem \ref {treilwickt}
a weaker result using these methods.

\begin {theorem}
\label {coronapdweak}
Suppose that $m > 1$, $k \in \mathbb N$, $N \in \mathbb N$ and
$$
F \in \hclass {\infty} {\mathbb D^m; \lsclassr {2} {N} \to \lsclassr {2} {k}}
$$
satisfies
$\delta^2 I \leqslant F F^* \leqslant I$ and $0 < \delta^2 < \frac 1 e$.
Let $f \in \hclass {\infty} {\mathbb D^m, \lsclassr {2} {k}}$.
Then there exists some function $g$ analytic in $\mathbb D^m$ such that $F g = f$ and the radial boundary values $g^*$ of $g$
satisfy $\|g^*\|_{\lclass {\infty} {\lsclassr {2} {N}}} \leqslant C \|f\|_{\hclass {\infty} {\lsclassr {2} {k}}}$
with a constant $C$ independent of $f$ and $N$.
\end {theorem}
Note that analytic functions with radial limits bounded a. e. may still have very nasty behavior.
As it appears from the construction below,
the solution may behave like a ratio of two functions from $\hclassg {p}$ not implicitly known, a far worse control
inside the polydisk than what is typically implied the results mentioned before
since the denominator may take arbitrarily small values inside the polydisk.
Nevertheless, solutions that are ``bounded in measure'', i. e. with radial values converging in measure to a bounded function,
appear to be new and may still be useful in some applications.

We now begin the proof of Theorem \ref {coronapdweak}.
Denote by $B$ the closed unit ball of $\hclass {2} {\lsclass {2}}$ with the weak topology.
The space $B$ is metrizable, and convergence of a sequence $f_n = \left\{f_j^{(n)}\right\}$
is equivalent to pointwise convergence of functions $f_j^{(n)}$ on the polydisk $\mathbb D^m$ for all $j$.
Fix any $\delta > 1$.
Let $P_r$, $0 \leqslant r < 1$, be the Poisson integral in $\mathbb D^m$.
We denote by $P_z \nu$ and $Q_z \nu$ convolutions of a measure $\nu$
with the Poisson kernel and the conjugate Poisson kernel for a point $z \in \mathbb D^m$ respectively.
The modulus sign $|\cdot|$, as usual,
denotes the norm in the sequence space $\lsclass {2}$; in this argument all such sequences are finite.
Let
\begin {multline*}
f_g (z) = \frac 1 {C (1 + \delta)} \exp \left[ (P_z + i Q_z) \left(\log \left(\delta + {|P_r g|}\right) dm - \right.\right.
\\
\left.\left.
[Y \log \left(\delta + {|P_r g|}\right)]d\mu\right) \right]
\end {multline*}
for $g \in B$ and $z \in \mathbb D^m$,
where $C$ is the constant in Theorem \ref {treilwickt} applied to functions $F$ and $f = g$ with $p = 2$,
singular Borel measure $\mu$ on $\mathbb T^n$ is the one defined in Section \ref {aafpd} and $Y$ is a mapping
defined by Theorem~\ref {ymapping}.

Note that every function $f_g$ is analytic, has no zeroes in $\mathbb D^m$ and satisfies estimates
$|f_g| \geqslant \frac {\delta} {C (1 + \delta)}$  and
$\|f_g\|_{\hclassg {2}} \leqslant \frac 1 C$.
Now define a set-valued map $\Phi : B \to 2^B$ by
$$
\Phi (g) = \left\{ h \in B \mid F h = f_g f \right\}.
$$
By Theorem \ref {treilwickt},
for all $g \in B$ the values $\Phi (g)$ are nonempty.
It is easy to see that $\Phi (g)$ is convex and closed in $B$.
Let us verify that the graph of $\Phi$ is closed.
Indeed, let $g_n \in B$, $h_n \in \Phi (g_n)$ and suppose that for some $g, h \in B$
sequences $g_n$ and $h_n$ converge in $B$ to some $g$ and $h$ respectively.
Then $f_{g_n}$ converges to $f_g$ pointwise in $\mathbb D^m$,
$h_n$ also converges to $h$ pointwise in $\mathbb D^m$, so we can pass to a limit in
the relation $F h_n = f_{g_n} f$.
Thus $B$ and $\Phi$ satisfy the conditions of the Ky Fan--Kakutani fixed point theorem.
\begin {theoremnn} [Ky Fan--Kakutani \cite {fanky1952}]
Suppose that $K$ is a compact set in a locally convex linear topological space.
Let $\Phi$ be a mapping from $K$ to the set of nonempty compact convex subsets of $K$.
If the graph
$\gr \Phi = \{ (x, y) \in K \times K \mid y \in \Phi (x) \}$ of $\Phi$ is closed in $K \times K$
then $\Phi$ has a fixed point, i. e. $x \in \Phi (x)$ for some $x \in K$.
\end {theoremnn}
It follows that for every $r$ there exists some $g_r \in B$ satisfying $g_r \in \Phi (g_r)$, i. e.
\begin {equation}
\label {coronafp1m}
F g_r = f_{g_r} f.
\end {equation}

Since $B$ is compact and metrizable,
there exists an increasing sequence $r_n \to 1$ such that $g_{r_n}$ converges to some $g \in B$ in $B$
and therefore $g_{r_n} \to g$ pointwise in $\mathbb D^m$.
For the same reasons, passing to a subsequence, we can assume that $f_{g_{r_n}}$ converges to some
analytic function $G$
uniformly on compact sets in $\mathbb D^m$.
The sequence $a_n = \log \left( \delta + |P_{r_n} g_{r_n}|\right)$ is bounded in $\lclassg {2}$;
therefore, passing to a subsequence, we can assume that $a_n$
converges weakly in $\lclassg {2}$ to some function~$a$.
Since $a_n \geqslant \log \delta$ almost everywhere, it is easy to conclude from the weak convergence that $a \geqslant \log \delta$
almost everywhere.
Now observe that $b \mapsto (P_z + i Q_z) (b\, dm - [Y b]d\mu)$ is a continuous linear functional on $\lclassg {2}$
(here we make use of the continuity of $Y : \lclass {1} {\mathbb T^n, dm} \to \lclass {1} {G, d\mu}$ established in Theorem \ref {ymapping}
and continuity of the inclusion $\lclassg {2} \subset \lclassg {1}$).
Therefore analytic functions $(P_z + i Q_z) (a_n dm - [Y a_n] d\mu)$
converge to $(P_z + i Q_z) (a\, dm - [Y a]d\mu)$ pointwise in $\mathbb D^m$,
and so functions $f_{g_{r_n}}$ converge pointwise in $\mathbb D^m$ to
\begin {equation}
\label {geeofzee}
G (z) = \frac 1 {C (1 + \delta)} \exp \left[ (P_z + i Q_z) (a\, dm - [Y a] d\mu) \right].
\end {equation}
Surely, \eqref {coronafp1m} implies that
\begin {equation}
\label {fggfeq1}
F g = G f
\end {equation}
on $\mathbb D^m$.

Passing to a subsequence once again, we can assume that the sequence $h_{r_n} = P_{r_n} g_{r_n}$ converges weakly to some
$h \in B$.
Then for every $0 < r < 1$ we have
$$
P_r h = \lim_{n \to \infty} P_r h_{r_n} = \lim_{n \to \infty} P_r g_{r_n} + \lim_{n \to \infty} (P_r P_{r_n} - P_r) g_{r_n} = P_r g
$$
almost everywhere because $\lim_{n \to \infty} \|P_r (P_{r_n} - I)\|_{\lclassg {1} \to \mathrm C (\mathbb T^m)} \to 0$.
Therefore $h = g$ everywhere on $\mathbb D^m$ and almost everywhere on $\mathbb T^m$.

Finally, since functions $\log |h_{r_n}|$ are plurisubharmonic and converge weakly to $h$,
\begin {multline}
\label {gpdest}
|g (z)| = \lim_n |h_{r_n} (z)| \leqslant
\lim_n \exp \left( P_z [\log |h_{r_n}|]\right) \leqslant
\\
\lim_n \exp \left( P_z [\log (\delta + |h_{r_n}|)]\right) = 
\lim_n \exp \left( P_z a_n \right) = \exp (P_z a) =
\\
C (1 + \delta) |G (z)| \exp \left( P_z [(Y a) d\mu] \right)
\end {multline}
for all $z \in \mathbb D^m$ using \eqref {geeofzee}.
Since measure $(Y a) d\mu$ is singular,
passing in \eqref {gpdest} to radial boundary values yields
$$
|g^*| \leqslant C (1 + \delta) |G^*|
$$
almost everhwhere.
By the construction $G$ has no zeroes in $\mathbb D^m$, so
function $b = \frac g G$
satisfies $F b = f$ by \eqref {fggfeq1} with the claimed estimates on $b^*$.
The proof of Theorem \ref {coronapdweak} is complete.

\section {Concluding remarks}

\label {concrem}

Semilinear spaces introduced in Section \ref {introduction}
have been studied to some extent and occasionally appear in various contexts; see, e. g.,
\cite {worth1970}, \cite {huh1974}.

Axiom 3 in our definition of a semilinear topological space given in Section \ref {cfac} and used throughout the paper
is actually a trick that allows us to entirely bypass
issues of lower semicontinuity and work exclusively with upper semicontinuity where it is needed.
Under appropriate restrictions upper semicontinuity implies the usual continuity of set-valued maps
defined by a Hausdorff topology (i. e. the topology corresponding to the Hausdorff metric); see, e. g., \cite {nikodem2003} and references contained therein.

The case of finite $\mathcal D$ in Theorem \ref {linselnd} of Section \ref {lsolm} roughly corresponds to the case of maps taking values
in bounded closed convex sets of a finite-dimensional space $\mathcal W = \mathbb R^n$
that are automatically compact, raising a question whether a sutable target space
$\mathcal V (\mathcal W)$
always has compact type if it admits a finite exhaustive set $\mathcal D$ of linear functionals that are defining for and consistent with
$\mathcal W$.  Nevertheless, the point is that one way or another the compact type of $\mathcal V$ is not explicitly needed in the proof if
the set $\mathcal D$ is finite.

It might seem that the method of obtaining linear slections of linear maps described in Section \ref {lslm},
not entirely dissimilar to other commonly used techniques (see, e. g., \cite [Chapter 9] {aubinfrankowska1990}),
may have further applications to other selection problems.
However, it is easy to construct an example of a nice continuous set-valued map having discontinuous tomographical selections.
Let $A : \mathbb R \to \mathcal C (\mathbb R^2)$ be defined by $A (x) = [(0, 0), (x, 1)]$,
i.~e. $A (x)$ is the segment connecting the origin with point $(x, 1)$, and $\mathcal D = \{f_x, f_y\}$ be the coordinate functionals
$f_x ((x, y)) = x$ and $f_y ((x, y)) = y$ in this order.  Then the corresponding tomographical selection for the point $(0, 0) \in A (0)$
is $a (x) = (0, 0)$ for $x \geqslant 0$ and $a (x) = (x, 1)$ for $x < 0$, a discontinuous function.

A natural analogue for additivity or superadditivity
in the case of maps defined on convex sets is the so-called \emph {midaffinity} or \emph {midconvexity}, respectively, i. e.
affinity or convexity with $\theta = \frac 1 2$, and, by extension, with any binary fraction $\theta = \frac {k} {2^n}$, $k, n \in \mathbb N$,
$0 < k < 2^n$.  It is easy to obtain an analogue of Proposition \ref {coaff} in Section \ref {ascm} for these properties
engage the results of Section~\ref {addandlin} to obtain conditions for which midaffinity implies affinity.

\begin {comment}
In the light of Theorems~\ref {latnecc} and~\ref {riesznecc} in Section~\ref {acocs} it is natural to inquire whether the Riesz decomposition property
is necessary for the existence of linear selections of a superlinear map in a more general setting, perhaps even in an arbitrary semilinear space
$\mathcal S$.  Unfortunately, the methods employed in Section \ref {acocs} are inapplicable in this setting.
One also wonders if there is a description similar to Theorem \ref {csselect} of all affine selections in the case of $K$ being merely a Choquet simplex.
\end {comment}

The examples showing that the Riesz decomposition property is not always necessary for the existence of a linear selection
in Theorem~\ref {linselg} presented in Section~\ref {asosas} and that not every linear set-valued map has a linear selection
presented in Section~\ref {lsolm} are simple but also unsatisfactory because they use the two-point semilinear space~
$\mathcal S$.  It is not clear whether such examples exist for more regular spaces $\mathcal S$ and $\mathcal V$,
in particular for a cone $\mathcal S$ in a linear space and $\mathcal V \subset \mathcal C (\mathcal W)$ for a linear space~$\mathcal W$.

Theorem \ref {everyls} and its extensions reveal a fairly simple necessary and sufficient condition to verify whether
a superlinear set-valued map has a linear selection with graph passing through a prescribed point in $\mathcal S \times \mathcal W$,
i. e. whether a linear selection on a one-dimensional subspace can be extended to a selection on all $\mathcal S$.
Higher order conditions, i. e. when there are several points in $\mathcal S \times \mathcal W$ to be passed, or a finite-dimensional subspace
of $\mathcal S$, are less clear and require further investigation.

Proposition \ref {csgsl} in Section~\ref {scm} is in many ways unsatisfactory, as it does not tell whether the results about linear selections
described in this paper are applicable to sets closed in measure for general Banach ideal lattices and not just those on the measurable
space $\mathbb Z$.  However, results of Section~\ref {ssvm} compensate to a large degree for this deficiency.

Despite Theorem \ref {linfcb}, results of Section \ref {conebases} for semilinear spaces with cone-bases may still be useful.
If $X$ is a suitable space of functions on a suitable topological manifold $\Omega$, for example
a space of smooth functions $X = C^{(m)} (\mathbb R^n)$, $0 \leqslant m \leqslant \infty$
or a Lebesgue space $X = \lclass {p} {\mathbb R^n}$, $0 < p < \infty$, we can construct a cone $K \subset X_+$ that is dense in $X_+$ and
has a cone-basis.  This is easily done using the Bernstein polynomials
$b_{n,k} (t) = \left( \substack{ k\\n } \right) (1 - t)^{n - k} t^k$, $t \in [0, 1]$,
and a standard partition of unity.
Thus we can construct at no cost a linear selection $S$ for a superlinear map $T$ on $K$.  If continuity of $S$ could be ensured somehow,
this would yield a continuous linear selection $S$ of $T$ on all $X$.  Results of Section \ref {sotctc}, however, show that this would
generally require some additional assumptions imposed on $T$ beyond upper semicontinuity.
There is a more direct way to see this subtle point, which is that representation of a function $f$ from $C_{0+} ([0, 1])$
(i. e. continuous and bounded away from $0$) as a sum $f = \sum_{n,k} \alpha_{n,k} b_{n,k}$, $\alpha_{n,k} > 0$,
is \emph {not unique}, as in fact an infinite number of such
decompositions are obtained by first choosing $n_1 \in \mathbb N$ and a decomposition $f = f_1 + f_2$, $f_1 = \sum_{k} \alpha_{n_1,k} b_{n_1,k}$,
$f_2 \in C_{0+} ([0, 1])$ and subsequently decomposing $f_2$ into a combination of $b_{n,k}$ with $n > n_1$.
Therefore arbitrary choice of values of $T$ on $b_{nk}$ does not guarantee in general that the map $S$ obtained in this manner will be continuous.

\begin {comment}
It is instructive to apply Theorem \ref {rit} from Section \ref {ationbo} to prove the discrete version of Wiener's Tauberian Theorem
(see, e. g., \cite [Chapter VIII, \S 6.4] {katznelson2004})
stating, in one of its forms, that the operator $T_h : \lsclass {1} \to \lsclass {1}$ of convolution $T_h f = f * h$, $f \in \lsclass {1}$ with a sequence
$h \in \lsclass {1}$ has dense range if and only if the Fourier transform $\hat h (z) = \sum_{n \in \mathbb Z} z^n h (n)$, $z \in \mathbb T$
of $h$ has no zeroes in the unit circle $\mathbb T$.
It is evident from $\sum_j (T_h g)_j = \left(\sum_j h_j\right) \left(\sum_j g_j\right)$
that if the range of $T_h$ is dense in $\lsclass {1}$ then $\sum_j h_j \neq 0$ and $T_h (\lsclass {1})$ coincides with the entire
$\lsclass {1}$.
Therefore, by Theorem \ref {rit} (which is trivial in this case)
there is a bounded right inverse $T_h^{-1}$ of $T_h$.  Since $T_h$ is translation invariant, $T_h^{-1}$
also must be translation invariant, and so $T_h^{-1} = T_g$ with some $g \in \lsclass {1}$, so $\hat {\strut h} \hat {\strut g} = 1$
and $\hat h \neq 0$ almost everywhere.
The converse is established the standard way by considering a sequence of convolutions with truncated Fourier coefficients of
$1/\hat h$.
\end {comment}

\begin {comment}
Weak compactness in $\lclassg {1}$ (and the Dunford-Pettis theorem) seems to be essential in the proof of Theorem \ref {ymapping}
in Section \ref {aafpd} despite dealing with sets of bounded functions.
A natural approach would be to take the weak* topology of $\lclass {\infty} {G}$; however, it is not clear if this would work
because the Fourier coefficients $\widehat {g d\mu}$ are not defined on all $\lclass {\infty} {G}$,
and it is not immediately clear whether the set $\lclass {\infty} {G} \cap \lclass {1} {d\mu}$ is closed in this topology.
\end {comment}

Singularities of the functions obtained with the help of Theorem~\ref {ymapping} are contained in the set $G$.
One would hypothesize that for any prescribed set $A \subset \mathbb T^n$ dense in $\mathbb T^n$ a suitable measure
supported in $A$ can be constructed.
Results similar to Theorem~\ref {ymapping} are likely possible for the unit ball of $\mathbb C^n$ and general strictly pseudoconvex domains;
see, e.~g., \cite {rudin1986}.
This would also make it possible to nearly effortlessly obtain extensions of
Theorem \ref {coronapdweak} for such domains, as the corresponding analogues of Theorem \ref {treilwickt} have been available
for some time; see, e.~g.,~\cite {anderssoncarlsson1996},~\cite {anderssoncarlsson2000}.

\begin {comment}
In the ``standard'' case $\mathcal S = \left[ \lclassg {1} \right]_+$, $\mathcal W$ being a Banach space and $\mathcal V (\mathcal W)$ being
the set of all nonempty compact convex sets,
the techniques for obtaining a linear selection of a superlinear map $T : \mathcal S \to \mathcal V (\mathcal W)$
described in this paper can reportedly be replaced by a suitable approximation construction.
Namely, suppose that $E_n$ is a sequence of
conditional expectations over finite $\sigma$-subalgebras of the measure space $\Omega$ such that $\left\|E_n f - f\right\|_{\lclassg {1}} \to 0$.
Then there is a sequence of linear operators $T_n$ satisfying $T_n E_n f \in T (E_n f)$ for all $f \in \mathcal S$, and if $T$ is upper semicontinuous
then a suitable limit of $T_n$ yields a linear selection of $T$.  Note that,
as opposed to the main results of this paper, such an argument depends on continuity of $T$.
\end {comment}








\bibliographystyle {plain}

\bibliography {bmora}

\end {document}

%% file: lwen.tex
\usepackage {hyperref}

\usepackage [utf8] {inputenc}

\usepackage{mathtext}

\usepackage {version}

\usepackage {amsmath}
\usepackage {amssymb}
\usepackage {amsthm}

\usepackage {graphicx}

\usepackage {appendix}

\usepackage{mathtext}
\usepackage{amsmath}
\usepackage{amsfonts}
\usepackage{amssymb}
\usepackage{mathrsfs}
\usepackage{amsthm}

\excludeversion {svntags}

\newtheorem {theorem} {Theorem}
\newtheorem* {theoremnn} {Theorem}

\newtheorem* {theoremchoquetmeyer} {Choquet-Meyer Theorem}

\newtheorem {proposition} [theorem] {Proposition}

\newcommand {\lclass} [2] {\ensuremath{\mathrm L_{#1} \left( #2 \right) }}
\newcommand {\lsclass} [1] {\ensuremath{\mathit l^{#1} }}
\newcommand {\lsclassr} [2] {\ensuremath{\mathit l^{#1}_{#2} }}

\newcommand {\hclass} [2] {\ensuremath{\mathrm H_{#1} \left( #2 \right) }}
\newcommand {\lclassg} [1] {\ensuremath{\mathrm L_{#1}}}

\newcommand {\hclassg} [1] {\ensuremath{\mathrm H_{#1}}}

\newcommand {\sclass} [2] {\ensuremath{\mathrm S \left( #1, #2 \right) }}

\newcommand {\BMO} {\ensuremath {\mathrm {BMO}}}

\DeclareMathOperator* {\supp} {supp}
\DeclareMathOperator* {\clos} {clos}
\DeclareMathOperator* {\tint} {int}

\DeclareMathOperator* {\sect} {sect}

\DeclareMathOperator* {\co} {co}

\DeclareMathOperator* {\gr} {gr}

\excludeversion {comment}

\newcommand {\weightw} {\ensuremath {\mathit w}}
